\documentclass[12pt,leqno]{amsart}
\usepackage{amssymb}
\usepackage{amsmath,amssymb}
\renewcommand{\baselinestretch}{1.1}

\newtheorem{theorem}{Theorem}[section]
\newtheorem{corollary}{Corollary}[section]
\newtheorem{lemma}{Lemma}[section]
\newtheorem{proposition}{Proposition}[section]
\newtheorem{condition}{Condition}[section]

\begin{document}
\theoremstyle{plain}
\newtheorem{MainThm}{Theorem}
\newtheorem{thm}{Theorem}[section]
\newtheorem{clry}[thm]{Corollary}
\newtheorem{prop}[thm]{Proposition}
\newtheorem{lem}[thm]{Lemma}
\newtheorem{deft}[thm]{Definition}
\newtheorem{hyp}{Assumption}
\newtheorem*{ThmLeU}{Theorem (J.~Lee, G.~Uhlmann)}

\theoremstyle{definition}
\newtheorem{rem}[thm]{Remark}
\newtheorem*{acknow}{Acknowledgments}
\numberwithin{equation}{section}
\newcommand{\eps}{{\varphi}repsilon}
\renewcommand{\d}{\partial}
\newcommand{\re}{\mathop{\rm Re} }
\newcommand{\im}{\mathop{\rm Im}}
\newcommand{\R}{\mathbf{R}}
\newcommand{\C}{\mathbf{C}}
\newcommand{\N}{\mathbf{N}}
\newcommand{\D}{C^{\infty}_0}
\renewcommand{\O}{\mathcal{O}}
\newcommand{\dbar}{\overline{\d}}
\newcommand{\supp}{\mathop{\rm supp}}
\newcommand{\abs}[1]{\lvert #1 \rvert}
\newcommand{\csubset}{\Subset}
\newcommand{\detg}{\lvert g \rvert}
\newcommand{\dd}{\mbox{div}\thinspace}
\newcommand{\www}{\widetilde}
\newcommand{\ggggg}{\mbox{\bf g}}
\newcommand{\ep}{\varepsilon}
\newcommand{\la}{\lambda}
\newcommand{\va}{\varphi}
\newcommand{\ppp}{\partial}
\newcommand{\ooo}{\overline}
\newcommand{\wwwx}{\widetilde{x}_0}
\newcommand{\sumkj}{\sum_{k,j=1}^n}
\newcommand{\walpha}{\widetilde{\alpha}}
\newcommand{\wbeta}{\widetilde{\beta}}
\newcommand{\weight}{e^{2s\va}}
\newcommand{\fdif}{\partial_t^{\alpha}}
\newcommand{\OOO}{\Omega}
\newcommand{\LLL}{L_{\lambda,\mu}}
\newcommand{\LLLL}{L_{\widetilde{\la},\widetilde{\mu}}}
\newcommand{\ppdif}[2]{\frac{\partial^2 #1}{{\partial #2}^2}}
\newcommand{\uu}{\mathbf{u}}
\renewcommand{\v}{\mathbf{v}}
\newcommand{\y}{\mathbf{y}}
\newcommand{\ddd}{\mbox{div}\thinspace}
\newcommand{\rrr}{\mbox{rot}\thinspace}
\newcommand{\Y}{\mathbf{Y}}
\newcommand{\w}{\mathbf{w}}
\newcommand{\z}{\mathbf{z}}
\newcommand{\G}{\mathbf{G}}
\newcommand{\f}{\mathbf{f}}
\newcommand{\F}{\mathbf{F}}
\newcommand{\dddx}{\frac{d}{dx_0}}
\newcommand{\CC}{_{0}C^{\infty}(0,T)}
\newcommand{\HH}{_{0}H^{\alpha}(0,T)}
\newcommand{\llll}{L^{\infty}(\Omega\times (0,t_1))}
\renewcommand{\baselinestretch}{1.5}
\renewcommand{\div}{\mathrm{div}\,}  
\newcommand{\grad}{\mathrm{grad}\,}  
\newcommand{\rot}{\mathrm{rot}\,}  

\title
[]
{Carleman estimate for linear viscoelasticity equations and
an inverse source problem}


\author{
O.~Yu.~Imanuvilov and \, M.~Yamamoto }

\thanks{
Department of Mathematics, Colorado State
University, 101 Weber Building, Fort Collins, CO 80523-1874, U.S.A.
e-mail: {\tt oleg@math.colostate.edu}.
Partially supported by NSF grant DMS 1312900}\,
\thanks{ Department of Mathematical Sciences, The University
of Tokyo, Komaba, Meguro, Tokyo 153, Japan
e-mail: {\tt myama@ms.u-tokyo.ac.jp}.
Partially supported by Grant-in-Aid for Scientific Research (S) 15H05740 of
Japan Society for the Promotion of Science}


\date{}

\maketitle

\begin{abstract}
We consider the linear system of viscoelasticity with
the homogeneous Dirichlet boundary condition.
First we prove a Carleman estimate with boundary values of solutions of
viscoelasticity  system.
Since a solution $\uu$ under consideration is not assumed to have
compact support, in the decoupling of the Lam\'e operator by introducing
div $\uu$ and rot $\uu$, we have no boundary condition for them, so that
we have to carry out arguments by a pseudodifferential operator.
Second we apply the Carleman estimate to an inverse source problem of
determining a spatially varying factor of the external source in the linear
viscoelastitiy by extra Neumann data on the lateral
subboundary over a sufficiently long time interval
and establish the  stability estimate.
\end{abstract}

\section{Introduction and main results}
Let $T$ be a positive constant, $x' = (x_1,..., x_n) \in \Bbb R^n$,
and $\Omega$ be a bounded domain in $\Bbb R^n$ with $\partial\Omega\in
C^\infty$, let $\vec\nu =\vec\nu(x')$ be the unit outward normal vector at $x'$ to
$\ppp\Omega$.
Let $x' \in \Bbb R^n$ be the spatial variable,
$x_0$ be the time variable, and we set
$
x = (x_0, x') = (x_0, x_1, ..., x_n).
$
Here and henceforth $\cdot^T$ denotes the transposes of vectors and matrices
under consideration, and $D= (D_0, D')$, $D_0=\frac{1}{i}
\partial_{x_0}$, $D'=(\frac{1}{i}\partial_{x_1},\dots,
\frac{1}{{i}}\partial_{x_n})$, $i=\root\of{-1}$, $\nabla' = (\ppp_{x_1}, \dots,
\ppp_{x_n})$, $\nabla = (\ppp_{x_0}, \nabla')$.

In the cylinder domain $Q:=(-T,T)\times \Omega$, we consider
the system of viscoelasticity
\begin{equation}\label{gora1A}
P(x,D)\mbox{\bf u} \equiv \rho\partial_{x_0}^2\mbox{\bf u}
- L_{\lambda,\mu}(x,D')\mbox{\bf u}
+\int_0^{x_0}L_{\widetilde \lambda, \widetilde\mu}
(x, \widetilde x_0, D')\mbox{\bf u}(\widetilde x_0, x')d\widetilde x_0
= \mbox{\bf F},
\end{equation}
\begin{equation}\label{gora2A} \mbox{\bf u}\vert_{\Sigma \triangleq (-T,T)\times \partial\Omega}=0,\quad
\mbox{\bf u}(T,\cdot) = \partial_{x_0}\mbox{\bf u}(T,\cdot)
= \mbox{\bf u}(-T,\cdot) = \partial_{x_0}\mbox{\bf u}(-T,\cdot)=0,
\end{equation}
where $\mbox{\bf u}(x) =(u_1(x),\dots , u_n(x))^T$ is the displacement  and $ \mbox{\bf F}(x)=(F_1(x),\dots,
F_n(x))^T$
is an external force. For the functions $\lambda(x)$ and $\mu(x)$,
the partial differential operator $L_{\lambda,\mu}(x,D')$ is defined by
\begin{eqnarray}\nonumber
L_{\lambda,\mu}(x,D')\mbox{\bf u} = \mu(x)\Delta\mbox{\bf u}
+ (\mu(x)+ \lambda(x))\nabla'\mbox{div}\,\mbox{\bf u}\nonumber\\
+ (\mbox{div}\,\mbox{\bf u})\nabla'\lambda + (\nabla'\mbox{\bf u}
+ (\nabla'\mbox{\bf u})^T) \nabla'\mu, \qquad x \in Q,
\end{eqnarray}
while $L_{\www\lambda,\www\mu}(x,\www{x}_0,D')$ is defined by
\begin{eqnarray}
L_{\www\lambda,\www\mu}(x, \wwwx,D')\mbox{\bf u}
= \www\mu(x,\wwwx)\Delta\mbox{\bf u}
+ (\www\mu(x,\wwwx)+ \www\lambda(x,\wwwx))\nabla'\mbox{div}\,\mbox{\bf u}
                                    \nonumber\\
+ (\mbox{div}\,\mbox{\bf u})\nabla'\lambda + (\nabla'\mbox{\bf u}
+ (\nabla'\mbox{\bf u})^T) \nabla'\mu, \qquad (x,\tilde x_0) \in Q\times (-T,T),  \nonumber
\end{eqnarray}

The coefficients  $\rho,$ $\lambda$, $\mu$, $\www\la$, $\www\mu$
are assumed to satisfy
\begin{equation}\label{pinok0}
\rho, \lambda, \mu \in C^2(\overline{Q}), \quad
 \rho(x)>0,\, \mu(x) > 0, \thinspace \lambda(x) + \mu(x) > 0 \quad \mbox{for}\,
x \in \overline{Q}
\end{equation}
and
\begin{equation}\label{pinok1}
\www\lambda, \www\mu\in C^2(\bar Q\times [-T,T]).
\end{equation}

The equation (\ref{gora1A}) is a model equation for the viscoelasticity.  The
viscoelasticity indicates a mixed physical property of
the viscosity and the elasticity.  The research for the
viscoelasticity  dates back to Maxwell, Boltzmann, Kelvin. There are some important applications of the viscoelasticity materials
related to modern technology.
For example for medical diagnosis, one has to take into consideration that
much of human tissues are viscoelastic (see works of  de Buhan \cite{dB},
Sinkus, Tanter, Xydeas, Catheline, Bercoff and Fink \cite{S} and
a monograph Lakes \cite{L}, and the references therein).
Finally we mention that
in the theory of  viscoelasticity some other equations was proposed (see the
monograph of Renardi, Hrusa and Nohel \cite{RHN}), and the
equation (\ref{gora1A}) is one of them and called the equation of linear
viscoelasticity.

The purpose of this paper is to establish the following for
(\ref{gora1A}):
\begin{enumerate}
\item
a Carleman estimate for functions without compact supports;
\item
the Lipschitz stability in an inverse source problem of determining
spatially varying factor of the external force ${\bf F}.$
\end{enumerate}

A Carleman estimate is an $L^2$-weighted estimate of solution $\uu$ to
(\ref{gora1A}) which holds uniformly in large parameter.
Carleman estimates have been well
studied for single equations (e.g., H\"ormander \cite{H}, Isakov \cite{Is}).
A Carleman estimate yields several important results such as
the unique continuation, the energy estimate called an observability
inequality and the stability in inverse problems.
However for systems whose principal part is coupled, for example even for
isotropic Lam\'e system (that is, (\ref{gora1A}) with $\www\lambda
= \www\mu \equiv 0$), the Carleman estimate is difficult to be proved for
functions $\mathbf{u}$ whose supports are not compact in $Q$.
In particular, for (\ref{gora1A}), no such Carleman estimates are not known.

In establishing a Carleman estimate for our system (\ref{gora1A}) for
non-compactly supported $\uu$, we can emphasize
the two main difficulties:
\begin{itemize}
\item
The principal part $\rho\ppp_{x_0}^2 - L_{\lambda,\mu}(x,D')$ is
strongly coupled.
\item
The Lam\'e operator $L_{\www\la,\www\mu}(x,\wwwx,D')$ appears as an integral.
\end{itemize}

In particular, because of the first difficulty, in all the existing papers
de Buhan \cite{dB}, de Buhan and Osses \cite{dBO},
Lorenzi, Messina and Romanov \cite{LMR}, Lorenzi and Romanov \cite{LR},
Romanov and Yamamoto \cite{RY}, it is assumed that functions under
consideration have compact supports or some special conditions are
satisfied in proving Carleman estimates, so that observation data
for the inverse problems have to be taken over the whole lateral
boundary $(0,T)\times \ppp\Omega$.

As for Carleman estimates for functions without compact supports and
applications to inverse problems for the Lam\'e system without the integral
terms, we refer to Imanuvilov and Yamamoto \cite{IY4} - \cite{IY7}.
For the Carleman estimate for with Lam\'e system with
$L_{\www\la,\www\mu}(x,\wwwx,D') = 0$, see Bellassoued, Imanuvilov and
Yamamoto \cite{BIY}, Bellassoued and Yamamoto \cite{BY},
Imanuvilov, Isakov and Yamamoto \cite{IIY}.
In this paper, we modify the arguments in those papers
and establish a Carleman estimate for (\ref{gora1A}) for $\mathbf{u}$ not
having compact supports.  Then we apply the Carleman estimate for an inverse
source problem by modifying the method in Imanuvilov and Yamamoto \cite{IY1} -
\cite{IY3} and Beilina, Cristofol, Li and M. Yamamoto \cite{BCLY}
discussing scalar hyperbolic equations.
As for the methodology for applying Carleman estimates to inverse
problems, we refer to a pioneering paper Bukhgeim and Klibanov
\cite{BK}.

\vspace{0.2cm}

For the statement of the Carleman estimate for (\ref{gora1A}), we need to
introduce notations and definitions.

We define the Poisson bracket by the formula
$$
\{\varphi,\psi\}=\sum^n_{j=0}\ppp_{\xi_j}\varphi\ppp_{x_j}\psi
- \ppp_{\xi_j}\psi\ppp_{x_j}\varphi.
$$
By $\overline{z}$ we denote the complex conjugate of $z \in \Bbb C$, and
we set $<a,b> = \sum^n_{k=0} a_k\overline b_k$ for
$a = (a_0, \dots, a_n), b = (b_0, \dots , b_n) \in \Bbb C^n$,
$\xi=(\xi_0,\dots,\xi_n), \xi'=(\xi_1,\dots, \xi_n), \widetilde \xi=(\xi_0,
\dots,
\xi_{n-1}), \www\nabla=(\partial_{y_0},\dots \partial_{y_{n-1}})$,
$\zeta=(\xi_0,\dots, \xi_{n-1},\tilde s),
\widetilde D=(D_0,\dots, D_{n-1})$.
\par
For $\beta \in C^2(\ooo Q)$, we introduce the symbol:
$$
p_\beta(x,\xi)=\rho(x)\xi_0^2 - \beta(x)\vert \xi'\vert^2.
$$

Let $\Gamma_0$ be some relatively open subset on $\partial\Omega$.
We set $\www\Gamma=\partial\Omega\setminus\Gamma_0$ and
$\Sigma_0=(-T,T)\times \Gamma_0,$ $\www\Sigma=(-T,T)\times \www\Gamma.$

In order to prove the Carleman estimate for the viscoelastic Lam\'e  system,
we assume the existence of the real-valued function $\psi$
which is pseudoconvex
with respect to the symbols $p_{\mu}(x,\xi)$ and $p_{\la+2\mu}(x,\xi)$.
More precisely,

\begin{condition}\label{A1}
{\it
There exists a function $\psi\in C^3(\overline Q)$ such that
$\nabla \psi(x)\ne 0$ on $\ooo{Q}$,
\begin{equation}
\{\text{p}_\beta(x,\xi-{i} s\nabla\psi(x)),
\text{p}_\beta(x,\xi+{i} s
\nabla\psi(x))\}/2{i} s>0,
      \quad\forall \beta\in\{\mu,\lambda+2\mu\}
\end{equation}
if $(x, \xi,s) \in \overline {Q}\times
(\Bbb R^{n+1} \setminus\{0\}) \times (\Bbb R^1_+\setminus \{ 0\})$
satisfies
$
\mbox{p}_\beta(x,\xi+ {i} s\nabla\psi(x))
=0.
$
On the lateral boundary, we assume
\begin{equation}\label{giga}
p_\mu(x,\nabla\psi)\vert
_{\overline{\Sigma_0}}<0, \quad \quad\mbox{and}\quad
\partial_{\vec \nu} \psi\vert
_{\overline{\Sigma_0}}<0.
\end{equation}
}
\end{condition}

Moreover we assume that
\begin{equation}\label{pinok22}
\partial_{x_0}\psi(x)<0 \quad \mbox{on}\,\,(0,T]\times \bar \Omega
\quad \mbox{and}\quad \partial_{x_0}\psi(x)>0 \quad \mbox{on}\,\,[-T,0)\times
\www\Omega.
\end{equation}

If the function $\psi$ satisfies Condition \ref{A1}, then for any constant $C$
the function $\psi+C$ also satisfies  Condition \ref{A1}.
Hence, without loss of generality, we can assume that $\psi$ is strictly positive on $\bar Q.$
Moreover, let us assume
\begin{equation}\label{giga1}
\nabla' \psi(x)\ne  0\quad \mbox{on}
\,\,\bar Q,\quad\psi(x)>0\quad\mbox{on}\,\,\bar Q.
\end{equation}
Using the function $\psi$, we introduce the function $\varphi(x)$ by
\begin{equation}\label{giga3}
\varphi(x)=e^{\tau\psi(x)},\quad \tau >1,
\end{equation}
where the parameter $\tau>0$ will be fixed below.

For any function $\mbox{\bf f}=(f_1,\dots ,f_n)$ we introduce the differential
form $\omega_{\mbox{\bf f}}=\sum_{j=1}^n f_jdx_j.$ Then $d\omega_{\mbox{\bf f}}
= \sum_{i<j}^n(\frac{\partial f_i}{\partial x_j}
- \frac{\partial f_j}{\partial x_i}) dx_i\wedge dx_j.$

We identify the differential form $d\omega_{\mbox{\bf f}}$ with the
vector-function:
$$
d\omega_{\mbox{\bf f}}=\left (\frac{\partial f_1}{\partial x_2}
-\frac{\partial f_2}{\partial x_1},\dots , \frac{\partial f_1}{\partial x_n}
-\frac{\partial f_n}{\partial x_1},\frac{\partial f_2}{\partial x_3}
-\frac{\partial f_3}{\partial x_2},\dots , \frac{\partial f_2}{\partial x_n}
-\frac{\partial f_n}{\partial x_2}, \dots,\frac{\partial f_{n-1}}{\partial x_n}
-\frac{\partial f_n}{\partial x_{n-1}}\right).
$$
Denote
$$
\Vert \mbox{\bf u}\Vert^2_{\mathcal B(\psi,s,\tau,Q)}=
\int_Q\Biggl(\sum^2_{\vert\alpha\vert=0}(s\varphi)^{4-2\vert\alpha\vert}
\tau^{6-2\vert\alpha\vert}\vert
\partial^{\alpha}_x\text{\bf u}\vert^2
+ s\varphi\tau^2\vert\nabla d\omega_{\mbox{\bf u}}\vert^2
$$
$$+(s\varphi)^3\tau^4\vert d\omega_{\mbox{\bf u}}\vert^2
+ s\varphi\tau^2\vert\nabla \mbox{div}\,\text{\bf u}\vert^2
+ (s\varphi)^3\tau^4\vert\mbox{div}\,\text{\bf u}\vert^2\Biggr)e^{2s\varphi}
dx,
$$
$$
\Vert \mbox{\bf u}\Vert^2_{\mathcal X(\psi,s,\tau,\Omega)}
= \int_\Omega\Biggl(\sum^2_{\vert\alpha'\vert=0}(s\varphi)
^{4-2\vert\alpha\vert}\tau^{5-2\vert\alpha'\vert}\vert
\partial^{\alpha'}_{x'}\text{\bf u}\vert^2
+ s\varphi\tau^2\vert\nabla' d\omega_{\mbox{\bf u}}\vert^2
$$
$$+(s\varphi)^3\tau^4\vert d\omega_{\mbox{\bf u}} \vert^2
+ s\varphi\tau^2\vert\nabla' \mbox{div}\,\text{\bf u}\vert^2
+ (s\varphi)^3\tau^4\vert\mbox{div}\,\text{\bf u}\vert^2\Biggr)e^{2s\varphi}
dx',
$$
where ${\Bbb N} := \{1,2,3,...\}$,
$\alpha=(\alpha_0,\alpha')=(\alpha_0,\alpha_1,\alpha_2,\dots,\alpha_n)$,
$\alpha'=(\alpha_1,\alpha_2,\dots,\alpha_n)$, $\alpha_j\in
{\Bbb N}\cup \{0\}$, $\partial^\alpha_x = \partial^{\alpha_0}_{x_0}
\partial^{\alpha_1}_{x_1}\partial^{\alpha_2}_{x_2}
\dots\partial^{\alpha_n}_{x_n}.$

Finally we introduce the norm
$$\Vert \text{\bf F}
e^{s\varphi}\Vert^2_{\mathcal Y(\psi,s,\tau,Q)}=\Vert\mbox{div}\, \text{\bf F}
e^{s\varphi}\Vert^2_{L^2(Q)}
+ \Vert d\omega_{\mbox{\bf F}}e^{s\varphi}\Vert^2_{L^2(Q)}+s\tau^2\Vert
\varphi^\frac 12\text{\bf F}
e^{s\varphi}\Vert^2_{L^2(\Sigma)}+\Vert\text{\bf F}
e^{s\varphi}\Vert^2_{L^2(Q)}.
$$

Now we are ready to state our first main result, a Carleman estimate
as follows.

\begin{theorem}\label{opa3}
Let $\mbox{\bf F}, \mbox{div}\, \F, d\omega_{\F}\in L^2(Q)$, $\mbox{\bf u}\in  H^1(Q)
\cap L^2(0,T; H^2(\Omega))$ satisfy (\ref{gora1A}), (\ref{gora2A}). Moreover let (\ref{pinok0})
- (\ref{giga1}) be satisfied
and let the function $\varphi$ be determined by (\ref{giga3}).
Then there exist $\tau_0>0$  and $s_0>0$ such
that for any $\tau > \tau_0$ and any  $s > s_0$  the following
estimate holds true:
\begin{eqnarray}\label{2.9'}
\Vert \mbox{\bf u}\Vert^2_{\mathcal B(\psi,s,\tau,Q)}
+ s\tau\Vert\varphi^\frac 12 \nabla\partial_{\vec \nu} \text{\bf u}
e^{s\varphi}\Vert_{{ L}^2(\Sigma)}^2 + s^3\tau^3\Vert\varphi^\frac 32 \partial_{\vec \nu} \text{\bf u}
e^{s\varphi}\Vert_{{ L}^2(\Sigma)}^2\\
\le C_1(\Vert \text{\bf F}
e^{s\varphi}\Vert^2_{\mathcal Y(\psi,s,\tau,Q)}+ s^3\tau^3\Vert\varphi^\frac 32 \partial_{\vec \nu} \text{\bf u}
e^{s\varphi}\Vert_{{ L}^2(\widetilde\Sigma)}^2+ s\tau\Vert\varphi^\frac 12
\nabla\partial_{\vec \nu} \text{\bf u}
e^{s\varphi}\Vert_{{ L}^2(\widetilde\Sigma)}^2\nonumber),
\end{eqnarray}
where the constant $C_1>0$ is independent of $s$ and $\tau.$
\end{theorem}


The proof of Theorem \ref{opa3} is given in Sections \ref{QQ1}-\ref{Q2}.
\\
\vspace{0.2cm}

Next we apply  Carleman estimate (\ref{2.9'}) to an inverse source problem of
determining a spatially varying factor of source term of the form
$\mathbf{F}(x):= R(x)\f(x')$.
Now we assume that $\rho, \la, \mu$ are independent of $x_0$:
$\rho(x) = \rho(x'), \la(x) = \la(x'), \mu(x) = \mu(x')$ for $x\in (0,T)\times \Omega$.
Let $\OOO \subset \Bbb R^n$ be a smooth bounded domain.
We consider
\begin{equation}\label{(1.1)}
\rho(x')\ppp_{x_0}^2\mbox{\bf y}
= L_{\lambda,\mu}(x,D') \mbox{\bf y} - \int^{x_0}_0 \LLLL(x,\www x_0,D')
\mbox{\bf y}(x,\www x_0) d\widetilde x_0 + R(x)\f(x') \quad
\mbox{in $(0,T)\times \OOO$},
\end{equation}
\begin{equation}\label{(1.222)}
\mbox{\bf y}(0,\cdot) = \ppp_{x_0}\mbox{\bf y}(0,\cdot) = 0, \quad \mbox{\bf y}\vert_{(0,T)\times \ppp\OOO} = 0.
\end{equation}
Here $R(x)$ is an $n\times n$ matrix function and $\f(x')$ is an
$\Bbb R^n$-valued function.

We further assume
\begin{equation}\label{(1.16)}
\www{\lambda}, \www{\mu}
\in C^2([0,T]\times \ooo{\OOO}\times [0,T]).
\end{equation}
\\

We consider\\
{\bf Inverse source problem.} {\it
Let the function $R$ be given and $\Gamma \subset \ppp\OOO.$
Determine $\f(x')$ by $\ppp_{\vec\nu}\mbox{\bf y}\vert
_{(0,T)\times \Gamma}$.}
\\

We state our main result on the inverse source problem.
\\
\begin{theorem}\label{theorem 1.2}{\it
We assume that there exists a function $\psi$  which satisfies Condition
\ref{A1}, (\ref{giga}) - (\ref{giga1}),
\begin{equation}\label{legion}
 \{x'; \thinspace
x\in [-T,T]\times \partial\Omega,\quad \partial_{\vec \nu}\psi(x)\ge 0\}\subset \Gamma
\end{equation} and
\begin{equation}\label{logoped}
\inf_{x'\in \Omega} \psi(0,x')> \max\{\inf_{x'\in \Omega} \psi(T,x'), \inf_{x'\in \Omega} \psi(-T,x')  \}.
\end{equation}
The Lam\'e coefficients $\rho,\lambda,\mu$ satisfy (\ref{pinok0}) and $\tilde\mu,\tilde\lambda$ satisfy (\ref{(1.16)}).
Let
\begin{equation}\label{(1.18)}
\vert \mbox{det} \thinspace R(0,x') \vert \ge \delta_0 > 0,
 \quad x' \in \ooo{\OOO},\quad R\in C^{2,1}([0,T]\times\bar\Omega), \quad \partial_{x_0}R(0,x')=0
                                               \end{equation}
with some constant $\delta_0 > 0$,
and $\mathbf{y}, \nabla '\mathbf{y}, \partial_{x_0}\nabla '\mathbf{y} \in H^1((0,T)\times \OOO)$ satisfy
(\ref{(1.1)}) and (\ref{(1.222)}), and $\partial^2_{x_0}\mathbf{y}
\in C([0,T_1];H^{\frac{3}{2}+\epsilon_0}(\Omega))$
with some positive $\epsilon_0$ and $T_1$.

Then there exists a constant $C_2>0$ such that
$$
\Vert \f\Vert_{H^1(\OOO)}
\le C_2\left( \sum_{j=0}^1 \Vert \nabla \ppp_{x_0}^j\ppp_{\vec\nu}
\mbox{\bf y}\Vert_{L^2((0,T)\times \Gamma)}
+ \Vert \ppp_{\vec\nu}\ppp_{x_0}^2\mbox{\bf y}(0,\cdot)\Vert
_{L^2(\Gamma)} \right)
$$
for all $\f \in H^1_0(\OOO)$.
}
\end{theorem}

Our viscoelasticity equation has a finite propagation speed and so for
determination of $\f$ over the whole domain $\OOO$, the observation time
$T$ must be longer than some critical value $T_0$

The proof of Theorem \ref{theorem 1.2} is provided in Section \ref{Q1}.

There are other works on inverse problems related to the viscoelasticity
and we refer to Cavaterra, Lorenzi and Yamamoto \cite{CLY},
Grasselli \cite{G}, Janno \cite{J}, Janno and von Wolfersdorf \cite{JW},
Loreti, Sforza and Yamamoto \cite{LSY}, von Wolfersdorf \cite{W}.
\section{Proof of Theorem \ref{opa3}}\label{QQ1}
Sections \ref{QQ1}-\ref{Q2} are devoted to the proof of Theorem \ref{opa3} and in Section \ref{Q!8}
we collect necessary lemmata for the proof.


For the function $\beta$ we introduce the operator                              $$
\square_{\rho, \beta}(x,D)=\rho(x)\partial^2_{x_0}- \beta(x)
\Delta.
$$
It is well known that the functions $d\omega_{\mbox{\bf u}}, \text{div}\,
\mbox{\bf u}$ satisfy the equations
\begin{eqnarray}\label{3.3}
\square_{\rho,\mu}(x,D) d\omega_{\mbox{\bf u}} -\int_0^{x_0}\widetilde
\mu(x,\widetilde x_0) \Delta d\omega_{\mbox{\bf u}}\, d\widetilde x_0
= \mbox{\bf q}_1 \quad \mbox{in}\,Q,\nonumber\\
\square_{\rho,\lambda+2\mu}(x,D)\text{div}\, \mbox{\bf u}-\int_0^{x_0}
(\widetilde \lambda+2\widetilde \mu)(x,\widetilde x_0) \Delta \text{div}\,
\mbox{\bf u}\, d\widetilde x_0= q_2\quad \mbox{in}\,Q,
\end{eqnarray}
\begin{eqnarray}\label{3.4}
\mbox{\bf q}_1=K_1 (x,D)d\omega_{\mbox{\bf u}} + K_2 (x,D) \text{div}\,
\mbox{\bf u}+\int_0^{x_0}(\widetilde K_1 (x,\widetilde x_0,D)d\omega
_{\mbox{\bf u}} +\widetilde K_2 (x,\widetilde x_0,D) \text{div}\,
\mbox{\bf u})d\widetilde x_0+\rho d\omega_{\mbox{\bf F}/\rho},\nonumber\\
\quad q_2=K_3 (x,D) d\omega_{\mbox{\bf u}} + K_4 (x,D) \text{div}\,
\mbox{\bf u}+\int_0^{x_0}(\widetilde K_3 (x,\widetilde x_0,D)
d\omega_{\mbox{\bf u}}+\widetilde K_4 (x,\widetilde x_0,D) \text{div} \,
\mbox{\bf u})d\widetilde x_0+\rho\text{div} \,(\mbox{\bf F}/\rho),    \nonumber
\end{eqnarray}
where $K_j (x,D),K_j (x,\widetilde x_0,D)$ are first-order differential
operators with
$C^1$ coefficients.
Now we introduce a new unknown function $\mbox{\bf v}=(\mbox{\bf v}_1,v_2)$
by formulae
\begin{equation}\label{victory}
\mbox{\bf v}_1=d\omega_{\mbox{\bf u}}- \int_0^{x_0}\frac{\widetilde
\mu(x,\widetilde x_0)}{\mu(x)} d\omega_{\mbox{\bf u}} \,d\widetilde x_0,
\quad v_2=\text{div}\, \mbox{\bf u}-\int_0^{x_0}\frac{(\widetilde \lambda
+2\widetilde \mu)(x,\widetilde x_0)}{(\lambda+2\mu)(x)} \text{div}\,
\mbox{\bf u}\, d\widetilde x_0.
\end{equation}
More specifically
$$\mbox{\bf v}_1=(v_{1,2}, \dots  v_{n-1,n}),\quad  v_{k,j}=\frac{\partial u_k}{\partial x_j}-\frac{\partial u_j}{\partial x_k}-\int_0^{x_0}\frac{\widetilde
\mu(x,\widetilde x_0)}{\mu(x)}(\frac{\partial u_k}{\partial x_j}-\frac{\partial u_j}{\partial x_k})(\tilde x_0,x')d\tilde x_0.
$$
Then
\begin{equation}\label{poko}
\mbox{\bf P}(x,D) \mbox{\bf v}=(
\square_{\rho,\mu} (x,D)\mbox{\bf v}_1, \square_{\rho,\lambda+2\mu} (x,D) v_2)=\mbox{\bf q} \quad \mbox{in}\,\, Q,
\end{equation}

where $\mbox{\bf q}=(\mbox{\bf q}_3,q_4)=(\mbox{\bf q}_1
+\widetilde{\mbox{\bf q}}_1,q_2
+\widetilde q_2):$
$$
\,\, \widetilde{\mbox{\bf q}}_1=-\int_0^{x_0}\left (2(\nabla'\frac{\widetilde
\mu(x,\widetilde x_0)}{\mu(x)},\nabla' d\omega_{\mbox{\bf u}})+ \Delta
\frac{\widetilde \mu(x,\widetilde x_0)}{\mu(x)}d\omega_{\mbox{\bf u}}\right)
\,d\widetilde x_0-\rho\partial^2_{x_0}\int_0^{x_0}\frac{\widetilde \mu(x,
\widetilde x_0)}{\mu(x)} d\omega_{\mbox{\bf u}} \,d\widetilde x_0,
$$
$$
 \widetilde q_2=-\int_0^{x_0}\left (2(\nabla'\frac{(\widetilde\lambda
+2\widetilde \mu)(x,\widetilde x_0)}{\mu(x)},\nabla ' \text{div} \,
\mbox{\bf u})+ \Delta\frac{\widetilde \mu(x,\widetilde x_0)}{\mu(x)}
\text{div} \, \mbox{\bf u}\right)\,d\widetilde x_0
$$
\begin{equation}\label{nokia}-\rho\partial^2_{x_0}\int_0^{x_0}\frac{(\widetilde\lambda+2\widetilde \mu)
(x,\widetilde x_0)}{(\lambda+2\mu)(x)}  \text{div} \, \mbox{\bf u}
\,d\widetilde x_0.
\end{equation}
We have
\begin{eqnarray}\label{nina}
\Vert \mbox{\bf q}e^{s\varphi}\Vert_{L^2(Q)}\le C_{1}( \Vert\nabla  d\omega_{\mbox{\bf u}} e^{s\varphi}\Vert
_{L^2(Q)}
+\Vert\nabla  \mbox{div}\, \mbox{\bf u}e^{s\varphi}\Vert_{L^2(Q)}+\Vert d\omega
_{\mbox{\bf u}}e^{s\varphi}\Vert_{L^2(Q)}\nonumber\\+\Vert  \mbox{div}\,\mbox{\bf u}e^{s\varphi}\Vert_{L^2(Q)}+\Vert d\omega_{\mbox{\bf F}}e^{s\varphi}
\Vert_{L^2(Q)}+\Vert  \mbox{div}\, \mbox{\bf F}e^{s\varphi}\Vert_{L^2(Q)}+\Vert \mbox{\bf F}e^{s\varphi}\Vert_{L^2(Q)}).
\end{eqnarray}

As for the  boundary conditions we have
\begin{equation}\label{poko1}
\mathcal B(x,D')\mbox{\bf v}=\mbox{\bf g}\quad \mbox{on}\,\Sigma,
\end{equation}
where
$$
\mbox{\bf g}=(g_1,\dots, g_{\frac{n^2+n}{2}}), \quad
\mathcal B(x,D')= (\widetilde {\mbox{B}}(x,D'), \mathcal C(x')),
$$
$$
\widetilde
{\mbox{B}}(x,D')=(b_{1}(x,D'),\dots ,b_{n}(x,D'))
$$
and
$$
b_{\hat j}(x,D')\mbox{\bf v}= -\sum_{j=1, j\ne \hat j}^n\mbox{sign}(\hat j-j)
\partial_{x_j} v_{j,\hat j}-\frac{(\lambda+2\mu)}{\mu}\partial_{x_{\hat j}}v_2
$$
for $1\le\hat j\le n$, and $\mathcal C(x)$ is the smooth matrix constructed in the  following   way: Consider a matrix of size $n\times n$ such that on main diagonal we have $\nu_n(x')$, the $n$th row is $(\nu_1(x'),\dots, \nu_n(x'))$ and the first $n-1$ elements of the last column are $-\nu_1(x'),\dots, -\nu_{n-1}(x'),$ otherwise the element of such a matrix is zero.
If $\nu_n(x')\ne 0 $, then the determinant of such a matrix is not equal zero. Denote the inverse to this matrix  as $\mathcal C_0$
and  set $\mbox{\bf r}=\mathcal C_0 \tilde{ \mbox{\bf v}},$  $\tilde{\mbox{\bf v}}=(v_{1,n},\dots, v_{2,n},\dots v_{n-1,n}, v_2).$  Then $\mathcal C(x')\mbox{\bf v}=\mbox{\bf v}-(\nu_2r_1-\nu_1 r_2, \dots , \nu_n r_1-\nu_1 r_n, \dots , \nu_n r_{n-1}-\nu_{n-1}r_n, \sum_{j=1}^n \nu_jr_j).$
Without loss of generality  we can assume that
\begin{equation}\label{normal}
\vec \nu(0)=-\vec e_n.
\end{equation}
Let all the components of the  function $\mbox{\bf g}$ starting from $n+1$
be equal to zero:
\begin{equation} \label{makaka} g_k=0\quad k\ge n+1.
\end{equation}

\begin{proposition}\label{zoopa} Let $\mbox{\bf v}\in H^1(Q)$ satisfy (\ref{poko}), (\ref{poko1}) and $\mbox{supp}\, \mbox{\bf v}\subset [-T+\epsilon_1,T-\epsilon_1]\times \bar\Omega$ for some positive $\epsilon_1$, $\mbox{\bf g}\in L^2(\Sigma)$.
There exist $\widehat\tau>1$  and $s_0>1$  such that
for any $\tau> \tau_0$ and all $s\ge s_0$
\begin{eqnarray}\label{3.2'1}
s\tau^2\Vert\varphi^\frac 12  \nabla \mbox{\bf v}\,
e^{s\varphi}\Vert^2_{L^2(Q)} +s^3\tau^4\Vert \varphi^\frac 32   \mbox{\bf v}\,
e^{s\varphi}\Vert^2_{L^2(Q)}
+\int_{ \Sigma}(s\tau\varphi \vert \nabla  \mbox{\bf v}\vert^2+s^3\tau^3\varphi^3\vert \mbox{\bf v}\vert^2) e^{2s\varphi}d\Sigma\\
\le C_2(\Vert \mbox{\bf P}(x,D)\mbox{\bf v} e^{s\varphi}\Vert^2_{L^2(Q)}
+\int_{\tilde \Sigma}(s\tau\varphi \vert \nabla  \mbox{\bf v}\vert^2+s^3\tau^3\varphi^3\vert \mbox{\bf v}\vert^2)d\Sigma+s\int_{ \Sigma}\tau\varphi\vert \mbox{\bf g}\vert^2) e^{2s\varphi}d\Sigma),\nonumber
\end{eqnarray}
where $C_2$ is independent of $s$ and $\tau.$
\end{proposition}

First, by an argument based on the partition of unity (e.g., Lemma 8.3.1
in \cite{H}), it suffices to prove the inequality (\ref{3.2'1})
locally, by assuming that
\begin{equation}\label{3.55}
\text{supp}\, \mbox{\bf v}\subset B (y^*,\delta),
\end{equation}
where $B(y^*,\delta)$ is the ball of the radius $\delta>0$
centered at some point $y^*$.

Otherwise, without loss of generality we may assume that $y^*=(y^*_0,0,
\dots, 0).$
Let $\theta\in C^\infty_0(\frac 12,2)$ be a nonnegative function such that
\begin{equation}\label{knight}\sum_{\ell=-\infty}^\infty\theta(2^{-\ell}t)=1
\quad  \mbox{for all}\,\, t\in \Bbb R^1.
\end{equation}

Set $\mbox{\bf v}_\ell(x)=\mbox{\bf v}(x)\kappa_\ell(x),\mbox{\bf g}_\ell(x)=\mbox{\bf g}(x)\kappa_\ell(x)$ where
\begin{equation}\label{knight0}
\kappa_\ell(x)=
\theta (2^{-\ell}e^{\tau^2\psi(x)}).
\end{equation}

Observe that it
suffices to prove the Carleman estimate (\ref{3.2'1}) for the function
$\mbox{\bf v}_\ell$ instead of $\mbox{\bf v}$ provided that the constants
$C_2,\tau_0$  and the function $s_0$ are independent of $\ell.$
Observe that if $G\subset \Bbb R^m$ is a bounded domain and $g\in L^2(G)$,
then there exist positive constants $C_3$ and $C_4$ which are independent of
$g$ (see e.g. \cite{Sog}) such that
\begin{equation}\label{gorokn}
C_3\sum_{\ell=-\infty}^{\infty}\Vert \kappa_\ell g\Vert^2_{L^2(G)}
\le \Vert g\Vert^2_{L^2(G)}\le C_4\sum_{\ell=-\infty}^{\infty}\Vert \kappa_
\ell g\Vert^2_{L^2(G)}.
\end{equation}
Denote the norm on the left-hand side of (\ref{3.2'1}) as $\Vert\cdot\Vert_*.$
By (\ref{knight}) and (\ref{knight0}), we have

$$
\Vert \mbox {\bf v} e^{s\varphi}\Vert_*
=\Vert \sum_{\ell=-\infty}^{+\infty} \mbox{\bf v}_\ell e^{s\varphi}\Vert_*\le
\sum_{\ell=-\infty}^{+\infty}\Vert \mbox {\bf v}_\ell e^{s\varphi}\Vert_*\le
C_5\sum_{\ell=-\infty}^{\infty}(\Vert \kappa_\ell \mbox{\bf P}(x,D)
\text{\bf v}
e^{s\varphi}\Vert^2_{L^2(Q)}
$$
$$+\Vert [\kappa_\ell,\mbox{\bf P}(x,D)]
\text{\bf v}
e^{s\varphi}\Vert^2_{L^2(Q)}
+ \int_{\tilde \Sigma}(s\tau\varphi \vert \nabla  \mbox{\bf v}_\ell\vert^2+s^3\tau^3\varphi^3\vert \mbox{\bf v}_\ell\vert^2+s\tau\varphi\vert \mbox{\bf g}_\ell\vert^2) e^{2s\varphi}d\Sigma
$$
\begin{equation}\label{gomuklus}
+ s\tau\Vert\varphi^\frac 12
[\kappa_\ell, \nabla] \text{\bf v}
e^{s\varphi}\Vert_{{ L}^2(\widetilde\Sigma)}^2)^\frac 12.
\end{equation}

By (\ref{gorokn}), from the above inequality we have
\begin{eqnarray}\label{o}\Vert \mbox {\bf v} e^{s\varphi}\Vert_*\le C_6(\Vert
\mbox{\bf P}(x,D)
\text{\bf v}
e^{s\varphi}\Vert^2_{L^2(Q)}+ \sum_{\ell=-\infty}^{+\infty} \int_{\tilde \Sigma}(s\tau\varphi \vert \nabla  \mbox{\bf v}_\ell\vert^2+s^3\tau^3\varphi^3\vert \mbox{\bf v}_\ell\vert^2) e^{2s\varphi}d\Sigma+
\nonumber
\\
\int_{\tilde \Sigma}s\tau\varphi\vert \mbox{\bf g}\vert^2 e^{2s\varphi}d\Sigma
+\sum_{\ell=-\infty}^{+\infty} \Vert [\kappa_\ell,\mbox{\bf P}(x,D)]
\text{\bf v}
e^{s\varphi}\Vert^2_{L^2(Q)}
+ s\tau\Vert\varphi^\frac 12
[\kappa_\ell, \nabla] \text{\bf v}
e^{s\varphi}\Vert_{{ L}^2(\widetilde\Sigma)}^2)^\frac 12
.
\end{eqnarray}
Let us estimate some terms on the right-hand side of  inequality (\ref{o}).
\begin{eqnarray}\label{oo}\sum_{\ell=-\infty}
^{\infty} s\tau\Vert\varphi
^\frac 12[\kappa_\ell, \nabla] \text{\bf v}
e^{s\varphi}\Vert_{{ L}^2(\widetilde\Sigma)}^2
\le C_7\sum_{\ell
=-\infty}^{\infty} s\tau^5
\Vert\varphi^\frac 12\chi_{\supp\kappa_\ell}\text{\bf v}
e^{s\varphi}\Vert_{{ L}^2(\widetilde\Sigma)}^2
\nonumber
\\
\le C_8 s\tau^5\Vert\varphi^\frac 12\text{\bf v}
e^{s\varphi}\Vert_{{ L}^2(\widetilde\Sigma)}^2.
\end{eqnarray}
Estimating the commutator $ [\kappa_\ell,\mbox{\bf P}(x,D)]$ we obtain
\begin{eqnarray}\label{ooo}
\sum_{\ell=-\infty}^{\infty}\Vert [\kappa_\ell,\mbox{\bf P}(x,D)]\text{\bf v}
e^{s\varphi}\Vert^2_{L^2(Q)}\\\le C_9\sum_{\ell=-\infty}
^{\infty}(\tau^4\Vert \chi_{\supp\kappa_\ell}\nabla\text{\bf v}
e^{s\varphi}\Vert^2_{L^2(Q)}
+\tau^8\Vert\chi_{\supp\kappa_\ell} \text{\bf v}
e^{s\varphi}\Vert^2_{L^2(Q)}).\nonumber
\end{eqnarray}
From (\ref{gomuklus}), (\ref{ooo}), (\ref{oo}) and (\ref{o}) we obtain (\ref{3.2'1}).

Now, thanks to (\ref{3.55}),  without loss of generality we assume that
\begin{equation}\label{3.55}
\text{supp}\, \mbox{\bf v}_\ell \subset B (y^*,\delta)\cap
\mbox{supp}\, \kappa_\ell,
\end{equation}
where $B(y^*,\delta)$ is the ball of the radius $\delta>0$
centered at some point $y^*=(y_0^*,0,\dots,0).$

Assume that near $(0,\dots,0)$,
the boundary $\partial \Omega$ is locally given by an equation
$x_n-\tilde\ell(x_1,\dots, x_{n-1})=0$ and
if $(x_1, \dots, x_n) \in \Omega$, then $x_n -\tilde\ell(x_1,\dots, x_{n-1}) > 0$
where $\tilde\ell\in C^3$
and $\tilde\ell(0) = 0$.  Since $\nu(0)=-\vec e_n$
\begin{equation}\label{napoleon}
\nabla'\tilde \ell(0)=0.
\end{equation}
Denote $F(x)=(x_0,\dots , x_{n-1},x_n -\tilde\ell(x_1,\dots, x_{n-1})) .$
We set
$$
\Delta_{\tilde\ell} u=\sum_{j=1}^{n-1}(\partial^2_{y_j} u
- 2\partial_{x_j}\tilde\ell\circ F^{-1}(y)\partial^2_{y_jy_n} u )
+ (1+\vert \nabla\tilde \ell\vert^2)\partial^2_{y_n} u.
$$
Henceforth we set $y = (y_0, y ')=(y_0,y_1,\dots, y_n)$.
After the change of variables, the  equations (\ref{poko}) have the forms
\begin{equation}\label{3.9}
\mbox{\bf P}(y,D) \mbox{\bf v}=(  \rho\partial^2_{y_0} \mbox{\bf v}_1-\mu \Delta_{\tilde \ell} \mbox{\bf v}_1,\rho\partial^2_{y_0} v_2
- (\lambda +2\mu)\Delta_{\tilde \ell} v_2)
= \mbox{\bf q}, \quad
\mbox{on}\, \mathcal Q\triangleq {\Bbb R}^n\times [0, \gamma] ,    \end{equation}
\begin{equation}\label{!poko1}
\tilde {\mathcal B}(y',D)\mbox{\bf v}=\mbox{\bf g},
\end{equation}
where $\gamma$ is some positive constant.
Here for  functions $\rho\circ F^{-1}(y), \mu\circ F^{-1}(y)$ and
$\lambda \circ F^{-1}(y)$, we use the notations $\rho,\mu,\lambda$
and by  $\mbox{\bf q}_{3}, q_4$ denote the functions $\mbox{\bf q}_3, q_4$
after the change of variables,
respectively. The operator $\tilde {\mathcal B}(y',D)$ is obtained form the
operator ${\mathcal B}(y',D)$ by the change of variables.

Now we introduce operators
\begin{equation}\label{3.18}
P_{\mu}(y,D,\tilde s,\tau) = e^{\vert s\vert\varphi}P_{\mu}(y,D)
e^{-\vert s\vert\varphi},
\quad
P_{\lambda+2\mu}(y,D,\tilde s,\tau)
= e^{\vert s\vert\varphi}P_{\lambda+2\mu} (y,D)
e^{-\vert s\vert\varphi},
\end{equation}
where
$$
\tilde s=s\tau\varphi(y^*).
$$

We denote the principal symbols
of the operators $P_\mu(y,D,\tilde s,\tau)$ and $P_{\lambda+2\mu}(y,D,\tilde s,
\tau)$ by \newline
$p_{\mu}(y,\xi,\tilde s,\tau)$ $=p_{\mu}(y,\xi+{i}\vert s\vert \nabla\varphi)$
and $p_{\lambda+2\mu}(y,\xi,\tilde s,\tau)=p_{\lambda+2\mu}(y,\xi
+{i}\vert s\vert \nabla\varphi)$ respectively.

The principal symbol of the operator $P_{\beta}(y,D,\tilde s,\tau)$ has the
form
\begin{eqnarray}\label{3.34}
p_\beta(y,\xi,\tilde s,\tau) = \rho(y)(\xi_0+{i}\vert s\vert\varphi_{y_0})^2
- \beta\biggl[
\sum_{j=1}^{n-1}(\xi_j+{i}\vert s\vert\varphi_{y_j})^2
                                               \nonumber\\
-2(\nabla'\tilde\ell, (\xi'+{i}\vert s\vert\nabla'\varphi))
(\xi_n+{i}\vert s\vert\varphi_{y_n})
+(\xi_n+{i}\vert s\vert\varphi_{y_n})^2 G\biggr],
\end{eqnarray}
where $ G(y_1,\dots,y_{n-1})=1+\vert\nabla\tilde \ell(y_1,\dots,y_{n-1})\vert^2.$
The zeros of the polynomial $p_{\beta}(y,\xi,\tilde s,\tau)$ with respect to variable $\xi_n$  for
$\vert (\widetilde \xi, \widetilde s)\vert \ge 1$ and $y\in B (y^*,\delta)\cap
\mbox{supp}\, \kappa_\ell$ are
\begin{equation}\label{3.35}
\Gamma^\pm_\beta(y,\xi_0,\dots,\xi_{n-1},\widetilde s,\tau)
= (-{i}\vert \widetilde s\vert\widetilde\mu_\ell(y)\psi_{n}(y)\kappa(\vert (\widetilde\xi,\widetilde s)\vert)
+\alpha^\pm_\beta(y,
\xi_0,\dots,\xi_{n-1},\widetilde s,\tau)),
\end{equation}
where $\vec \psi=(\psi_0,\dots \psi_n),\psi_j= \frac{\varphi(y)}{\varphi(y^*)}
\psi_{y_j}(y)$, \begin{equation}\label{begemot12}
\widetilde\mu_\ell=\eta_*\sum_{k=\ell-2}^{\ell+2}\kappa_\ell,\quad
\eta_*\in C_0^\infty (B(y^*,2\delta)),\quad \eta_*\vert_{B(y^*,\delta)}=1,
\end{equation}
the function $\kappa_\ell$ is given by (\ref{knight0}),
\begin{equation}\label{3.36}
\alpha^\pm_\beta(y,\widetilde \xi,\tilde s,\tau)=\widetilde \mu_\ell
\left(\frac{-\kappa(\vert (\widetilde\xi,\widetilde s)\vert)\sum_{j=1}^{n-1}(\xi_j+{i}\vert \widetilde s\vert
\psi_j)\partial_{y_j}\tilde\ell(y_1,\dots,y_{n-1})}{\vert G\vert}\pm\root
\of{r_\beta(y,\tilde\xi,\tilde s,\tau)}\right),                          \end{equation}
\begin{equation}\label{3.37}
r_\beta(y,\widetilde \xi,\tilde s,\tau)=\kappa(\vert (\widetilde\xi,\widetilde s)\vert)
\frac{(\rho(\xi_0+{i}\vert\widetilde
s\vert\psi_0)^2-\beta\sum_{j=1}^{n-1}(\xi_j+{i}
\vert \widetilde s\vert\psi_j)^2)
G+\beta(\xi'+{i}\vert \widetilde s\vert\vec {\psi},\nabla'\tilde\ell)
^2}
{\beta G^2},
\end{equation}
$\kappa\in C^\infty(\Bbb R^1), \kappa(t)\ge 0, \kappa(t)=1$ for
$t\ge 1$ and $\kappa(t)=0$ for $t\in [0,1/2],$ $\tilde \xi=(\xi_0,\dots, \xi_{n-1}).$

Let $\chi_\nu \in C^\infty_0(\Bbb S^{n})$ on
$\Bbb S^{n}=\{\zeta\triangleq(\widetilde \xi,\widetilde s); \thinspace
\vert (\widetilde \xi, \widetilde s)\vert =1\}$ such that
$\chi_\nu$ is identically equal to $1$ in some neighborhood of
$\zeta^*\in \Bbb S^{n}$ and $\mbox{supp}\, \chi_{\nu}\subset \mathcal O(\zeta^*,\delta).$
Assume that
\begin{equation}\label{book}
\kappa(\vert (\widetilde\xi,\widetilde s)\vert)\vert_{\mbox{supp}\, \chi_\nu}=1,\quad \mbox{supp}\, \kappa(\vert (\widetilde\xi,\widetilde s)\vert)
\subset \mathcal O(\zeta^*,2\delta).
\end{equation}

We extend the function $\chi_\nu$ on $\Bbb R^n$ as follows :
$\chi_\nu(\zeta/\vert(\widetilde\xi,\widetilde s)\vert)$  for
$\vert(\widetilde\xi,\widetilde s)\vert>1$ and \newline
$\chi_\nu(\zeta/\vert(\widetilde\xi, \widetilde s)\vert)\kappa(
\vert (\widetilde\xi,s)\vert)$  for $\vert(\widetilde\xi, \widetilde s)\vert
<1$. Denote by $\chi_\nu(y,\tilde D,\widetilde s)$ the pseudodifferential operator
with the symbol $\eta_\ell(y)\chi_\nu$ and $\eta_\ell(y)=\kappa_{\ell-1}(y)
+\kappa_\ell(y)+\kappa_{\ell+1}(y),$ where $\kappa_\ell$ is
given by (\ref{knight0}).
We set $\mbox{\bf w}=\mbox{\bf v}_\ell e^{s\varphi}$ and $\mbox{\bf w}_\nu
=\chi_\nu(y,\tilde D,\tilde s)\mbox{\bf w}, w_{i,j,\nu}=\chi_\nu(y,\tilde D,\widetilde s)(v_{i,j}\circ F^{-1}(y)).$

Let ${\mathcal O}$ be a domain in ${\Bbb R}^n.$

{\bf Definition.} {\it We say that the symbol
$a(\widetilde y,\widetilde \xi,\tilde  s)\in C^k(\bar{\mathcal O}
\times
{\Bbb R}^{n+1})$ belongs to the class
$C^k_{cl}S^{\kappa,\tilde  s}({\mathcal O})$ if

{\bf A}) There exists a compact set $K\subset\subset{\mathcal O}$
such that $a(\widetilde y,\widetilde \xi,\tilde s)\vert_{{\mathcal O}
\setminus K}=0;$

{\bf B}) For any $\beta=(\beta_0,\dots,\beta_{n})$ there exists a
constant $C_\beta$ such that
$$
\left\Vert \partial^{\beta_0}_{\xi_0}\cdots
\partial^{\beta_{n-1}}_{\xi_{n-1}}
\partial^{\beta_{n}}_{\tilde s}
a(\cdot,\widetilde\xi,\tilde s)\right\Vert_{C^k(\bar{\mathcal O})}\le
C_\beta\left(\tilde s^2+\sum_{i=0}^{n-1}\xi^2_i\right)
^{\frac{\kappa-\vert\beta\vert}{2}}\quad,
$$
where $\vert \beta\vert=\sum_{j=0}^{n}\beta_j$ and
$\vert (\widetilde\xi,\tilde s)\vert \ge 1$;

{\bf C}) For any $N\in \Bbb N$ the symbol $a$ can be represented
as
$$
a(\widetilde y,\widetilde \xi,\tilde s)=\sum_{j=1}^Na_j(\widetilde y,
\widetilde \xi,\tilde s)
+ R_N(\widetilde y,\widetilde \xi,\tilde s),
$$
where the functions $a_j$ have the following properties
$$
a_j(\widetilde y,{\lambda}\widetilde \xi,{\lambda}
\tilde s)={\lambda}^{\kappa-j}a_j(\widetilde y,\widetilde \xi,\tilde  s) \quad
\forall{\lambda}>1,\,\,
\forall(\widetilde y,\widetilde \xi,\tilde  s)\in\{(\widetilde y,
\widetilde \xi,\tilde  s)
\vert \widetilde y\in K, \,\vert (\widetilde\xi,\tilde s)\vert >1\}
$$
$$
\left\Vert  \partial^{\beta_0}_{\xi_0}\cdots
\partial^{\beta_{n-1}}_{\xi_{n-1}}
\partial^{\beta_{n}}_{\tilde s}
a_j(\cdot,\widetilde \xi,\tilde s)\right\Vert_{C^k(\bar{\mathcal O})}\le
C_\beta\left(\tilde s^2+\sum_{i=0}^{n-1}\xi^2_i\right)
^{\frac{\kappa-j-\vert\beta\vert}{2}}
$$
for all $\beta$ and $(\widetilde \xi,\tilde  s)$ satisfying
$\vert(\widetilde\xi,\tilde s)\vert\ge 1$
and the term $R_N$ satisfies the estimate
$$
\Vert R_N(\cdot,\widetilde\xi,\tilde s)\Vert_{C^k(\bar{\mathcal O})}\le
C_N(\tilde s^2+\sum_{i=0}^{n-1}\xi^2_i)^{\frac{\kappa-N}{2}}\quad
\forall (\widetilde\xi,\tilde  s)\,\, \mbox{satisfying}\,\, \vert
(\widetilde\xi,\tilde s)
\vert \ge 1.
$$}
Obviously
\begin{equation}
\pi_{C^k(B(0,\delta(y^*)))}(\chi_\nu)\le C_{10}\tau^{2k}\quad \forall
k\in\Bbb N_+.
\end{equation}
Obviously the  pseudodifferential operators with the symbols
$\Gamma_\beta^\pm$ belongs \\
to the class $C^2_{cl}S^{1,s}(B(0,\delta(y^*)))$
and
\begin{equation}
\pi_{C^2(B(0,\delta(y^*)))}(\Gamma_\beta^\pm)\le C_{11}\tau^4.
\end{equation}

We have
\begin{proposition} Let $w\in H^1(\mathcal Q),$ $\mbox{supp}\, w\subset
B(y^*,\delta)\cap \mbox{supp}\,\eta_\ell$  and $P_\beta(y,D,\tilde s,\tau)
\chi_\nu w\in L^2(\mathcal Q).$ Then
there exist positive constants $\delta(y^*), C_{12}, C_{13}$
independent of $s$
and $\tau$  and  independent constants $s_0,\tau_0$ independent of $s$ such that for all $\tau\ge \tau_0$
and $s\ge s_0$ we have
\begin{eqnarray}\label{klop}
C_{12}\int_{\mathcal Q}(\vert s\vert\tau^2\varphi\vert \nabla\chi_\nu w\vert^2
+ \vert s\vert^3\tau^4\varphi^3\vert \chi_\nu w\vert^2)dy
+ \varXi_\beta(\chi_\nu w)\\
\le  \Vert P_\beta(y,D,\tilde s,\tau)\chi_\nu w\Vert^2_{L^2(\mathcal Q)}
+ C_{13} \epsilon(
\tau_0)
\int_{\Bbb R^n}(\vert s\vert\tau\varphi(y^*)\vert \nabla \chi_\nu w\vert^2
+\vert s\vert^3\tau^3\varphi^3(y^*)\vert\chi_\nu w\vert^2 )(\widetilde y,0)
d\widetilde y,\nonumber
\end{eqnarray}
where $\epsilon({\tau_0})\rightarrow +0$ as ${\tau_0}\rightarrow
+\infty$ and
$$
\varXi_\beta(w)=\sum_{j=1}^3 \frak I_j(\beta,w),\quad \frak I_1(\beta,w)
=\int_{\Bbb R^n}
(\vert\widetilde  s\vert\beta^2(y^*)\psi_{y_n}(y^*)\vert\partial_{y_n}
w\vert^2 +\vert\widetilde
s\vert^3\beta^2(y^*)\psi^3_{y_n}(y^*)\vert w\vert^2)d\widetilde y,
$$
\begin{equation}\label{01}\frak I_2(\beta,w)=- \frac 12 Re
\int_{\Bbb R^n} 2\vert \widetilde s\vert\beta(y^*)\partial_{y_n} w
\overline{(\nabla_{\widetilde \xi}p_\beta(y^*,\widetilde \nabla w,0),
\widetilde \nabla \psi(y^*))} d\widetilde y,\end{equation}
\begin{equation}\label{02} \frak I_3(\beta,w)=\int_{\Bbb R^n}\vert
\widetilde s\vert\beta(y^*)\psi_{y_n}(y^*)(p_\beta(y^*,\widetilde \nabla w,0)
-\widetilde s^2p_\beta(y^*,\widetilde\nabla\psi(y^*),0)\vert
 w\vert^2) d\widetilde y.
\end{equation}
\end{proposition}

{\bf Proof.} It suffices to prove the statement of the lemma separately for
$\mbox{Re}\, w_\nu, \mbox{Im}\, w_\nu.$  Let  $v_\nu=\mbox{Re}\, w_\nu$ or
$v_\nu=\mbox{Re}\, w_\nu.$   For simplicity instead of $p_\beta (y,\xi)$
we use the notation
$p(y,\xi)=\sum_{k,j=0}^np_{kj}(y)\xi_k\xi_j,$ where $p_{kj}=p_{jk}$ for all
$k,j\in\{0,\dots, n\}$.
We set
$$
p_j(y,\xi)=\partial_{y_j}p(y,\xi),
p^{(j)}(y,\xi)=\partial_{\xi_j}p(y,\xi),
p^{(j,m)}(y,\xi)=\partial^2_{\xi_j \xi_m}p(y,\xi), p(y,\eta,\xi)=\sum_{k,j=0}^n
p_{kj}(y)
\eta_i\xi_j.
$$
Observe that since $y\in  B(y^*,\delta)\cap \mbox{supp}\, \eta_\ell$ then either
$\frac 12\le 2^{-\ell}e^{\tau^2\psi(y)}\le 2$ or $\frac 12\le 2^{-\ell-1}
e^{\tau^2\psi(y)}\le 2$ or $\frac 12\le 2^{-\ell+1}e^{\tau^2\psi(y)}\le 2.$

This is equivalent to
\begin{equation}\label{koko}
(\ell-2)\ln 2/\tau^2 \le \psi(y)\le(\ell+2)\ln 2/\tau^2\quad \forall y
\in  B(y^*,\delta)\cap \mbox{supp}\, \eta_\ell.
\end{equation}
Hence, for any $\epsilon$ there exists $\tau_0(\epsilon)$ such that for all
$\tau\ge \tau_0$
\begin{equation}\label{kokom}
\vert\varphi(y)-\varphi(y^*)\vert \le \epsilon \varphi(y^*)\quad \forall y
\in  B(y^*,\delta)\cap \mbox{supp}\, \eta_\ell.
\end{equation}
Indeed, since $\psi(y)>C_{14}>0$ on $B(0,2\delta),$  by (\ref{koko}), there exists an independent
constant $C_{15}$ such that
\begin{equation}\label{monarch}
\vert \ell\vert /\root\of{\tau}\le C_{15}\quad \mbox{if}\,\, B(y^*,\delta)\cap \mbox{supp}\, \eta_\ell\ne \{\emptyset\}.
\end{equation}
Then estimate (\ref{kokom}) follows  from the inequality
\begin{equation}\label{nonsense}
\vert\varphi(y)-\varphi(y^*)\vert=e^{\tau \psi(y^*)}\vert 1
-e^{\tau (\psi(y^*)-\psi(y))}\vert\le\varphi(y^*)\vert 1-e^\frac {(2\vert\ell\vert+4)
\ln 2}{\tau} \vert\le \varphi(y^*)\vert 1-e^\frac {C_{16}}{\root\of{\tau}}
\vert .
\end{equation}
In order to get the last two inequalities in (\ref{nonsense}) we used (\ref{monarch}) and (\ref{koko}).

We introduce the operators
\begin{equation}\label{poker}
L_1(y,D,\tilde s,\tau) v_{\nu}=-\sum_{k=0}^n
s\varphi_{y_k}p^{(k)}(y,\nabla v_{\nu}), \quad L_2(y,D,\tilde s,\tau)v_{\nu}
= P(y,D)v_{\nu}
+s^2p(y,\nabla \varphi,\nabla\varphi) v_{\nu}.
\end{equation}
Denote $\mbox{\bf f}_s=P(y,D,\tilde s,\tau)v_{\nu}-sv_{\nu}P(y,D)\varphi.$ Then we have
\begin{eqnarray}\label{okop}\Vert \mbox{\bf f}_s\Vert^2
_{L^2(\mathcal Q)}
=\Vert L_1(y,D,\tilde s,\tau)v_{\nu}\Vert^2_{L^2(\mathcal Q)}\nonumber\\
+ \Vert L_2(y,D,\tilde s,\tau)w_
\nu\Vert^2_{L^2(\mathcal Q)}+2\mbox{Re}\,(L_1(y,D,\tilde s,\tau)v_{\nu},
L_2(y,D,\tilde s,\tau)v_{\nu})_{L^2(\mathcal Q)}.
\end{eqnarray}
The following equality is proved in Imanuvilov \cite{Im}:
\begin{eqnarray}\label{LK}
\mbox{Re}\,(L_1(y,D,\tilde s,\tau)v_{\nu},L_2(y,D,\tilde s,\tau)v_{\nu})
_{L^2(\mathcal Q)}
                                       \nonumber\\
= -\mbox{Re}\,\int_{\partial\mathcal Q}p(y,\vec e_n,\nabla v_{\nu})
{L_1(y,D,\tilde s,\tau)v_{\nu}} d\Sigma-s\int_{\partial\mathcal Q}
p(y,\vec e_n,\nabla
\varphi)p(y,\nabla v_{\nu},\nabla v_{\nu}) d\Sigma\nonumber\\
+s^3\int
_{\partial\mathcal Q}p(y,\nabla\varphi,\nabla \varphi)p(y,\vec e_n,\nabla \varphi)
\vert v_{\nu}\vert^2 d\Sigma
+ s\int_{\mathcal Q}\mathcal G(y,\tilde s,\tau,v_{\nu})dy+\\
 \int_{\mathcal Q}\frac{s}{2}(\sum_{k,m=0}^n p^{(k)}_k(y,\nabla v_{\nu})
\partial_{y_m}\varphi p^{(m)}(y,\nabla v_{\nu})-\theta(y)(p(y,\nabla v_{\nu},\nabla v_{\nu})
-s^2p(y,\nabla\varphi,\nabla\varphi)\vert v_{\nu}\vert^2))dy,\nonumber
\end{eqnarray}
where
$$
\mathcal G(y,\tilde s,\tau,w)=\{p,\{p,\varphi\}\}(y,\nabla w)
+ s^2\sum_{k,j=0}^np_{j}
(y,\nabla \varphi) w^2+s^2\sum_{k,j=0}^n\partial^2_{y_ky_j}\varphi p^{(k)}
(y,\nabla\varphi)p^{(j)}(y,\nabla\varphi) w^2
$$
and
$$
\theta(y)=\sum_{k,m=0}^n(\partial^2_{y_ky_j}\varphi p^{(k,m)}(y,\nabla \varphi)
+\partial_{y_k}\varphi p_m^{(k,m)}(y,\nabla v_{\nu})).
$$
Observe that the function $\theta(y)$ is independent of $v_\nu$ and
\begin{equation}\label{ok}
\sup_{y\in \mbox{supp}\, \eta_\ell}\vert \theta(y)\vert
\le C_{16}\tau^2\varphi(y^*).
\end{equation}
Without loss of generality we may assume that $\partial_{y_n} \psi(y^*)\ne 0.$
We introduce the form $\frak G(\tilde s,\tau,\nabla v_{\nu})$ in the following
way:
In the function $\mathcal G(y^*,\tilde s,\tau,\nabla v_{\nu})$ we replace
$\partial_{y_n}
v_{\nu}$ \newline by $-\frac{1}{\sum_{j=0}^np_{jn}(y^*)
\partial_{y_j}\varphi(y^*)}
\sum_{j,k=0}^{n-1}\partial_{y_j}\varphi(y^*)p_{jk}(y^*)\partial_{y_k}v_{\nu}.$
Since $p_{kj}(y^*)=0$ for $k\ne j$
\begin{equation}\label{prizrak}
\sum_{j=0}^np_{jn}(y)\partial_{y_j}\varphi(y)
\ne 0\quad \forall y\in B(y^*,\delta)\cap \mbox{supp}\, \eta_\ell.
\end{equation}
By (\ref{kokom}) we have the inequality
\begin{equation}\label{papuas1}
\left\vert \int_{\mathcal Q} \mathcal G(y,\tilde s,\tau,\nabla v_{\nu})
- \frak G(\tilde s,\tau,\nabla v_{\nu}) dy\right\vert
\le \frac{C_{17}}{s^2}\Vert L_1(y,D,\tilde s,\tau)v_{\nu}\Vert^2
_{L^2(\mathcal Q)}
+\epsilon(\tau_0)\Vert v_{\nu}\Vert^2_{H^{1,\widetilde s}(\mathcal Q)},
\end{equation}  where $\epsilon(\tau_0)\rightarrow +0$ as $\tau_0\rightarrow +\infty.$

Let
\begin{equation}\label{propoganda}
\xi = \left(\widetilde\xi,-\frac{1}{\sum_{j=0}^np_{jn}(y^*)\partial_{y_j}
\varphi(y^*)}\sum_{j,k=0}^{n-1}\partial_{y_j}\varphi(y^*)p_{jk}(y^*)\xi_k
\right).
\end{equation}
We set
$$
q(\widetilde \xi,s)=\sum_{k,j=0}^n \partial^2_{y_ky_j}\varphi(y^*)p^{(k)}(y^*,\xi+{i} s\nabla\varphi(y^*))
\overline{p^{(j)}(y^*,\xi+{i} s\nabla\varphi(y^*))}
$$
$$+\frac 1s \mbox{Im}\sum_{k=0}^np_k(y^*,\xi+{i} s\nabla\varphi(y^*))
\overline{p^{(k)}(y^*,\xi+{i} s\nabla\varphi(y^*))}.
$$
Observe that
\begin{equation}
\int_{\Bbb R^{n+1}_+}q(\widetilde \xi,s)\vert \mbox{F}_{\widetilde y
\rightarrow \widetilde \xi}v_{\nu}\vert^2d\widetilde \xi dy_n=\int_{\mathcal Q}
\frak G(\tilde s,\tau,\nabla v_{\nu}) dy.
\end{equation}
Denote $\widetilde w(\tilde \xi, y_n)=(2\pi)^{-\frac{n}{2}}\int_{\Bbb R^n}sign(\mbox{Re}\,
p(y^*,\xi+{i} s\nabla \varphi(y^*))\mbox {F}_{\widetilde y\rightarrow
\widetilde \xi}
v_{\nu} e^{{i}<\widetilde \xi,\widetilde y>}d\widetilde \xi$,
where $\mbox {F}_{\widetilde y\rightarrow \widetilde \xi}$
is the Fourier transform given by
$$
F_{\widetilde y\rightarrow
\widetilde \xi}u=\frac{1}{(2\pi)^\frac{n}{2}}\int_{{\Bbb R}^{n}}
e^{-{i}\sum_{j=0}^{n-1}y_j\xi_j}u(y_0,\dots,y_{n-1})d\widetilde y.
$$

Taking the scalar product of the function $L_2(y,D,\tilde s,\tau)v_{\nu}$ and
$\overline{\widetilde w}$  in  $L^2(\mathcal Q)$ we have
\begin{eqnarray}\label{begemot1}
\int_{\mathcal Q}(-p(y,\nabla v_{\nu},\overline{\nabla\widetilde  w})
+ s^2p(y,\nabla \varphi,\nabla \varphi)v_{\nu}\overline{\widetilde w}dy
=-\int_{\Bbb R^n}\partial_{\vec\nu_p} v_{\nu}(\widetilde y,0)
\overline{\widetilde w(\widetilde y,0)}d\widetilde y\nonumber\\
-\sum_{k,j=0}^n\int_{\mathcal Q}\partial_{y_k} p_{kj}
\partial_{y_j} v_{\nu}\overline{\widetilde w}dy
+ (L_2(y,D,\tilde s,\tau)v_{\nu},\overline{\widetilde w})_{L^2(\mathcal Q)},
\end{eqnarray}
where $\partial_{\vec\nu_p} w
=-\sum_{j=0}^np_{nj}\partial_{y_j}w.$

By (\ref{kokom}) for any  $\epsilon \in (0,1)$ there exists $\tau_0(\epsilon)$
such that  for all $\tau\ge \tau_0$ and all $s\ge 1$
\begin{eqnarray}\label{begemot2}
\int_{\mathcal Q}\widetilde s \tau s^2\vert p(y^*,\nabla\varphi (y^*),\nabla\varphi
(y^*))-p(y,\nabla\varphi ,\nabla\varphi)\vert \vert v_{\nu}\vert^2dy
                                           \nonumber\\
+\int_{\mathcal Q}\tilde s\tau \vert p(y,\nabla v_{\nu},\overline{\nabla \widetilde w})
- p(y^*,\nabla v_{\nu},\overline{\nabla \widetilde w})\vert dy
\le \epsilon \widetilde s \tau\Vert v_{\nu}\Vert^2_{H^{1,\tilde s}(\mathcal Q)}.
\end{eqnarray}
The inequalities  (\ref{begemot1})  and (\ref{begemot2}) imply that for any
positive $\epsilon\in (0,1)$ there exists  a constant $C_{18}(\epsilon)$ such that
\begin{eqnarray}
\widetilde s\tau\left\vert \int_{\mathcal Q}(-p(y^*,\nabla v_{\nu},
\overline{\nabla\widetilde  w})+s^2p(y^*,\nabla \varphi(y^*),
\nabla \varphi(y^*))v_{\nu}\overline{\widetilde w}dy\right\vert\\
\le C_{19}\left\vert\widetilde s\tau \int_{\Bbb R^n}\partial_{\vec\nu_p} v_{\nu}(\widetilde y,0)\overline{\widetilde w}(\widetilde y,0)d
\widetilde y\right\vert
+ \epsilon\Vert L_2(y,D,\tilde s,\tau)v_{\nu}\Vert^2_{L^2(\mathcal Q)}
+ (\epsilon \tilde s\tau +C_{18} \tau^2)\Vert v_{\nu}\Vert^2_{H^{1,\widetilde s}(\mathcal Q)}.
                                    \nonumber
\end{eqnarray}
We set
$\widehat
\nabla \widetilde w=(\partial_{y_0}\widetilde w, \dots, \partial_{y_{n-1}}
\widetilde w,\frac{-1}{\sum_{j=0}^np_{jn}(y^*)\partial_{y_j}\varphi(y^*)}
\sum_{j,k=0}^{n-1}\partial_{y_j}\varphi(y^*)p_{jk}(y^*)\partial_{y_k}
\widetilde w).$
Hence, if $\xi$ is given by (\ref{propoganda}), then we have
\begin{eqnarray}\label{LL0}
\int_{\Bbb R^n_+}\widetilde s\tau\vert p(y^*, \xi+{i} s\nabla \varphi(y^*))
\vert^2 \vert\mbox{F}_{\widetilde y\rightarrow \widetilde \xi} v_{\nu}\vert^2
d\widetilde \xi dy_n\\
= \left\vert \widetilde s\tau\int_{\mathcal Q}(-p(y^*,\widehat\nabla v_{\nu},
\overline{\widehat \nabla\widetilde  w})+s^2p(y^*,\nabla \varphi(y^*),\nabla
\varphi(y^*))v_{\nu}\overline{\widetilde w}dy\right\vert
                                                  \nonumber\\
\le C_{20}\widetilde s\tau \int_{\Bbb R^n}\left\vert\partial_{\vec\nu_p} v_{\nu}(\widetilde y,0)\overline{\widetilde w(\widetilde y,0)}
\right\vert
d\widetilde y+\epsilon\Vert L_2(y,D,\tilde s,\tau)v_{\nu}\Vert^2
_{L^2(\mathcal Q)}
+(\epsilon\tilde s\tau +C_{21} \tau^2)\Vert v_{\nu}\Vert^2_{H^{1,\widetilde s}(\mathcal Q)}.
                                                      \nonumber
\end{eqnarray}

Observe that pseudoconvexity Condition \ref{A1} implies that there exists a
positive constant $C_{22}$ such that
\begin{equation}\label{LL1}
q(\widetilde \xi,s)+\tau^2\varphi(y^*)\vert p(y^*,\widetilde \xi+{i} s\nabla
\varphi(y^*))\vert\ge C_{22}\tau^2\varphi(y^*)\vert
(\widetilde\xi,\widetilde s)
\vert^2 \quad \forall (\tilde\xi,\widetilde s)\in \Bbb S^{n}.
\end{equation}
Therefore, from (\ref{LL1}) and (\ref{LL0}),
for some positive  constant $C_{23}$  for all $\tau\ge \tau_0$ and all $s\ge 1$ we have the inequality
\begin{eqnarray}\label{papuas}
C_{23}\widetilde s\tau\int_{\mathcal Q} (\vert \widetilde \nabla v_{\nu}
\vert^2
+\widetilde s^2\vert v_{\nu}\vert^2)dy\le C_{24}(\widetilde s\tau \int
_{\Bbb R^n}
\vert \partial_{\vec\nu_p} v_{\nu}(\widetilde y,0)
\overline{\widetilde w(\widetilde y,0)}\vert d\widetilde y\nonumber\\
+\epsilon\Vert L_2(y,D,\tilde s,\tau)v_{\nu}\Vert^2_{L^2(\mathcal Q)})
+ \int_{\mathcal Q}  s
\frak G(\tilde s,\tau,\nabla v_{\nu}) dy+(\epsilon\tilde s\tau +C_{25} \tau^2) \Vert v_{\nu}\Vert^2
_{H^{1,\widetilde s}(\mathcal Q)}.
\end{eqnarray}
We take $\epsilon\in (0, \mbox{min}\{\frac{C_{23}}{4}, \frac{1}{8C_{24}}\})$.
Thanks to (\ref{LK}), (\ref{papuas1}), (\ref{papuas}) for some positive  constant $C_{26}$  for all $\tau\ge \tau_0$ and all $s\ge 1$ we have the inequality
\begin{eqnarray}\label{mina1}
C_{26}\widetilde s\tau\int_{\mathcal Q}
(\vert \widetilde \nabla v_{\nu}\vert^2
+\widetilde s^2\vert v_{\nu}\vert^2)dy
\le C_{27}\widetilde s\tau \int_{\Bbb R^n}
\vert\partial_{\vec\nu_p} v_{\nu}(\widetilde y,0)\overline{
\widetilde w(\widetilde y,0)}\vert d\widetilde y
+ \frac 18\Vert L_2(y,D,\tilde s,\tau)w_
\nu\Vert^2_{L^2(\mathcal Q)}              \nonumber\\
+2(L_1(y,D,\tilde s,\tau)v_{\nu}, L_2(y,D,\tilde s,\tau)v_{\nu})
_{L^2(\mathcal Q)}+\mbox{Re}\,
\int_{\partial\mathcal Q}p(y,\vec e_n,\nabla v_{\nu}){L_1(y,D,\tilde s,\tau)
v_{\nu} }
d\Sigma\nonumber\\-s\int_{\partial\mathcal Q}p(y,\vec e_n,\nabla \varphi)p(y,\nabla
v_{\nu},{\nabla v_{\nu}}) d\Sigma+s^3\int_{\partial\mathcal Q}
p(y,\nabla\varphi,\nabla \varphi)p(y,\vec e_n,\nabla \varphi)
\vert v_{\nu}\vert^2 d
\Sigma                                                 \nonumber\\
-\int_{\mathcal Q}\frac{s}{2}(\sum_{k,m=0}^n p^{(k)}_k(y,\nabla v_{\nu})
p^{(m)}(y,\nabla v_{\nu})-\theta(y)(p(y,\nabla v_{\nu})
-s^2p(y,\nabla\varphi)
\vert v_{\nu}\vert^2))dy.
\end{eqnarray}

Thanks to (\ref{prizrak}) and (\ref{poker}), the  following identity holds true
$$
\partial_{y_n} v_{\nu}=\frac{L_1(y,D,\tilde s,\tau) v_{\nu}/s-\sum_{k=0}^n
\sum_{j=0}^{n-1}p_{k,j}\partial_{y_k}\varphi \partial_{y_j} v_{\nu}}
{\sum_{k=0}^n p_{j,n}\partial_{y_k}\varphi}.
$$
Therefore there exists a constant $C_{28}$ independent of $\widetilde s,\tau$
such that
\begin{equation}\label{mina}
\widetilde s\tau\int_{\mathcal Q} \vert \partial_{y_n} v_{\nu}\vert^2 dy
\le C_{28}(\int_{\mathcal Q} \widetilde s\tau (\vert \widetilde \nabla v_{\nu}\vert^2
+\widetilde s^2
\vert v_{\nu}\vert^2)dy+\Vert L_1(y,D,\tilde s,\tau)v_{\nu}\Vert^2
_{L^2(\mathcal Q)}).
\end{equation}

Using (\ref{okop}), from (\ref{mina}) and (\ref{mina1}) for any $\tau\ge \tau_0$ and $s>1$ we obtain the estimate
\begin{eqnarray}\label{logotip}
C_{29}\widetilde s\tau\int_{\mathcal Q} (\vert  \nabla v_{\nu}\vert^2
+\widetilde s^2
\vert v_{\nu}\vert^2)dy
\le C_{30}\widetilde s\tau\left\vert \int_{\Bbb R^n}
\partial_{\vec\nu_p} v_{\nu}(\widetilde y,0)
\overline{\widetilde w(\widetilde y,0)}d\widetilde y\right\vert
                                         \nonumber\\
+
\int_{\partial\mathcal Q}p(y,\vec e_n,\nabla v_{\nu}){L_1(y,D,\tilde s,\tau)
v_{\nu}}
d\Sigma\\-s\int_{\partial\mathcal Q}p(y,\vec e_n,\nabla \varphi)
p(y,\nabla v_{\nu},{\nabla v_{\nu}}) d\Sigma+s^3\int_{\partial\mathcal Q}
p(y,\nabla\varphi,\nabla \varphi)p(y,\vec e_n,\nabla \varphi)\vert v_{\nu}
\vert^2
d\Sigma                       \nonumber\\
-\int_{\mathcal Q}\frac{s}{2}\left(\sum_{k,m=0}^n p^{(k)}_
k(y,\nabla v_{\nu})\partial_{y_m}\varphi
p^{(m)}(y,\nabla v_{\nu})-\theta(y)(p(y,\nabla v_{\nu},\nabla  v_{\nu})
-s^2p(y,\nabla\varphi,\nabla\varphi)\vert v_{\nu}\vert^2)\right)dy
                                            \nonumber\\
+\Vert \mbox{\bf f}_s\Vert^2_{L^2(\mathcal Q)}-\frac 12\Vert
L_1(y,D,\tilde s,\tau)
v_{\nu}\Vert^2_{L^2(\mathcal Q)}-\frac{1}{2}\Vert L_2(y,D,\tilde s,\tau)v_{\nu}\Vert^2
_{L^2(\mathcal Q)}.\nonumber
\end{eqnarray}
Now we estimate some integrals on the right hand side of (\ref{logotip}):
There exist a constant $C_{31}$ such that
\begin{equation}\label{puk1}
\vert \int_{\mathcal Q}\frac{s}{2}\sum_{k,m=0}^n p^{(k)}_k(y,\nabla v_{\nu})\varphi_{y_m}
p^{(m)}(y,\nabla v_{\nu})dy\vert\le C_{31}\int_{\mathcal Q}\vert\nabla
v_{\nu}
\vert^2dy+\frac 18\Vert L_1(y,D,\tilde s,\tau)v_{\nu}\Vert^2
_{L^2(\mathcal Q)}.
\end{equation}
Integrating by parts we have
\begin{eqnarray}
\int_{\mathcal Q}\theta(y)(p(y,\nabla v_{\nu},\nabla v_{\nu})-s^2p(
y,\nabla\varphi,\nabla\varphi)\vert v_{\nu}\vert^2)dy=-\int _{\mathcal Q}
\theta(y)L_2(y,D,\tilde s,\tau)v_{\nu}{v_{\nu}}dy\nonumber\\ +\int_{\Bbb R^n}
\theta\partial_{\vec\nu_p} v_{\nu}(\widetilde y,0){v_{\nu}(
\widetilde y,0)}d\widetilde y+\sum_{j,k=0}^n\int_{\mathcal Q}(\theta p_k^{(j)}
(y,\nabla v_{\nu}){v_{\nu}}+p^{(j)}(y,
\nabla v_{\nu})\partial_{y_k}\theta{v_{\nu}})dy.\nonumber
\end{eqnarray}
Therefore (\ref{ok}) yields
\begin{eqnarray}\label{puk}
s\left\vert \int_{\mathcal Q}\theta(p(y,\nabla v_{\nu},\nabla  v_{\nu})
-s^2p(y,\nabla\varphi)\vert v_{\nu}\vert^2)dy\right\vert
\le \frac 18\Vert L_2(y,D,\tilde s,\tau)
v_{\nu}\Vert^2_{L^2(\mathcal Q)}\nonumber\\
+ C_{32}\tau^3\Vert v_{\nu}\Vert^2_{H^{1,\widetilde s}(\mathcal Q)}
+ C_{33}\Vert (\partial_{y_n} v_{\nu}(\cdot,0),v_
\nu(\cdot,0))\Vert^2
_{L^2(\Bbb R^n)\times H^{1,\widetilde s}(\Bbb R^n)}.
\end{eqnarray}
Using (\ref{puk}) and (\ref{puk1}), from (\ref{logotip}) for any $\tau\ge \tau_0$ and $s>1$ we obtain the estimate
\begin{eqnarray}\label{logotip1}
C_{34}\widetilde s\tau\int_{\mathcal Q} (\vert  \nabla v_{\nu}\vert^2
+\widetilde s^2\vert v_{\nu}\vert^2)dy
\le C_{35}\Vert (
\partial_{y_n} v_{\nu}(\cdot,0),v_{\nu}(\cdot,0)
)\Vert^2
_{L^2(\Bbb R^n)\times H^{1,\widetilde s}(\Bbb R^n)}\nonumber\\+
\int_{\partial\mathcal Q}p(y,\vec e_n,\nabla v_{\nu}){L_1(y,D,\tilde s,\tau)
v_{\nu}}
d\Sigma\\-s\int_{\partial\mathcal Q}p(y,\vec e_n,\nabla \varphi)
p(y,\nabla v_{\nu},{\nabla v_{\nu}}) d\Sigma+s^3\int_{\partial\mathcal Q}
p(y,\nabla\varphi,\nabla \varphi)p(y,\vec e_n,\nabla \varphi)\vert v_{\nu}
\vert^2 d\Sigma
\nonumber\\
+\Vert P(y,D,\tilde s,\tau)v_{\nu}\Vert^2_{L^2(\mathcal Q)}
-\frac 34\Vert L_1(y,D,\tilde s,\tau)v_{\nu}
\Vert^2_{L^2(\mathcal Q)}-\frac 12\Vert L_2(y,D,\tilde s,\tau)v_{\nu}\Vert^2
_{L^2(\mathcal Q)}.\nonumber
\end{eqnarray}
Next we estimate the difference between the boundary integrals and
$\sum_{j=1}^3 \frak I_j(\beta,v_{\nu})$. Using (\ref{kokom}) we have
\begin{eqnarray}\label{logotip2}
\biggl\vert -\int_{\partial\mathcal Q}p(y,\vec e_n,\nabla v_{\nu})
{L_1(y,D,\tilde s,\tau)v_{\nu}} d\Sigma-s\int_{\partial\mathcal Q}
p(y,\vec e_n,\nabla
\varphi)p(y,\nabla v_{\nu},{\nabla v_{\nu}}) d\Sigma\nonumber\\
+s^3\int_{\partial\mathcal Q}p(y,\nabla\varphi,\nabla \varphi)p(y,\vec e_n,\nabla
\varphi)\vert v_{\nu}\vert^2 d\Sigma-\sum_{j=1}^3 \frak I_j(\beta,v_{\nu})
\biggr\vert                                             \nonumber\\
\le C_{36} \epsilon(\tau_0)
\int_{\Bbb R^n}(\vert s\vert\tau\varphi(y^*)\vert \nabla v_{\nu}\vert^2
+\vert s\vert^3\tau^3\varphi^3(y^*) \vert v_{\nu}\vert^2 )(\widetilde y,0)
d\widetilde y,
\end{eqnarray} where $\epsilon({\tau_0})\rightarrow +0$ as ${\tau_0}\rightarrow
+\infty.$
From (\ref{logotip1}) and (\ref{logotip2}), we obtain (\ref{klop}).
$\blacksquare$

In some cases, we can represent the operator $
P_{\beta}(y,D,\tilde s,\tau)$ as a product of two first order
pseudodifferential operators.
\begin{proposition}\label{gorokx1} Let $\beta\in \{\mu,\lambda+2\mu\},$ $r_\beta(y^*,\zeta^*,\tau)\ne 0$
and $ supp\, {\chi_\nu}\subset \mathcal O(\delta_1).$
Then we can
factorize the operator $P_{\beta}(y,D,\tilde s,\tau)$ into the product of two
first order pseudodifferential operators:
\begin{equation} \label{min}
P_{\beta}(y,D,\tilde s,\tau)w_{\nu}=\beta G(\frac 1i\partial_{y_n}-\Gamma^-_
\beta(y,{\widetilde D},\tilde s,\tau))
(\frac 1i\partial_{y_n}-\Gamma^+_\beta(y,{\widetilde D},\tilde s,\tau))
w_{\nu}+T_\beta w,
\end{equation}
where $T_\beta : H^{1,\widetilde s}(\mathcal Q)
\rightarrow L^2(\mathcal Q)$ is a continuous operator and there
exists a constant $C_{37}$ independent of $\tilde s$ and $\tau$ such that
\begin{equation}\label{lokom}
\Vert T_\beta w_{\nu}\Vert_{L^2(\mathcal Q)}\le C_{37}\tau^2\Vert w\Vert_{ H^{1,\widetilde s}(\mathcal Q)}.
\end{equation}
\end{proposition}
{\bf Proof.} Let
$$
\widetilde R(y,\tilde D,\tilde s,\tau)=\left [\frac{\rho(\frac 1 i\partial_{y_0}+{i}
\vert
\widetilde s
\vert\psi_0)^2}{\beta G}-\frac{\sum_{j=1}^{n-1}(\frac 1i\partial_{y_j}+{i}\vert
\widetilde s\vert\psi_j)^2}{G}\right]
$$
and  $\Gamma(y,{\widetilde D},\tilde s,\tau)$ is the operator with symbol
$\Gamma^-_\beta(y,{\widetilde \xi},\tilde  s,\tau)\Gamma^+_\beta
(y,{\widetilde \xi},\tilde s,\tau):$
$$
\Gamma(y,\widetilde \xi,\tilde s,\tau)=(-\vert \widetilde s\vert^2(
\widetilde \mu_\ell\psi_n\kappa)^2+\alpha^+_\beta {i}\vert \widetilde s\vert
\psi_n\kappa+\alpha_\beta^-{i}\vert \widetilde s\vert\psi_n\kappa+\alpha_\beta^+
\alpha_\beta^-)
$$
$$
= \kappa^2(\vert(\widetilde \xi,\widetilde s)\vert)\left [\frac{\rho(\xi_0
+{i}
\vert\widetilde s\vert\psi_0)^2}{\beta G}-\frac{\sum_{j=1}^{n-1}
(\xi_j+{i}\vert\widetilde s\vert\psi_j)^2}{G}\right].
$$ Functions $\tilde \mu_\ell$ and $\kappa$ are given by (\ref{begemot12}) and (\ref{book}) respectively.
We set $\Upsilon_\ell=B(y^*,2\delta)\cap\text{supp}\,\widetilde \mu_\ell.$
Then
\begin{eqnarray}\label{lesopoval1}
\Gamma(y,{\widetilde D},\tilde s,\tau)w_{\nu}=[\Gamma,\eta_\ell]
\chi_\nu(\tilde D,\tilde s)w
+\eta_\ell\Gamma(y,\tilde D,\tilde s,\tau)\chi_\nu(\tilde D,\tilde s)w               \\
= [\Gamma,\eta_\ell]\chi_\nu(\tilde D,\tilde s)w + \eta_\ell \widetilde
R(y,\tilde D,\tilde  s,\tau)\chi_\nu(\tilde D,\tilde  s)w
                                                \nonumber\\
=[\Gamma,\eta_\ell]\chi_\nu(\tilde D,\tilde s)w+\widetilde R(y,\tilde D,\tilde s,\tau)
w_{\nu}+[\eta_\ell,
\widetilde R(y,\tilde D,\tilde s,\tau)]\chi_\nu(\tilde D,\tilde s)w.\nonumber
\end{eqnarray}
In order to obtain the second equality in
(\ref{lesopoval1}) we used (\ref{book}).  The short computations imply
\begin{equation}\label{lecopoval}
(\frac 1i\partial_{y_n}-\Gamma^-_\beta(y,{\widetilde D},\tilde s,\tau))
(\frac 1i\partial_{y_n}-\Gamma^+_
\beta(y,{\widetilde D},\tilde s,\tau))
\end{equation}
\begin{eqnarray*}
= &&-\partial_{y_n}^2-\frac 1i[\partial_{y_n},\Gamma^+_
\beta(y,{\widetilde D},\tilde s,\tau)]+\Gamma^-_\beta(y,{\widetilde D},
\tilde s,\tau)\Gamma^+_\beta(
y,{\widetilde D},\tilde s,\tau)  \\
+&& {i}\Gamma^-_\beta(y,{\widetilde D},\tilde s,\tau)\partial_{y_n}
+{i}\Gamma^+_\beta(y,{\widetilde D},\tilde s,\tau)\partial_{y_n}.
\end{eqnarray*}
By Lemma \ref{Fops3}, we have
$$
\Gamma^-_\beta(y,{\widetilde D},\tilde s,\tau)\Gamma^+_\beta(y,{\widetilde D},
\tilde s,\tau)
= \Gamma(y,{\widetilde D},\tilde s,\tau)+R_0,
$$
where
$$\Vert R_{0}\Vert_{{\mathcal L}(H_0^{1,\tilde s}({\Upsilon_\ell}),L^2({\Upsilon_\ell}))}\le C_{38}\pi_{C^1(\Upsilon_\ell)}
(\Gamma^+_\beta)\pi_{C^1(\Upsilon_\ell)}(\Gamma^-_\beta)
\le C_{39}\tau^4.
$$
%
The commutator
$
[\partial_{y_n},\Gamma^+_\beta(y,{\widetilde D},\tilde s,\tau)]
$ is the  pseudodifferential operator with the symbol \newline
$\partial_{y_n}\Gamma^+_\beta(y,{\widetilde \xi},\tilde s,\tau)$.
By Lemma \ref{Fops0} we have
$$
\Vert [\partial_{y_n},\Gamma^+_\beta(y,{\widetilde D},\tilde s,\tau)]\Vert
_{\mathcal L (H^{1,\tilde s}_0(\Upsilon_\ell),L^2(\Upsilon_\ell))}
\le C_{40}\pi_{C^1(\Upsilon_\ell)} ([\partial_{y_n},\Gamma^+_\beta
(y,{\widetilde D},\tilde s,\tau)])\le C_{41}\tau^4.
$$
Denote
$$
R(y,\tilde D,\tilde s,\tau) = \left(2\vert \widetilde s\vert \psi_n+
\frac{\sum_{j=1}^{n-1}\partial_{y_j}\tilde\ell(y_1,\dots,y_{n-1})(\partial_{y_j}
-\vert \widetilde s\vert \psi_j)}{ G}\right).
$$
By (\ref{3.34})-(\ref{3.36}), (\ref{book}) and the fact that
$\tilde\mu_\ell\eta_\ell=\eta_\ell$ the following is true:
\begin{eqnarray}\label{lecopoval3}
({i}\Gamma^-_\beta(y,{\widetilde D},\tilde s,\tau)\partial_{y_n}
+{i}\Gamma^+_\beta(y,{\widetilde D},\tilde s,\tau)\partial_{y_n}) w_{\nu}
=\tilde\mu_\ell R(y,\tilde D,\tilde s,\tau)\partial_{y_n}\kappa^2(\tilde D,\widetilde s)
w_{\nu}\\
=\tilde\mu_\ell  [R(y,\tilde D,\tilde s,\tau)\partial_{y_n}\kappa^2,\eta_\ell]
\chi_\nu(\tilde D,\widetilde s)w
+\tilde\mu_\ell\eta_\ell
R(y,\tilde D,\tilde s,\tau)\partial_{y_n}\chi_\nu(\tilde D,\widetilde s)w      \nonumber\\
=\tilde\mu_\ell  [R(y,\tilde D,\tilde  s,\tau)\partial_{y_n}\kappa^2,\eta_\ell]
\chi_\nu(\tilde D,\widetilde s)w
+[\eta_\ell ,R(y,\tilde D,\tilde s,\tau)\partial_{y_n}]\chi_\nu(\tilde D,\widetilde s)w
+R(y,\tilde D,\tilde s,\tau)\partial_{y_n}w_{\nu}.
                             \nonumber
\end{eqnarray}
Since $-\partial^2_{y_n}w_{\nu}+R(y,\tilde D,\tilde s,\tau)\partial_{y_n}w_{\nu}
+ \widetilde R(y,\tilde D,\tilde s,\tau)
\partial_{y_n}w_{\nu}=\frac{1}{\beta G}P_\beta(y,D,\tilde s,\tau)w_{\nu}$,
setting
$$
T_\beta =-R_0+[\partial_{y_n},\Gamma^+_\beta(y,{\widetilde D},\tilde s,\tau)]
\chi_\nu(y,\tilde D,\tilde s)
-[\Gamma,\eta_\ell]\chi_\nu(\tilde D,\tilde s)+[\eta_\ell, \widetilde R(y,\tilde D,\tilde s,
\tau)]\chi_\nu(\tilde D,\tilde s)
$$
$$
-\tilde\mu_\ell  [R(y,\tilde D,\tilde  s,\tau)\partial_{y_n}\kappa^2,\eta_\ell]\chi_\nu
(\tilde D,\widetilde s)
-[\eta_\ell ,R(y,\tilde D,\tilde  s,\tau)\partial_{y_n}]\chi_\nu(\tilde D,\widetilde s)
$$
and using (\ref{lesopoval1}) - (\ref{lecopoval3}), we obtain (\ref{min}).
Now we prove (\ref{lokom}). Lemma \ref{Fops2} yields
\begin{eqnarray}\label{AQ1}
\Vert [\Gamma,\eta_\ell]\Vert
_{\mathcal L (H^{1,\tilde  s}_0(\Upsilon_\ell),L^2(\Upsilon_\ell))}
\le C_{42}(\pi
_{C^0(\Upsilon_\ell)}(\Gamma)\pi_{C^0( \Upsilon_\ell)}(\eta_\ell)+\pi
_{C^0(\Upsilon_\ell)}(\Gamma)\pi_{C^1(\Upsilon_\ell)}(\eta_\ell)\nonumber\\+\pi
_{C^1(\Upsilon_\ell)}(\Gamma)\pi_{C^0(\Upsilon_\ell)}(\eta_\ell))\le C_{43} \tau^2.
\end{eqnarray}
For differential operators $R$ and $\tilde R$,  we obtain the
estimates
\begin{equation}\label{LQ2}
 \Vert [\mu_\ell, \widetilde R(y,\tilde D,\tilde s,\tau)]\Vert
_{\mathcal L (H^{1,\tilde s}
_0(\Upsilon_\ell),L^2(\Upsilon_\ell))}\le C_{44}\tau^2,
\end{equation}
and
\begin{eqnarray}\label{Q3}
\Vert\mu_\ell  [R(y,\tilde D,\tilde s,\tau)\partial_{y_n}\kappa^2,\eta_\ell]\Vert
_{\mathcal L (H^{1,\tilde s}_0(\Upsilon_\ell),L^2(\Upsilon_\ell))}
                                     \nonumber\\
\le \Vert
  [R(y,\tilde D,\tilde s,\tau)\partial_{y_n}\kappa^2,\eta_\ell]\Vert
_{\mathcal L (H^{1,\tilde s}_0(\Upsilon_\ell),L^2(\Upsilon_\ell))}
\le C_{45}\tau^2,
\end{eqnarray}
and
\begin{equation}\label{LQ1}
\Vert [\eta_\ell ,R(y,\tilde D,\tilde s,\tau)\partial_{y_n}]\Vert _{\mathcal L
 (H^{1,\tilde  s}_0(\Upsilon_\ell),L^2(\Upsilon_\ell))}
\le C_{46} \Vert \eta_\ell\Vert
_{C^2(\Upsilon_\ell)}\le C_{45}\tau^2.
\end{equation}
Form (\ref{AQ1})-(\ref{LQ1}) we obtain (\ref{lokom}).
The proof of the proposition is complete. $\blacksquare$

Now we apply the Lemma \ref{Fops4} to estimate $(1-\eta_\ell)\chi_\nu(\tilde D,
\tilde s) \mbox{\bf w}$. Set $\mathcal O_1
=\mbox{supp}\, \kappa_\ell\cap B(y^*,\delta) $, $\mathcal O=B(0,2R)$ there radius $R$
is sufficiently large, $\mathcal O_2=supp(1-\eta_\ell)\cap B(0,R).$ By (\ref{3.55}) and definition of the function $\mbox{\bf w}$ $\mbox{supp}\, \mbox{\bf w}\subset \mathcal O_1.$ Since $\eta_\ell=1$ on $\mbox{supp }\,\kappa_\ell$ we have $\mathcal O_1\cap \mathcal O_2=\emptyset.$
Then  for some positive constant $C_{55}$ $dist (\mathcal O_1,\mathcal O_2)\ge \frac {C_{46}}{\tau^2}.$
By (\ref{gabon})
\begin{equation}\label{po}
\vert \tilde s\vert\Vert (1-\eta_\ell)\chi_\nu(\tilde D, \tilde s) \mbox{\bf w}
\Vert_{H^{1,\tilde s} (B(0,R))}\le C_{47}\tau^{4n+6}\Vert \mbox{\bf w}\Vert
_{L^2(\Bbb R^n)}.
\end{equation}

By arguments, same as in proposition \ref{Fops4} we obtain
\begin{equation}\label{po1}
\vert \tilde s\vert\Vert (1-\eta_\ell)\chi_\nu(\tilde D, \tilde s)
\mbox{\bf w}\Vert_{H^{1,\tilde s} (\Bbb R^n\setminus B(0,R))}\le{ C_{48}}
\Vert \mbox{\bf w}\Vert_{L^2(\Bbb R^n)}.
\end{equation}
Hence the inequalities (\ref{po}) and (\ref{po1}) imply
\begin{equation}\label{nikolas}
\vert \tilde s\vert\Vert (1-\eta_\ell)\chi_\nu(\tilde D, \tilde s) \mbox{\bf w}
\Vert_{H^{1,\tilde s} (\Bbb R^n)}\le C_{49}\tau^{4n+6}\Vert
\mbox{\bf w}\Vert_{L^2(\Bbb R^n)}.
\end{equation}

Denote
$$
V^\pm_{\mu}(k,j)=(\frac 1i \partial_{y_n}-\Gamma^\pm_{\mu}(y,{\widetilde D},
\tilde s,\tau))
w_{k,j,\nu}, \quad
V^\pm_{\lambda+2\mu}=(\frac 1i \partial_{y_n}-\Gamma^\pm_{\lambda+2\mu}
(y,{\widetilde D},\tilde s,\tau))w_{2,\nu}.
$$

Let us consider the equation
\begin{equation}\label{barmalei0}
(\frac 1i \partial_{y_n}-\Gamma^-_\beta(y,{\widetilde D},\tilde s,\tau))V=q,
\quad V\vert_{y_n=\gamma}=0.
\end{equation}
For solutions of this problem, we can prove an a priori estimate.
\begin{proposition}\label{gorokx} Let $\beta\in \{\mu,\lambda+2\mu\},$ $r_\beta(y^*,\zeta^*,\tau)\ne 0$
and $V=V^+_{\mu}(k,j)$ if $\beta =\mu$ and $V=V^+_{\lambda+2\mu}$ if $\beta
=\lambda+2\mu.$ There exists a
constant $C_{50}>0$ such that
\begin{equation} \label{min1}
\root\of{\vert \widetilde s\vert}\Vert V(\cdot, 0)\Vert_{L^2({\Bbb R}^n)}\le
C_{50}(\tau^{4n+6}\Vert \mbox{\bf w}\Vert_{H^{1,\tilde s}(\mathcal Q)}+\Vert q\Vert_{L^2(\mathcal Q)}).
\end{equation}
\end{proposition}

${\bf Proof.}$ We taking the scalar product of the equation (\ref{barmalei0})
and the function $-i\bar V$ in $L^2(\mathcal Q)$ and integrating by parts
we obtain
\begin{equation}\label{fosset}
\frac 12\Vert V(\cdot, 0)\Vert^2_{L^2(\Bbb R^n)}+i\int_0^\gamma(\Gamma^-_\beta
(y,{\widetilde D},\tilde s,\tau)V,\bar V)_{L^2(\Bbb R^n)}dy_n=-i\int_0^\gamma(q,\bar V)
_{L^2(\Bbb R^n)}dy_n.
\end{equation}
If $\tilde s^*\ne 0$ or $\mbox{Im}\,\root \of{ r_\beta(y^*,\zeta^*,\tau)}\ne 0$.
Then  for some positive constant $C_{51}$
\begin{equation}\label{3.38}
-\text{Im}\, \vert s\vert \Gamma^-_\beta(y,\tilde \xi,\widetilde s,\tau)
\ge \widehat C_{51}\vert(\tilde \xi, \widetilde s)\vert^2
\quad \forall (y,\tilde \xi,\widetilde s)\in \Upsilon_\ell\times
\mathcal O(y^*,\delta_1(y^*)).
\end{equation}
By (\ref{nikolas}) and G\aa rding's inequality (\ref{goblin}) there exists a positive constant $C_{52}$
such that
$$
-\int_0^\gamma\vert \tilde s\vert\mbox{Im}\,(\Gamma^-_\beta(y,{\widetilde D},\tilde s,\tau)V,\bar V)
_{L^2(\Bbb R^n)}dy_n\ge C_{52}\int_0^\gamma\Vert V\Vert^2_{H^{1,\tilde s}(\Bbb R^n)}
dy_n-C_{53} \tau^{4n+6}  \Vert V\Vert^2_{L^2(\mathcal Q)}.
$$
This inequality, (\ref{fosset}) and Proposition \ref{gorokx1} imply
(\ref{min1}).
Now let $\tilde s^*=\mbox{Im}\,\root \of{ r_\beta(y^*,\zeta^*,\tau)}= 0.$
Denote $\vert \tilde s\vert b(y,\zeta,\tau)
= \text{Im}\,r_\beta(y,\zeta,\tau),$
$a(y,\zeta,\tau) = \text{Re}\,r_\beta(y,\zeta,\tau)$,
$\zeta=(\xi_0,\dots,\xi_{n-1},\tilde s)$.
If $\text{Im}\,\root\of{r_\beta(y^*,\zeta^*,\tau)}= 0$, then we have
$a(y^*, \zeta^*,\tau)>0$.
In that case, since near $(y^*, \zeta^*)$, we have
$r_{\beta}(y,\zeta,\tau) = a(y,\zeta,\tau)\left(
1 + \frac{{i\vert\tilde s\vert} b(y,\zeta,\tau)}
{a(y,\zeta,\tau)}
\right)$ and $\left\vert \frac{\vert \tilde s\vert b(y,\zeta,\tau)}
{a(y,\zeta,\tau)}\right\vert < 1$ on $\Upsilon_\ell \times
\mathcal O(y^*,\delta(y^*))$
we may define the function
$\root\of{r_\beta(y,\zeta,\tau)}$ by the infinite series for the
function $(1+t)^\frac 12=\sum_{n=0}^\infty c_n t^n$, $c_n
=\frac{\frac12(\frac{1}{2}-1)(\frac{1}{2}-2)\dots (\frac{1}{2}-(n-1))}{n!}$
which holds true for $\vert t\vert<1$.
On $\Upsilon_\ell\times \mathcal O(y^*,\delta(y^*))$
we set
\begin{equation}\label{3.39}
\root\of{r_\beta(y,\zeta,\tau)}=\root\of{a}\sum_{n=0}^\infty c_n
\left (\frac{{i}\vert \tilde s\vert b}{a}\right)^n
= \root\of{a} + \frac{{i}}{2}\tilde s\left(\frac{b}{\root\of{a}}\right)
+\root\of{a}\sum_{n=2}^\infty c_n
\left (\frac{{i} \vert s\vert b}{a}\right )^n.
\end{equation}
Since $\partial_{y_n}\varphi(y^*)<0$ the exists a constant $C_{54}$ such that
\begin{eqnarray}\label{Leva}
\mbox{Re}\left\{i\int_0^\gamma(\Gamma^-_\beta(y,{\widetilde D},\tilde s,\tau)V,\bar V)
_{L^2(\Bbb R^n)}dy_n\right\}\nonumber \\
\ge \mbox{Re}\left\{i\int_0^\gamma(\alpha^-_\beta
(y,{\widetilde D},\tilde s,\tau)V,\bar V)_{L^2(\Bbb R^n)}dy_n\right\}
\nonumber\\
+ C_{54}\vert\tilde s\vert \Vert V\Vert^2_{L^2(\mathcal Q)}.
\end{eqnarray}
Let $A$ be the pseudodifferential operator with the symbol $a(y,\zeta,\tau)$
and $A^\frac 12 $ the pseusodifferential operator with the  symbol
$a^\frac 12(y,\zeta,\tau).$ By Lemma \ref{Fops3}  $A=A^\frac12 A^\frac 12
+R_0 $ where $\Vert R_0\Vert_{\mathcal L(L^2(\mathcal Q),L^2(\mathcal Q))}\le \tau^2C_{55}.$ By Lemma \ref{Fops1}  $(A^\frac 12)^*
=A^{\frac 12}+R_1$ where $\Vert R_1\Vert_{\mathcal L(L^2(\mathcal Q),L^2(\mathcal Q))}\le \tau^2C_{56}.$
Hence
\begin{eqnarray}\label{Leva1}
\mbox{Re}\left\{i\int_0^\gamma(A(y,{\widetilde D},\tilde s,\tau)V,\bar V)
_{L^2(\Bbb R^n)}dy_n\right\}=\mbox{Re}\left\{
i\int_0^\gamma(A^\frac 12(y,{\widetilde D},\tilde s,\tau)V,A^\frac 12
(y,{\widetilde D},\tilde s,\tau)\bar V)_{L^2(\Bbb R^n)}dy_n\right\}
\nonumber\\
+\mbox{Re}\left\{i\int_0^\gamma(A^\frac 12(y,{\widetilde D},\tilde s,\tau)V,
(2R_0+R_1)\bar V)_{L^2(\Bbb R^n)}dy_n\right\}\\
= \mbox{Re}\left\{i\int_0^\gamma(A^\frac 12(y,{\widetilde D},\tilde s,\tau)V,
(2R_0+R_1)\bar V)_{L^2(\Bbb R^n)}dy_n\right\}\ge -\tau^2C_{57}  \Vert V\Vert^2_{L^2(\mathcal Q)}.\nonumber
\end{eqnarray}
Let $B$ be the pseudodifferential operator with the symbol $\frac{{i}\vert
\tilde s\vert b}{\root\of{a}}+\root\of{a}\sum_{n=2}^\infty c_n
\left (\frac{{i} \vert s\vert b}{a}\right )^n.$
By G\aa rding's inequality (\ref{goblin})  for any positive $\epsilon$ there exists a constant $C_{58}(\epsilon)$
such that
\begin{equation}\label{Leva2}
\mbox{Re}\left\{\vert \tilde s\vert i\int_0^\gamma(B(y,{\widetilde D},\tilde s,\tau)
V,\bar V)_{L^2(\Bbb R^n)}dy_n\right\}
\ge- \epsilon \vert\tilde s\vert^2 \Vert V\Vert^2_{L^2(\mathcal Q)}
-C_{58}\tau^{4n+6} \Vert V\Vert^2_{L^2(\mathcal Q)}.
\end{equation}
Let $C(x,D,\tilde s,\tau)$ be the pseudodifferential operator with the symbol $\widetilde \mu_\ell
\frac{-\kappa(\vert (\widetilde\xi,\widetilde s)\vert)\sum_{j=1}^{n-1}(\xi_j+{i}\vert \widetilde s\vert
\psi_j)\partial_{y_j}\tilde\ell}{\vert G\vert}.$  By (\ref{napoleon})  and  G\aa rding's inequality (\ref{goblin})  for any positive $\epsilon$  taking the raius of the ball $\delta$ sufficiently small
one can find a constant $C_{59}(\epsilon,\delta)$ such that
\begin{equation}\label{Leva3}
\mbox{Re}\left\{\vert \tilde s\vert i\int_0^\gamma(C(y,{\widetilde D},\tilde s,\tau)
V,\bar V)_{L^2(\Bbb R^n)}dy_n\right\}
\ge- \epsilon \vert\tilde s\vert^2 \Vert V\Vert^2_{L^2(\mathcal Q)}
-C_{59}\tau^{4n+6} \Vert V\Vert^2_{L^2(\mathcal Q)}.
\end{equation}

Inequalities (\ref{Leva1}), (\ref{Leva}), (\ref{Leva2}) and (\ref{Leva3}) imply (\ref{min1}).
$\blacksquare$


We will separately consider the two cases
$r_\mu(y^*,\zeta^*,\tau)=0$ in Section \ref{section3} and $r_{\la+2\mu}(y^*,\zeta^*,\tau)=0$
in Section \ref{section4}.
\bigskip

\section{ Case $r_\mu(y^*,\zeta^*,\tau)=0.$}\label{section3}

In this section, we mainly treat the case when
$\text{supp}\,\chi_\nu\subset \mathcal O(y^*,\delta_1(y^*)),$ and
$(y^*,\zeta^*)$ be a point on $\Bbb R^{n+1}\times \Bbb S^{n}$ such that
$r_\mu(y^*,\zeta^*,\tau)=0.$
By (\ref{3.35})-(\ref{3.37}) there exists $C_1>0$ such that
\begin{eqnarray}\label{2.7} \vert p_\mu(y^*,\widetilde \xi,0)- \widetilde
s^2 p_\mu(y^*,\widetilde \nabla \psi(y^*),0)\vert+\vert\widetilde  s(\nabla
_{\widetilde \xi}p_\mu(y^*,\widetilde \xi,0 ), \widetilde \nabla \psi(y^*))
\vert\nonumber\\ \le \delta_1C_1\vert(\widetilde \xi,\widetilde s)\vert^2,
\quad \forall\zeta\in\mathcal
O(\zeta^*,\delta_1(y^*)).
\end{eqnarray}
Hence, by (\ref{2.7}) and (\ref {nikolas}) for some independent constants
$C_2,C_3$
\begin{eqnarray}\label{2.99}
\vert \frak J_3(\mu, w_{k,j,\nu})\vert\le C_2\delta_1\vert \tilde s\vert\Vert
(\partial_{y_n} \mbox{\bf w}_{\nu}(\cdot,0),\mbox{\bf  w}
_{\nu}(\cdot,0))\Vert^2
_{L^2(\Bbb R^n)\times H^{1,\widetilde s}(\Bbb R^n)}\nonumber\\+C_3\tau^{4n+6}
\Vert
(\partial_{y_n} \mbox{\bf w}(\cdot,0),\mbox{\bf  w}(\cdot,0))\Vert^2
_{L^2(\Bbb R^n)\times H^{1,\widetilde s}(\Bbb R^n)}.
\end{eqnarray}
We recall that by (\ref{klop}) there exist $C_4>0$ and $C_5>0$ such
that
\begin{eqnarray}\label{2.1}C_4(\vert \widetilde s\vert \tau\Vert  w_{k,j,\nu}
\Vert^2_{H^1(\mathcal Q)}+\vert \widetilde s\vert^3\tau\Vert
w_{k,j,\nu}\Vert^2_{L^2(\mathcal Q)})+\varXi_\mu(w_{k,j,\nu})\le C_5\Vert
P_\mu(y,D,\tilde s,\tau) w_{k,j,\nu} \Vert^2_{L^2(\mathcal Q)}\nonumber\\
+\epsilon(\delta)\vert \widetilde s\vert\Vert
(\partial_{y_n} \mbox{\bf w}_{\nu}(\cdot,0),\mbox{\bf  w}
_{\nu}(\cdot,0))\Vert^2
_{L^2(\Bbb R^n)\times H^{1,\widetilde s}(\Bbb R^n)},
\end{eqnarray}
where $\epsilon(\delta)\rightarrow 0$ as $\delta\rightarrow +0.$
If $1\le k<j<n$, then we have
\begin{eqnarray}\label{yoyo}
\varXi_\mu(w_{k,j,\nu})\ge\int_{{\Bbb R^n}}\left(
\vert\widetilde  s\vert\mu^2(y^*)\psi_{y_n}(y^*)
\left\vert\partial_{y_n} w_{k,j,\nu}\right\vert^2
+ \vert\widetilde s\vert^3\mu^2(y^*)
\psi_{y_n}^3(y^*)\vert
w_{k,j,\nu}\vert^2 \right)(\tilde y,0) d\widetilde y\nonumber\\
-(\epsilon(\delta)\vert\widetilde s\vert +C_6)\Vert
(\partial_{y_n} \mbox{\bf w}_{\nu}(\cdot,0),\mbox{\bf  w}
_{\nu}(\cdot,0))\Vert^2
_{L^2(\Bbb R^n)\times H^{1,\widetilde s}(\Bbb R^n)}.
\end{eqnarray}

Indeed, making the change of variables in the function $v_{kj}(x)$ and using the fact that all components of the function $\mbox{\bf g}$ starting from $n+1$ are zeros we have
$$
\tilde v_{k,j}(y)=v_{k,j}\circ F^{-1}(y)=(\nu_k r_j-\nu_j r_k)\circ F^{-1}(y).
$$
Then

$$
w_{k,j}(y)=v_{kj}\circ F^{-1}(y)=((\nu_k r_j-\nu_j r_k)e^{\vert s\vert \varphi})\circ F^{-1}(y).
$$
So
$$
\Vert w_{k,j}\Vert_{H^{1,\tilde s}(\Bbb R^n)}\le C_7\sum_{j=1}^{n-1}\Vert \nu_j\Vert_{C^0(\overline{B(y^*,\delta)})}(\Vert \nabla'\mbox{\bf w}_1\Vert_{L^2(\Bbb R^n)}+s\Vert \varphi  \mbox{\bf w}_1\Vert_{L^2(\Bbb R^n)})
$$
$$+C_8
\sum_{j=1}^{n-1}\Vert \nu_j\Vert_{C^1(\overline{B(y^*,\delta)})}\Vert \mbox{\bf w}_1\Vert_{L^2(\Bbb R^n)}.
$$
Since $\nu(0)=-\vec e_n$ this inequality and Lemma \ref{Fops0} implies (\ref{yoyo}).

{\bf Case A.}  Assume that $ \widetilde s^*\ne 0.$
The parameter $\delta_1>0$ will be fixed later.

Inequality (\ref{2.7}) yields
\begin{eqnarray}\label{2.9}
\vert \frak J_2(\mu, w_{j,n,\nu})\vert
\le \frac{C_9\delta_1\mu(y^*)\vert\widetilde s\vert}{\vert\widetilde s^*\vert}
\Vert
(\partial_{y_n} \mbox{\bf w}_{\nu}(\cdot,0),\mbox{\bf  w}
_{\nu}(\cdot,0))\Vert^2
_{L^2(\Bbb R^n)\times H^{1,\widetilde s}(\Bbb R^n)}\nonumber\\+C_
{10}\tau^{4n+6}
\Vert
(\partial_{y_n} \mbox{\bf w}(\cdot,0),\mbox{\bf  w}(\cdot,0))\Vert^2
_{L^2(\Bbb R^n)\times H^{1,\widetilde s}(\Bbb R^n)}.
\end{eqnarray}

Applying the Cauchy-Bunyakovskii inequality and using (\ref{2.9}) and
(\ref{2.99}), we see that there exists some positive
constants $C_{11},C_{12}$ such that
\begin{eqnarray}\label{2.4}
\varXi_\mu(w_{j,n,\nu})\ge
\int_{{\Bbb R^n}} \left(\vert\widetilde s\vert\mu^2(y^*)\psi_{y_n}(y^*)
\vert\partial_{y_n}w_{j,n,\nu}\vert^2
+ \vert\widetilde s\vert^3\mu(y^*)\psi_{y_n}(y^*)\vert
w_{j,n,\nu}\vert^2\right)(\tilde y,0) d\widetilde y         \nonumber        \\
- \frac{C_{11} \delta_1\mu(y^*)}{\vert\widetilde s^*\vert}\vert\widetilde s\vert
\Vert
(\partial_{y_n} \mbox{\bf w}_{\nu}(\cdot,0),\mbox{\bf  w}
_{\nu}(\cdot,0))\Vert^2
_{L^2(\Bbb R^n)\times H^{1,\widetilde s}(\Bbb R^n)}\\-C_{12}\tau^{4n+6}
\Vert(\partial_{y_n} \mbox{\bf w}(\cdot,0),\mbox{\bf  w}(\cdot,0))\Vert^2
_{L^2(\Bbb R^n)\times H^{1,\widetilde s}(\Bbb R^n)}.\nonumber
\end{eqnarray}
If in addition $r_{\lambda+2\mu}(y^*,\zeta^*,\tau)=0$, then
$$
\xi_0^*=\psi_{y_0}(y^*)=0\quad\mbox{and}\quad \widetilde s^*\ne 0.
$$

Hence, similarly to (\ref{2.4}), we obtain
\begin{equation} \label{2.444}
\varXi_{\lambda+2\mu}(w_{2,\nu}) \ge C_{13}\int_{{\Bbb R^n}} (\vert \widetilde
s \vert
\vert\nabla w_{2,\nu}\vert^2+\vert \widetilde s\vert^3
\vert w_{2,\nu}\vert^2)(\tilde y,0) d\widetilde y
\end{equation}
$$
- \epsilon(\delta)\vert\widetilde s\vert \Vert
(\partial_{y_n} \mbox{\bf w}_{\nu}(\cdot,0),\mbox{\bf  w}
_{\nu}(\cdot,0))\Vert^2
_{L^2(\Bbb R^n)\times H^{1,\widetilde s}(\Bbb R^n)}-C_{14}\tau^{4n+6}\Vert
(\partial_{y_n} \mbox{\bf w}(\cdot,0),\mbox{\bf  w}(\cdot,0))\Vert^2
_{L^2(\Bbb R^n)\times H^{1,\widetilde s}(\Bbb R^n)}.
$$

Combining estimates  (\ref{2.4}) and (\ref{2.444}), we have
\begin{eqnarray} \label{pulemet}
\vert\widetilde s\vert \Vert
(\partial_{y_n} \mbox{\bf w}_{\nu}(\cdot,0),\mbox{\bf  w}
_{\nu}(\cdot,0))\Vert^2
_{L^2(\Bbb R^n)\times H^{1,\widetilde s}(\Bbb R^n)} + \vert\widetilde s\vert
\tau\Vert
\mbox{\bf w}_{\nu}\Vert^2_{H^{1,\widetilde s}(\mathcal Q)}\le
C_{15}(\tau^{4n+6}\Vert \text{\bf w}
\Vert^2_{H^{1,\widetilde s}(\mathcal Q)}                   \nonumber\\
+ \vert\widetilde s\vert\Vert
\mbox{\bf g}e^{\vert s\vert\varphi}\Vert^2_{L^2({\Bbb R^n})}
+\Vert P_{\mu}(y,D,\tilde s,\tau)\mbox{\bf
w}_{1,\nu}\Vert^2_{L^2(\mathcal Q)}+\Vert P_{\lambda+2\mu}(y,D,\tilde s,\tau)
w_{2,\nu}\Vert^2_{L^2(\mathcal Q)}\nonumber\\
+\tau^{4n+6}\Vert
(\partial_{y_n} \mbox{\bf w}(\cdot,0),\mbox{\bf  w}(\cdot,0))\Vert^2
_{L^2(\Bbb R^n)\times H^{1,\widetilde s}(\Bbb R^n)}).
\end{eqnarray}

If $r_{\lambda+2\mu}(y^*,\zeta^*,\tau)\ne 0$, then by Proposition \ref{gorokx}
there exists a constant $C_{16}$ independent of $s$ and $\tau$ such that
\begin{equation}\label{ix}
\root\of{1+\vert \widetilde s\vert}\Vert (\frac 1i \partial_{y_n} w_{2,\nu}
-\Gamma^+
_{\lambda+2\mu}(y,{\widetilde D},\tilde s,\tau))w_{2,\nu}\vert_{y_n=0}\Vert
_{L^2(\Bbb R^n)}
\end{equation}
$$
\le C_{16}(\Vert P_{\lambda+2\mu}(y,D,\tilde s,\tau)w_{2,\nu}\Vert
_{L^2(\mathcal Q)}
+ \tau^{2n+3}\Vert w_{2,\nu}\Vert_{H^{1,\widetilde s}(\mathcal Q)}).
$$

On the other hand, on $\Bbb R^n$ from the boundary condition (\ref{poko1}), $ b_n(y,D)\mbox{\bf  w}=g_n$ we have
\begin{equation}\label{ix1}
((\lambda +2\mu)(y^*)\left(\partial_{y_n} w_{2,\nu}
-\vert \widetilde s\vert\psi_{y_n}(y^*)w_{2,\nu}\right)
- \mu(y^*)\sum_{j=1}^{n-1} \left(\partial_{y_j}
w_{j,n,\nu}-\vert\widetilde s\vert\psi_{y_j}(y^*)w_{j,n,\nu}
\right))(\cdot,0)
= \mbox{\bf r},
\end{equation}
where the function $\mbox{\bf r}$ satisfies
\begin{eqnarray}\label{garmodon}
\Vert \mbox{\bf r}\Vert^2_{L^2(\Bbb R^n)}
\le \epsilon(\delta)\Vert
(\partial_{y_n} \mbox{\bf w}_{\nu}(\cdot,0),\mbox{\bf  w}
_{\nu}(\cdot,0))\Vert^2
_{L^2(\Bbb R^n)\times H^{1,\widetilde s}(\Bbb R^n)} \nonumber\\+ \frac{C_{17}
\tau^{4n+6}}{1+\vert \tilde s\vert}\Vert
(\partial_{y_n} \mbox{\bf w}(\cdot,0),\mbox{\bf  w}(\cdot,0))\Vert^2
_{L^2(\Bbb R^n)\times H^{1,\widetilde s}(\Bbb R^n)}
+ C_{18}\Vert \mbox{\bf g}e^{\vert s\vert\varphi}\Vert^2_{L^2(\Bbb R^n)}
\end{eqnarray}
with some constants $C_{17}, C_{18}$ independent of $\widetilde s$ and $\tau$.

From (\ref{ix1}) and (\ref{2.4}), it follows that
\begin{eqnarray}\label{02.4}
\vert\widetilde s\vert\left\Vert(\lambda +2\mu)(y^*)
\left(\partial_{y_n} w_{2,\nu}-\vert\widetilde s\vert
\psi_{y_n}(y^*)w_{2,\nu}\right)(\cdot,0)\right\Vert^2_{L^2(\Bbb R^n)} \\
\le C_{19}\sum_{j=1}^{n-1}\varXi_\mu(w_{j,n,\nu}) + \epsilon(\delta)
\vert\widetilde s\vert \Vert
(\partial_{y_n} \mbox{\bf w}_{\nu}(\cdot,0),\mbox{\bf  w}
_{\nu}(\cdot,0))\Vert^2
_{L^2(\Bbb R^n)\times H^{1,\widetilde s}(\Bbb R^n)}\nonumber\\
+C_{20}\tau^{4n+6}\Vert
(\partial_{y_n} \mbox{\bf w}(\cdot,0),\mbox{\bf  w}
(\cdot,0))\Vert^2
_{L^2(\Bbb R^n)\times H^{1,\widetilde s}(\Bbb R^n)}
+ C_{21}\vert\widetilde s\vert\Vert \mbox{\bf g}
e^{\vert s\vert\varphi}
\Vert_{L^2(\Bbb R^n)}^2.\nonumber
\end{eqnarray}
Then this estimate, the G\aa rding inequality (\ref{goblin}) and (\ref{ix})
imply
\begin{eqnarray}\label{3.8}
\vert\widetilde s\vert^3\Vert \psi_{y_n}(y^*)w_{2,\nu}(\cdot,0)\Vert^2
_{L^2(\Bbb R^n)}
+ \vert\widetilde s\vert\Vert \psi_{y_n}(y^*)\widetilde \nabla w_{2,\nu}
(\cdot,0)\Vert_{L^2(\Bbb R^n)}^2 \\
\le C_{22}\sum_{j=1}^{n-1}\varXi_\mu(w_{j,n,\nu})+\epsilon(\delta)\vert
\widetilde s\vert
\Vert (\partial_{y_n} \mbox{\bf w}_{\nu}(\cdot,0),\mbox{\bf  w}
_{\nu}(\cdot,0))\Vert^2
_{L^2(\Bbb R^n)\times H^{1,\widetilde s}(\Bbb R^n)}
+ C_{23}\vert\widetilde s\vert\Vert \mbox{\bf g}
e^{\vert s\vert\varphi}
\Vert_{L^2(\Bbb R^n)}^2\nonumber\\+C_{24}\tau^{4n+6}\Vert
(\partial_{y_n} \mbox{\bf w}(\cdot,0),\mbox{\bf  w}
(\cdot,0))\Vert^2
_{L^2(\Bbb R^n)\times H^{1,\widetilde s}(\Bbb R^n)}\nonumber\\
+C_{25}(\Vert P_{\lambda+2\mu}(y,D,\tilde s,\tau)w_{2,\nu}\Vert^2
_{L^2(\mathcal Q)}
+\tau^2\Vert w_{2,\nu}\Vert^2_{H^{1,\widetilde s}(\mathcal Q)}). \nonumber
\end{eqnarray}
Inequalities (\ref{2.4}) and (\ref{02.4}) imply (\ref{pulemet}).

{\bf Case B.}  Let $\tilde s^*=0.$
Then $\xi_0^*\ne 0$ and therefore $r_{\lambda+2\mu}(y^*,\zeta^*,\tau)
\ne 0.$
By Proposition \ref{gorokx} there exists a constant $C_{16}$ independent of
$s$ and $\tau$ such that estimate (\ref{ix}) holds true.
Let $\delta_1$ be sufficiently small, so that
there exists a constant $C_{26}>0$ such that
\begin{equation}\label{2.8}
\vert\xi_0\vert^2\le C_{26}(\sum_{j=1}^{n-1}\vert\xi_j\vert^2+\widetilde s^2
) \quad \forall \zeta\in\mathcal
O(\zeta^*,\delta_1(y^*)). \end{equation}

Now we need again to estimate $\varXi_\mu(j,n).$ We start from
the term $\frak J_2(\mu, w_{\widehat j,n,\nu})  $. By (\ref{ix1}) we have
\begin{eqnarray}\label{2.12}
\frak J_2(\mu, w_{\widehat j,n,\nu})                  \nonumber\\
= -\frac 12\text{Re}\int_{{\Bbb R^n}}2\vert\widetilde
s\vert(\lambda+2\mu)(y^*)\left(\partial_{y_{\widehat j}} w_{2,\nu}
- \vert\widetilde
s\vert\psi_{y_{\widehat j}}(y^*)w_{2,\nu}\right)
\overline{(\nabla_{\widetilde \xi}p_\mu(y^*,\nabla w_{\widehat j,n,\nu},0),
\widetilde \nabla \psi(y^*))} d\widetilde y                 \nonumber\\
- \frac 12\text{Re}\int_{{\Bbb R^n}}2\vert \widetilde s\vert(\lambda+2\mu)(y^*)
(\vert \widetilde s\vert\psi_{y_n}(y^*)w_{\widehat j,n,\nu}+r_{\widehat j})
\overline{(\nabla_{\widetilde \xi}p_\mu(y^*,\nabla w_{\widehat j,n,\nu},0),
\widetilde \nabla \psi(y^*))} d\widetilde y               \nonumber\\
= -\frac 12\text{Re}\int_{{\Bbb R^n}}2\vert\widetilde
s\vert(\lambda+2\mu)(y^*)\left(\partial_{ y_{\widehat j}} w_{2,\nu}
-\vert
\widetilde s\vert\psi_{y_{\widehat j}}(y^*)w_{2,\nu}\right)
\overline{(\nabla_{\widetilde \xi}p_\mu(y^*,\nabla w_{\widehat j,n,\nu},0),
\widetilde \nabla \psi(y^*))}d\widetilde y\nonumber\\
-\frac 12\text{Re}\int_{{\Bbb R^n}}2\vert s\vert(\lambda+2\mu)(y^*)
r_{\widehat j}\overline{(\nabla_{\widetilde \xi}p_\mu(y^*,\nabla w_{\widehat j,
n,\nu},0), \widetilde \nabla \psi(y^*))} d\widetilde y .
\end{eqnarray}
Integrating by parts we have
\begin{eqnarray}\label{2.12}
I_2(\widehat j)=-\text{Re}\int_{{\Bbb R^n}}\vert\widetilde
s\vert(\lambda+2\mu)(y^*)\left(
\partial_{y_{\widehat j}} w_{2,\nu}-\vert\widetilde s\vert
\psi_{y_{\widehat j}}(y^*) w_{2,\nu}\right)
\overline{(\nabla_{\widetilde \xi}p_\mu(y^*,\nabla w_{\widehat j,n,\nu},0),
\widetilde
\nabla \psi(y^*))} d\widetilde y                      \nonumber\\
=\text{Re}\int_{{\Bbb R^n}}\vert\widetilde
s\vert(\lambda+2\mu)(y^*)\partial_{y_{\widehat j}} w_{\widehat j,n,\nu}
\overline{(\nabla_{\widetilde \xi}p_\mu(y^*,
\nabla w_{2,\nu},0), \widetilde \nabla \psi(y^*))} d\widetilde y \\
+\text{Re}\int_{{\Bbb R^n}}\vert\widetilde
s\vert(\lambda+2\mu)(y^*)\vert\widetilde s\vert\psi_{y_{\widehat j}}(y^*)
w_{2,\nu}\overline{(\nabla_{\widetilde \xi}p_\mu(y^*,\nabla w_{\widehat j,n,
\nu},0), \widetilde \nabla \psi(y^*))}d\widetilde y .\nonumber
\end{eqnarray}

Denote
$$
\mathcal M_{\hat j}=\text{Re}\int_{{\Bbb R^n}}2\vert
s\vert(\lambda+2\mu)(y^*)\vert\widetilde s\vert\psi_{y_{\widehat j}}(y^*)
w_{2,\nu}\overline{(\mu(y^*)\sum_{j=1}^{n-1}\partial_{y_j}
w_{\widehat j,n,\nu}\psi_{y_j}(y^*)-\partial_{y_0}
w_{\widehat j,n,\nu}(\rho\psi_{y_0})(y^*))} d\widetilde y.
$$

The simple computations imply
\begin{eqnarray}\label{ono}
 \vert \mathcal M(\hat j)\vert
\le C_{27}\int_{{\Bbb R^n}}\widetilde s\vert p_\mu(y^*,\widetilde \xi + {i}
\widetilde s\widetilde \nabla \psi,0)\vert \vert \widehat w_{2,\nu}\widehat
w_{\widehat j,n,\nu}
\vert d\widetilde \xi
\\
\le C_{28}\delta_1 \vert\widetilde s\vert  \Vert
(\partial_{y_n} \mbox{\bf w}_{\nu}(\cdot,0),\mbox{\bf  w}
_{\nu}(\cdot,0))\Vert^2
_{L^2(\Bbb R^n)\times H^{1,\widetilde s}(\Bbb R^n)}
 +C_{29} \tau^{2n+6}\Vert
(\partial_{y_n} \mbox{\bf w}(\cdot,0),\mbox{\bf  w}(\cdot,0))\Vert^2
_{L^2(\Bbb R^n)\times H^{1,\widetilde s}(\Bbb R^n)}.\nonumber
\end{eqnarray}
Here in order to obtain the last equality, we used (\ref{2.7}).

From (\ref{ix}) and (\ref{ix1}) on $\{y_n=0\}$, for some function $\mbox{\bf r}$ we obtain
\begin{equation}\label{ono1}
{i}(\lambda +2\mu)(y^*)\alpha^+_{\lambda+2\mu}(\widetilde y,0,
{\widetilde D},\tilde s,\tau)w_{2,\nu}
- \mu(y^*)\sum_{j=1}^{n-1}\left(\partial_{y_j}
w_{j,n,\nu}-\widetilde s\psi_{y_j}(y^*)w_{j,n,\nu}
\right) = \mbox{\bf r},
\end{equation}
where $\mbox{\bf r}$ satisfies estimate (\ref{garmodon}).

Using this equation, we transform $\sum_{j=1}^{n-1} I_2(j)$ as
\begin{eqnarray}\label{moloko}
\sum_{j=1}^{n-1} I_2(j)\\
= \text{Re}\int_{{\Bbb R^n}}2\vert
\widetilde s\vert(\lambda+2\mu)(y^*){i}\frac{(\lambda +2\mu)}
{\mu}(y^*)\alpha^+_{\lambda+2\mu}(\widetilde y,0,{\widetilde D},\tilde s,\tau)
w_{2,\nu}
\overline{(\nabla_{\widetilde \xi}p_\mu(y^*,\widetilde \nabla\psi,0),
\widetilde \nabla w_{2,\nu})} d\widetilde y \nonumber\\
-\text{Re}\int_{{\Bbb R^n}}2\vert
\widetilde s\vert\frac{(\lambda +2\mu)}{\mu}(y^*)\mbox{\bf r}
\overline{(\nabla_{\widetilde \xi}p_\mu(y^*,\widetilde \nabla\psi,0),
\widetilde \nabla w_{2,\nu})} d\widetilde y\nonumber\\
+\text{Re}\int_{{\Bbb R^n}}2\vert
\widetilde s\vert(\lambda+2\mu)(y^*) \sum_{j=1}^{n-1}\vert \widetilde s\vert
\psi_{y_j}(y^*)w_{j,n,\nu}\overline{(\nabla_{\widetilde \xi}p_\mu(y^*,
\widetilde \nabla\psi,0),\widetilde \nabla w_{2,\nu})} d\widetilde y
+\sum_{j=1}^{n-1} \mathcal M_{j}.\nonumber
\end{eqnarray}

Since $r_\mu(y^*,\zeta^*,\tau)=0$, we have
$r_{\lambda+2\mu}(y^*,\zeta^*,\tau)=-\frac{(\lambda+\mu)(y^*)}
{\mu(y^*)(\lambda+2\mu)(y^*)}
(\xi_0^*+{i}\vert \widetilde s^*\vert\psi_{y_0}(y^*))^2$.
This implies that $Re\{(\nabla_{\widetilde \xi}p_\mu(y^*,\widetilde
\nabla\psi,0),\xi^*)\alpha^+_{\lambda+2\mu}(y^*,\zeta^*,\tau)\}=0.$

Therefore by (\ref{moloko}), (\ref{ono}), (\ref{garmodon}) and the G\aa rding
inequality  (\ref{goblin}), we have
 \begin{eqnarray}\label{zavik}
\sum_{j=1}^{n-1} I_2(j)\ge -\epsilon(\delta)\vert\widetilde s\vert
  \Vert
(\partial_{y_n} \mbox{\bf w}_{\nu}(\cdot,0),\mbox{\bf  w}
_{\nu}(\cdot,0))\Vert^2
_{L^2(\Bbb R^n)\times H^{1,\widetilde s}(\Bbb R^n)}\\-C_{30}\tau^{4n+6}\Vert
(\partial_{y_n} \mbox{\bf w}(\cdot,0),\mbox{\bf  w}(\cdot,0))\Vert^2
_{L^2(\Bbb R^n)\times H^{1,\widetilde s}(\Bbb R^n)}-C_{31}
\vert\widetilde s\vert\Vert \mbox{\bf g}
e^{\vert s\vert\varphi}\Vert_{L^2(\Bbb R^n)}^2.\nonumber
\end{eqnarray}
By (\ref{zavik}), (\ref{2.99}), (\ref{2.12}), (\ref{garmodon}) we obtain
\begin{eqnarray}\label{resonans}
\sum_{j=1}^{n-1}\varXi_\mu(w_{j,n,\nu})
\ge  C_{32}\int_{{\Bbb R^n}} (\vert \widetilde s\vert
\sum_{j=1}^{n-1}\left\vert\partial_{y_n} w_{j,n,\nu}
\right\vert^2 +\vert \widetilde s\vert^3
\vert w_{j,n,\nu}\vert^2)(\tilde y,0) d\widetilde y     \nonumber\\
- \epsilon(\delta)\vert\widetilde s\vert \Vert
(\partial_{y_n} \mbox{\bf w}_{\nu}(\cdot,0),\mbox{\bf  w}
_{\nu}(\cdot,0))\Vert^2
_{L^2(\Bbb R^n)\times H^{1,\widetilde s}(\Bbb R^n)}
\\-C_{33}\tau^{4n+6}\Vert
(\partial_{y_n} \mbox{\bf w}(\cdot,0),\mbox{\bf  w}(\cdot,0))\Vert^2
_{L^2(\Bbb R^n)\times H^{1,\widetilde s}(\Bbb R^n)}
-C_{34}\vert\widetilde s\vert\Vert \mbox{\bf g}
e^{\vert s\vert\varphi}\Vert_{L^2(\Bbb R^n)}^2.\nonumber
\end{eqnarray}

By (\ref{ix1}), (\ref{garmodon}) and (\ref{resonans}), we have
\begin{eqnarray}\label{resonans1}
\sum_{j=1}^{n-1}\varXi_\mu(w_{j,n,\nu})\ge  C_{35}\int_{{\Bbb R^n}}
(\vert \widetilde s\vert
\sum_{j=0}^{n-1}\vert
\partial_{y_j} w_{2,\nu}\vert^2
+ \vert\widetilde s\vert\left\vert\partial_{y_n} w_{2,\nu}
-\vert\widetilde s\vert\psi_{y_n}w_{2,\nu}\right\vert^2)(\tilde y,0)
d \widetilde y                                   \nonumber\\
- \epsilon(\delta)\vert\widetilde s\vert  \Vert
(\partial_{y_n} \mbox{\bf w}_{\nu}(\cdot,0),\mbox{\bf  w}
_{\nu}(\cdot,0))\Vert^2
_{L^2(\Bbb R^n)\times H^{1,\widetilde s}(\Bbb R^n)}
- C_{36}\vert\widetilde s\vert\Vert \mbox{\bf g}
e^{\vert s\vert\varphi}
\Vert_{L^2(\Bbb R^n)}^2            \\
\ge  C_{37}\int_{{\Bbb R^n}}\left(\vert \widetilde s\vert
\sum_{j=0}^{n}\vert\partial_{y_j} w_{2,\nu}\vert^2
+ \vert\widetilde s\vert^3\vert w_{2,\nu}\vert^2\right)d \widetilde y
- \epsilon(\delta)\vert\widetilde s\vert  \Vert
(\partial_{y_n} \mbox{\bf w}_{\nu}(\cdot,0),\mbox{\bf  w}
_{\nu}(\cdot,0))\Vert^2
_{L^2(\Bbb R^n)\times H^{1,\widetilde s}(\Bbb R^n)}
                                                \nonumber\\
-C_{38}\tau^{4n+6}\Vert
(\partial_{y_n} \mbox{\bf w}(\cdot,0),\mbox{\bf  w}(\cdot,0))\Vert^2
_{L^2(\Bbb R^n)\times H^{1,\widetilde s}(\Bbb R^n)}
+C_{39}\vert\widetilde s\vert\Vert \mbox{\bf g}
e^{\vert s\vert\varphi}\Vert_{L^2(\Bbb R^n)}^2.\nonumber                              \nonumber
\end{eqnarray}

Now we estimate the tangential derivatives of $w_{j,n,\nu}.$
By (\ref{ono1}) and (\ref{resonans1}), we see that there  exists a function
$p$ such that
\begin{equation}\label{hhh}
\sum_{j=1}^{n-1}\partial_{y_j}
w_{j,n,\nu}(\widetilde y,0)
= p (\tilde y)  \quad \mbox{in}\,\,\Bbb R^n,\nonumber
\end{equation}
\begin{eqnarray}\label{hhhh}
\vert\widetilde s\vert\Vert p\Vert^2_{H^{1,\widetilde s}(\Bbb R^n)}
\le \epsilon(\delta)\vert\widetilde s\vert \Vert
(\partial_{y_n} \mbox{\bf w}_{\nu}(\cdot,0),\mbox{\bf  w}
_{\nu}(\cdot,0))\Vert^2
_{L^2(\Bbb R^n)\times H^{1,\widetilde s}(\Bbb R^n)}+C_{40}\sum_{j=1}^{n-1}
\varXi_\mu(w_{j,n,\nu})
\nonumber\\+C_{41}\tau^{4n+6}\Vert
(\partial_{y_n} \mbox{\bf w}(\cdot,0),\mbox{\bf  w}(\cdot,0))\Vert^2
_{L^2(\Bbb R^n)\times H^{1,\widetilde s}(\Bbb R^n)}
+ C_{42}\vert\widetilde s\vert\Vert \mbox{\bf g}
e^{\vert s\vert\varphi}\Vert_{L^2(\Bbb R^n)}^2.
\end{eqnarray}
Taking the Fourier transform of the first equality  in (\ref{hhh}) we have
\begin{equation}\label{pop1}
\sum_{j=1}^{n-1}i\xi_j
\widehat w_{j,n,\nu}(\widetilde \xi,0)=\widehat p\quad \forall \widetilde
\xi\in \Bbb R^n.
\end{equation}

By (\ref{ix1}) for $1\le k,j\le n-1$, there exist $p_{kj}(\widetilde \xi)$
such that
\begin{eqnarray}\label{pop2}
\xi_k\widehat w_{j,n,\nu}(\widetilde \xi,0)-\xi_j \widehat w_{k,n,\nu}
(\widetilde \xi,0)
= p_{kj},\\
\vert\widetilde s\vert\Vert p_{kj}\Vert^2_{H^{1,\widetilde s}(\Bbb R^n)}
\le \epsilon(\delta)\vert \widetilde s\vert  \Vert
(\partial_{y_n} \mbox{\bf w}_{\nu}(\cdot,0),\mbox{\bf  w}
_{\nu}(\cdot,0))\Vert^2
_{L^2(\Bbb R^n)\times H^{1,\widetilde s}(\Bbb R^n)}\nonumber
\\+C_{43}\tau^{4n+6}\Vert
(\partial_{y_n} \mbox{\bf w}(\cdot,0),\mbox{\bf  w}(\cdot,0))\Vert^2
_{L^2(\Bbb R^n)\times H^{1,\widetilde s}(\Bbb R^n)}
+ C_{44}\vert\widetilde s\vert\Vert \mbox{\bf g}
e^{\vert s\vert\varphi}\Vert_{L^2(\Bbb R^n)}^2.\nonumber\nonumber
\end{eqnarray}
If $\xi_{\widehat j}^*\ne 0$  from (\ref{pop2}), then we have
\begin{equation}
\widehat w_{k,n,\nu}(\widetilde \xi,0)=\frac{\xi_k}{\xi_{\widehat j}}
\widehat w_{\widehat j,n,\nu}(\widetilde \xi,0)
- \frac{p_{k,\widehat j}(\widetilde \xi)}{\xi_{\widehat j}}\quad
                    \forall \zeta\in\mathcal
O(\zeta^*,\delta_1(y^*)).
\end{equation} Substituting this equality into (\ref{pop1}), we obtain
\begin{equation}\label{pop4}
\widehat w_{\widehat j,n,\nu}(\widetilde \xi,0)\sum_{k=1}^{n-1}\xi_k^2
=\xi_{\widehat j}(\hat p(\widetilde \xi)+\sum_{k=1,j\ne \widehat j}^{n-1}
p_{k,\widehat j}(\widetilde \xi))\quad   \forall  \zeta\in\mathcal
O(\zeta^*,\delta_1(y^*)).\nonumber
\end{equation}
Inequalities (\ref{pop1}), (\ref{hhhh}) and (\ref{pop2}) yield
\begin{eqnarray}\label{pop3}
\vert\widetilde s\vert\int_{{\Bbb R^n}}\sum_{j=1}^{n-1}
\vert\partial_{y_j} w_{\widehat j,n,\nu}(\tilde y,0)\vert^2
d \widetilde y
\le \epsilon(\delta)\vert\widetilde s\vert \Vert
(\partial_{y_n} \mbox{\bf w}_{\nu}(\cdot,0),\mbox{\bf  w}
_{\nu}(\cdot,0))\Vert^2
_{L^2(\Bbb R^n)\times H^{1,\widetilde s}(\Bbb R^n)}
+
\\C_{45}\sum_{j=1}^{n-1}\varXi_\mu(w_{j,n,\nu})+C_{46}\tau^{4n+6}\Vert
(\partial_{y_n} \mbox{\bf w}(\cdot,0),\mbox{\bf  w}(\cdot,0))\Vert^2
_{L^2(\Bbb R^n)\times H^{1,\widetilde s}(\Bbb R^n)}
+C_{47}\vert\widetilde s\vert\Vert \mbox{\bf g}
e^{\vert s\vert\varphi}\Vert_{L^2(\Bbb R^n)}^2.\nonumber
\end{eqnarray}
By (\ref{pop1}), (\ref{pop2}) and (\ref{pop3}), we obtain from (\ref{pop4})
\begin{eqnarray}\label{pop3}
\vert\widetilde s\vert\int_{{\Bbb R^n}}\sum_{j,k=1}^{n-1}
\left\vert \partial_{y_j} w_{ k,n,\nu}(\tilde y,0)\right\vert^2d
\widetilde y
\le \epsilon(\delta)\vert\widetilde s\vert
\Vert
(\partial_{y_n} \mbox{\bf w}_{\nu}(\cdot,0),\mbox{\bf  w}
_{\nu}(\cdot,0))\Vert^2
_{L^2(\Bbb R^n)\times H^{1,\widetilde s}(\Bbb R^n)}
\\+C_{48}(\sum_{j=1}^{n-1}\varXi_\mu(w_{j,n,\nu})+\tau^{4n+6}\Vert
(\partial_{y_n} \mbox{\bf w}(\cdot,0),\mbox{\bf  w}(\cdot,0))\Vert^2
_{L^2(\Bbb R^n)\times H^{1,\widetilde s}(\Bbb R^n)}
+\vert\widetilde s\vert\Vert \mbox{\bf g}
e^{\vert s\vert\varphi}\Vert_{L^2(\Bbb R^n)}^2).\nonumber
\end{eqnarray}
Form this inequality and (\ref{2.8}), we have
\begin{eqnarray}\label{pop4}\vert\widetilde
s\vert\int_{{\Bbb R^n}}\sum_{k=1}^{n-1}\vert \nabla' w_{ k,n,\nu}(\tilde y,0)\vert^2
d \widetilde y
\le \epsilon(\delta)\vert\widetilde s\vert
\Vert
(\partial_{y_n} \mbox{\bf w}_{\nu}(\cdot,0),\mbox{\bf  w}
_{\nu}(\cdot,0))\Vert^2
_{L^2(\Bbb R^n)\times H^{1,\widetilde s}(\Bbb R^n)}
+
\\C_{49}(\sum_{j=1}^{n-1}\varXi_\mu(w_{j,n,\nu})+\tau^{4n+6}\Vert
(\partial_{y_n} \mbox{\bf w}(\cdot,0),\mbox{\bf  w}(\cdot,0))\Vert^2
_{L^2(\Bbb R^n)\times H^{1,\widetilde s}(\Bbb R^n)}
+\vert\widetilde s\vert\Vert \mbox{\bf g}
e^{\vert s\vert\varphi}\Vert_{L^2(\Bbb R^n)}^2).\nonumber
\end{eqnarray}
Inequalities (\ref{pop4}), (\ref{resonans}) and (\ref{resonans1})
imply (\ref{pulemet}).

\bigskip
\section{ Case $r_{\lambda+2\mu}(y^*,\zeta^*,\tau)=0.$}\label{section4}
\bigskip
Let $(y^*,\zeta^*)$ be a point on ${\Bbb R}^{n+1}\times \Bbb S^{n}$ such
that $r_{\lambda+2\mu}(y^*,\zeta^*,\tau)=0$ and $\text{supp}\,
\chi_\nu\subset \mathcal O(y^*,\delta_1(y^*)).$  We note that if
$r_\mu(y^*,\zeta^*,\tau)=0$ then $\widetilde s^*=0,
\widetilde\nabla\psi(y^*)=0$.
This case was treated in the previous section. Therefore we may
assume $r_\mu(y^*,\zeta^*,\tau)\ne 0$. By (\ref{3.35})-(\ref{3.37})
there exists
$\delta_0>0$ and $C_1>0$ such that for all $\delta_1\in
(0,\delta_0)$ we have
\begin{equation}\label{3.1}
\vert \xi_0\vert^2
\le C_1(\sum_{j=1}^{n-1}\xi_j^2+\widetilde s^2)
\quad \forall \zeta\in \mathcal O(y^*,\delta_1(y^*)).\end{equation}

By (\ref{3.35})-(\ref{3.37}) there exists $C_2>0$ such that
\begin{eqnarray}\label{2.777} \vert p_{\lambda+2\mu}(y^*,\widetilde \xi,0)
- \widetilde
s^2 p_{\lambda+2\mu}(y^*,\widetilde \nabla \psi(y^*),0)\vert+\vert\widetilde
s(\nabla
_{\widetilde \xi}p_{\lambda+2\mu}(y^*,\widetilde \xi,0 ),
\widetilde \nabla \psi(y^*))
\vert\nonumber\\ \le \delta_1C_2\vert(\widetilde \xi,\widetilde s)\vert^2,
\quad \forall\zeta\in\mathcal
O(\zeta^*,\delta_1(y^*)).
\end{eqnarray}

By (\ref{klop}) there exists $C_{3}>0$ independent of $s,\tau$ such that
\begin{equation}\label{3.19}
\varXi_{\lambda+2\mu}(w_{2,\nu})+C_{3}({\vert \widetilde s\vert}\tau\Vert
w_{2,\nu}\Vert^2_{H^1(\mathcal Q)}+ {\vert \widetilde s\vert}^3\tau\Vert
w_{2,\nu}\Vert^2_{L^2(\mathcal Q)})
\end{equation}
\begin{eqnarray*}
\le &&C_{4}(\Vert P_{\lambda+2\mu}(y,D,\tilde s,\tau) w_{2,\nu} \Vert^2
_{L^2(\mathcal Q)}
+ \tau^{4n+6}\Vert\text{\bf w}\Vert^2_{H^{1,\widetilde s}(\mathcal Q)}) \\
+ &&\epsilon {\vert \widetilde s\vert} \Vert
(\partial_{y_n} \mbox{\bf w}_{\nu}(\cdot,0),\mbox{\bf  w}
_{\nu}(\cdot,0))\Vert^2
_{L^2(\Bbb R^n)\times H^{1,\widetilde s}(\Bbb R^n)},
\end{eqnarray*}
where $\epsilon(\delta)\rightarrow 0$ as $\delta\rightarrow +0.$

We consider several cases.

{Case \bf A.} Let $  \widetilde s^*=0.
$
Since $\widetilde s^*=0$, by decreasing the parameter $\delta_1$ we can
assume that for some constant $C_5>0$
\begin{equation}\label{3.2}
\vert \xi_0\vert^2+\widetilde s^2\le
C_5\sum_{j=1}^{n-1}\xi_j^2\quad \forall \zeta\in \mathcal O(y^*,\delta_1(y^*)).
\end{equation}
By (\ref{2.777})
\begin{equation}\label{i.xxx}
p_\mu(y,\widetilde \xi,0)-{\vert \widetilde s\vert}^2p_\mu(y,\widetilde
\nabla \psi,0)\ge-\epsilon(\delta,\delta_1) \vert (\tilde \xi,\tilde s)\vert^2 \quad (y,\zeta)
\in \Upsilon_\ell\times\mathcal O(y^*,\delta_1(y^*)).
\end{equation}
If $\lim_{\zeta\rightarrow \zeta^*}Im\, r_\mu(y^*,\zeta,\tau)/\vert\widetilde
s\vert\ne 0$, then we set $\text{\bf I}_\mu=sign\,\lim_{\zeta\rightarrow
\zeta^*}Im\,
r_\mu(y^*,\zeta,\tau)/\vert \widetilde s\vert.$
Note that $Re\, r_\mu(y^*,\zeta^*,\tau)>0$
by $\widetilde s^*=0$. For all $(y,\zeta)\in
\Upsilon_\ell\times\mathcal O(y^*,\delta_1(y^*))$ we have
\begin{equation}\label{3.9}
\Gamma_\mu^+(y^*,\zeta^*,\tau)=\text{\bf I}_\mu\,\root\of{Re\,
r_\mu(y^*,\zeta^*,\tau)}\quad \mbox{if} \lim_{\zeta\rightarrow \zeta^*}
Im\, r_\mu(y^*,\zeta,\tau)/\vert\widetilde
s\vert\ne  0
\end{equation}
and
\begin{equation}\label{3.999}
\Gamma_\mu^+(y^*,\zeta^*,\tau)=0\quad \mbox{if} \lim_{\zeta\rightarrow \zeta^*}
Im\, r_\mu(y^*,\zeta,\tau)/\vert\widetilde
s\vert= 0.
\end{equation}
By (\ref{3.9}) and (\ref{3.999}), there exists $\epsilon(\delta,\delta_1)>0$
such that
\begin{equation}\label{3.11}
-(\nabla_{\widetilde \xi}p_\mu(y^*,\widetilde\xi,0), \widetilde \nabla
\psi(y^*))\Gamma_
\mu^+( y,\zeta,\tau)\ge \epsilon(\delta,\delta_1) \vert\zeta\vert^2
\quad
\forall
(\widetilde y,\zeta)\in  \Upsilon_\ell\times\mathcal O(y^*,\delta_1(y^*)).
\end{equation}

Since $r_\mu(y^*,\zeta^*,\tau)\ne 0$ by  Proposition \ref{gorokx}, we have
\begin{equation}\label{ixx}
\Vert \root\of{{\vert \widetilde s\vert}}(\frac 1i\partial_{y_n}-\Gamma^+_{\mu}
(y,{\widetilde D},\tilde s,\tau))w_{k,j,\nu}\vert_{y_n=0}\Vert
_{L^2({{\Bbb R^n}})}\le
C_6(\Vert P_{\mu}(y,D,\tilde s,\tau)w_{k,j,\nu}\Vert_{L^2(\mathcal Q)}
+\tau^{2n+3}\Vert
\text{\bf w}\Vert_{H^{1,\widetilde s}(\mathcal Q)}).
\end{equation}
Let us consider formula
(\ref{2.1}) from the previous section. By  (\ref{3.11}), Lemmata \ref{Fops3}
and \ref{Fops5}, we have
\begin{eqnarray}\label{3.12}
\frak J_2(\mu, w_{k,n,\nu})=-\text{Re}\int_{{\Bbb R^n}}2\vert\widetilde
s\vert\mu(y^*)\partial_{y_n} w_{k,n,\nu}
\overline{(\nabla_{\widetilde \xi}p_\mu(y^*,\widetilde\nabla w_{k,n,\nu},0),
\widetilde \nabla \psi(y^*))}d \widetilde y         \nonumber\\
=-\text{Re}\int_{{\Bbb R^n}}2\vert\widetilde
s\vert\mu(y^*) {i}  \Gamma^+_{\mu}(\widetilde y,0,{\widetilde D},\tilde s,\tau)
w_{k,n,\nu}\overline{(\nabla_{\widetilde \xi}p_\mu(y^*,\widetilde
\nabla w_{k,n,\nu},0), \widetilde \nabla \psi(y^*))}d \widetilde y                          \nonumber\\
+ \text{Re}\int_{{\Bbb R^n}}\vert\widetilde s\vert\mu(y^*)
{i} V^+_\mu(k,n)(\cdot,0)\overline{(\nabla_{\widetilde \xi}p_\mu(y^*,
\widetilde\nabla w_{k,n,\nu},0), \widetilde \nabla \psi(y^*))}d \widetilde y
                                                      \nonumber\\
\ge -C_7(\epsilon(\delta,\delta_1){\vert \widetilde s\vert+\tau^{4n+6})}
\Vert (
\partial_{y_n} w_{k,n,\nu}(\cdot,0),
w_{k,n,\nu}(\cdot,0))\Vert^2_{L^2(\Bbb R^n)\times
H^{1,\widetilde s}(\Bbb R^n)}                               \nonumber\\
+ \text{Re}\int_{{\Bbb R^n}}2\vert\widetilde
s\vert\mu(y^*) {i} V^+_\mu(k,n)(\cdot,0)\overline{
(\nabla_{\widetilde \xi}p_\mu(y^*,\widetilde \nabla w_{k,n,\nu},0), \widetilde
\nabla \psi(y^*))}d \widetilde y.
\end{eqnarray}
Therefore from (\ref{3.12}) and (\ref{ixx}), we obtain
\begin{eqnarray}\label{3.13}
\frak J_2(\mu, w_{ k,n,\nu})\ge
-C_8(\epsilon(\delta,\delta_1){\vert \widetilde s\vert}+\tau^{4n+6})\Vert (
\partial_{y_n} w_{k,n,\nu}(\cdot,0),
w_{k,n,\nu}(\cdot,0))\Vert^2_{L^2(\Bbb R^n)\times
H^{1,\widetilde s}
(\Bbb R^n)} \nonumber\\
-C_{9}(\delta,\delta_1)(\Vert P_{\mu}(y,D,\tilde s,\tau)
w_{k,n,\nu}\Vert^2 _{L^2(\mathcal Q)}+\Vert \text{\bf
w}\Vert^2_{H^{1,\widetilde s}(\mathcal Q)}).
\end{eqnarray}
Inequalities (\ref{3.13}) and (\ref{i.xxx}) imply
\begin{eqnarray}\label{3.14}
\varXi_{\mu}(w_{k,n,\nu})
\ge\int_{\Bbb R^n}\left({\vert \widetilde s\vert}\vert\partial_{y_n}
w_{k,n,\nu}\vert^2
+ {\vert \widetilde s\vert}^3\vert w_{k,n,\nu}\vert^2\right)(\tilde y,0)
d\widetilde y
                                    \nonumber\\
- C_{10}\epsilon(\delta,\delta_1)\vert s\vert\Vert
(\partial_{y_n} \mbox{\bf w}_{\nu}(\cdot,0),\mbox{\bf  w}
_{\nu}(\cdot,0))\Vert^2
_{L^2(\Bbb R^n)\times H^{1,\widetilde s}(\Bbb R^n)}\nonumber\\
-C_{11}\tau^{4n+6}\Vert
(\partial_{y_n} \mbox{\bf w}(\cdot,0),\mbox{\bf  w}(\cdot,0))\Vert^2
_{L^2(\Bbb R^n)\times H^{1,\widetilde s}(\Bbb R^n)} \nonumber\\
-C_{12}(\delta,\delta_1)(\Vert P_{\mu}(y,D,\tilde s,\tau)
w_{k,n,\nu}\Vert^2 _{L^2(\mathcal Q)}+\tau^{4n+6}\Vert \text{\bf w}\Vert^2
_{H^{1,\widetilde s}(\mathcal Q)}).
\end{eqnarray}
By (\ref{poko1}) for any $k\in\{1,\dots, n-1\}$, we obtain
\begin{equation}\label{barmalei}
\mu(y^*)\left(\partial_{y_n} w_{ k,n,\nu}-{\vert \widetilde s
\vert} \psi_{y_n}(y^*) w_{k,n,\nu}\right)(\cdot,0)
\end{equation}
$$
= (\lambda+\mu)(y^*)\left(
\partial_{y_k} w_{2,\nu}-{\vert \widetilde s\vert}
\psi_{y_{k}}(y^*) w_{2,\nu}\right)(\cdot,0)
+ \mbox{\bf r}\quad \mbox{in}\,\,\Bbb R^n,
$$
where the function $\mbox{\bf r}$ satisfies estimate (\ref{garmodon}).
By (\ref{barmalei}), (\ref{3.14}) and (\ref{3.1}), we derive
\begin{eqnarray}\label{3.14!}
\sum_{j=1}^{n-1}\varXi_{\mu}(w_{j,n,\nu})
\ge \int_{\Bbb R^n}({\vert \widetilde s\vert}\vert
\partial_{y_n} w_{j,n,\nu}\vert^2
+ {\vert \widetilde s\vert}^3\vert\mbox{\bf w}_{1,\nu}\vert^2+{\vert
\widetilde s\vert}
\vert \widetilde\nabla w_{2,\nu}\vert^2)(\tilde y,0)d\widetilde y
                                                     \nonumber\\
- C_{13}(\epsilon(\delta,\delta_1){\vert \widetilde s\vert}+1)\left\Vert
\left(
\partial_{y_n} \mbox{\bf  w}_{\nu}(\cdot,0),
\mbox{\bf  w}_{\nu}(\cdot,0)\right)\right\Vert^2_{L^2(\Bbb R^n)\times
H^{1,\widetilde s}(\Bbb R^n)}                      \nonumber\\
-C_{14}\tau^{4n+6}\Vert
(\partial_{y_n} \mbox{\bf w}(\cdot,0),\mbox{\bf  w}(\cdot,0))\Vert^2
_{L^2(\Bbb R^n)\times H^{1,\widetilde s}(\Bbb R^n)}\nonumber\\
- C_{15}(\delta,\delta_1)(\Vert P_{\mu}(y,D,\tilde s,\tau)
w_{j,n,\nu}\Vert^2 _{L^2(\mathcal Q)}+\tau^{4n+6}\Vert \text{\bf w}\Vert^2
_{H^{1,\widetilde s}(\mathcal Q)}).
\end{eqnarray}

By (\ref{ixx}) and Lemma 8.6, (\ref{goblin}),
there exists a constant $C_{16}>0$ such that
\begin{eqnarray}
\sum_{k,j=1,k<j}^n\varXi_\mu (w_{k,j,\nu})\ge
C_{16}{\vert \widetilde s\vert}\sum_{k,j=1,k<j}^n\Vert (
\partial_{y_n} w_{k,j,\nu}(\cdot,0),w_{k,j,\nu}(\cdot,0))\Vert^2
_{L^2(\Bbb R^n)\times H^{1,\widetilde s}(\Bbb R^n)}       \nonumber\\
-C_{17}\tau^{4n+6}\Vert
(\partial_{y_n} \mbox{\bf w}(\cdot,0),\mbox{\bf  w}(\cdot,0))\Vert^2
_{L^2(\Bbb R^n)\times H^{1,\widetilde s}(\Bbb R^n)}\nonumber\\
-C_{18}(\delta,\delta_1)(\Vert P_{\mu}(y,D,\tilde s,\tau)\mbox{\bf
w}_{1,\nu}\Vert^2 _{L^2(\mathcal Q)}+\tau^{4n+6}\Vert \text{\bf
w}\Vert^2_{H^{1,\widetilde s}(\mathcal Q)}).\nonumber
\end{eqnarray}
This inequality and  (\ref{ix1}) imply that
\begin{eqnarray}\label{3.17}
\sum_{k,j=1,k<j}^n\varXi_\mu(w_{k,j,\nu})
\ge C_{19}{\vert \widetilde s\vert}
\Vert
(\partial_{y_n} \mbox{\bf w}_{\nu}(\cdot,0),\mbox{\bf  w}
_{\nu}(\cdot,0))\Vert^2
_{L^2(\Bbb R^n)\times H^{1,\widetilde s}(\Bbb R^n)}         \nonumber\\
-C_{20}\tau^{4n+6}\Vert
(\partial_{y_n} \mbox{\bf w}(\cdot,0),\mbox{\bf  w}(\cdot,0))\Vert^2
_{L^2(\Bbb R^n)\times H^{1,\widetilde s}(\Bbb R^n)}\\
- C_{21}(\delta,\delta_1)(\Vert P_{\mu}(y,D,\tilde s,\tau)
\mbox{\bf w}_{1,\nu}\Vert^2 _{L^2(\mathcal Q)}
+ {\vert \widetilde s\vert}\Vert
\mbox{\bf g}e^{\vert s\vert\varphi}\Vert^2_{L^2({\Bbb R^n})}
+ \tau^{4n+6}\Vert \text{\bf w}\Vert^2_{H^{1,\widetilde s}(\mathcal Q)})
\nonumber
\end{eqnarray}
with some positive constant $C_{19}.$
From (\ref{3.17}) and (\ref{klop}), we obtain (\ref{pulemet}).

{\bf Case B.} Let $\widetilde s^*\ne 0.$ If $\delta_1>0$ is
small enough, then there exists a constant $C_{22}>0$ such that for all $(\tilde\xi,\tilde s)\in \mathcal O(y^*,\delta_1(y^*))$
\begin{equation}\label{3.18}
\vert
\rho(y^*)\xi_0 \psi_{y_0}(y^*)-(\lambda+2\mu)(y^*)\sum_{j=1}^{n-1}
\xi_j\psi_{y_j}(y^*)\vert^2
\le \delta_1^2 C_{22}(\sum_{j=1}^{n-1}\vert
\xi_j\vert^2+\widetilde s^2).
\end{equation}

By (\ref{01}), (\ref{02}) and (\ref{3.18}), we have
\begin{eqnarray}\label{3.21}
\vert \frak J_2(\lambda+2\mu,w_{2,\nu})+\frak J_3(\lambda+2\mu,w_{2,\nu})\vert
\le C_{23}\delta_1{\vert \widetilde s\vert} \Vert (\partial_{y_n}
w_{2,\nu}(\cdot,0),w_{2,\nu}(\cdot,0))\Vert^2
_{L^2(\Bbb R^n)\times H^{1,\widetilde s}(\Bbb R^n)}\nonumber\\
+C_{24}\tau^{4n+6}\Vert (\partial_{y_n}
w_{2}(\cdot,0),w_{2}(\cdot,0))\Vert^2
_{L^2(\Bbb R^n)\times H^{1,\widetilde s}(\Bbb R^n)}.
\end{eqnarray}

By (\ref{3.21}), there exists a constant $C_{25}>0$ such that
\begin{eqnarray}\label{3.22}
\varXi_{\lambda+2\mu}(w_{2,\nu})\ge C_{25}\int_{{\Bbb R^n}} \left(
\vert\widetilde s\vert(\lambda+2\mu)^2(y^*)\psi_{y_n}(y^*)
\vert\partial_{y_n} w_{2,\nu}\vert^2
+ {\vert \widetilde s\vert}^3(\lambda+2\mu)^2(y^*)\psi^3_{y_n}(y^*) \vert
w_{2,\nu}\vert^2\right)d \widetilde y        \nonumber\\
-\epsilon{\vert \widetilde s\vert}
\Vert (\partial_{y_n}
w_{2,\nu}(\cdot,0),w_{2,\nu}(\cdot,0))\Vert^2
_{L^2(\Bbb R^n)\times H^{1,\widetilde s}(\Bbb R^n)}
-C_{26}\tau^{4n+6}\Vert (\partial_{y_n}
w_{2}(\cdot,0),w_{2}(\cdot,0))\Vert^2
_{L^2(\Bbb R^n)\times H^{1,\widetilde s}(\Bbb R^n)}.\nonumber
\end{eqnarray}
Since $\widetilde s^*\ne 0$, we have
\begin{eqnarray}\label{3.22}
\varXi_{\lambda+2\mu}(w_{2,\nu})\ge C_{27}\int_{{\Bbb R^n}} (
{\vert \widetilde s\vert}(\lambda+2\mu)^2(y^*) \vert\nabla
w_{2,\nu}\vert^2+{\vert \widetilde s\vert}^3
(\lambda+2\mu)^2(y^*)\psi^3_{y_n}(y^*)
\vert w_{2,\nu}\vert^2)(\cdot,0)d \widetilde y\\
- \epsilon{\vert \widetilde s\vert}\Vert (\partial_{y_n}
w_{2,\nu}(\cdot,0),w_{2,\nu}(\cdot,0))\Vert^2
_{L^2(\Bbb R^n)\times H^{1,\widetilde s}(\Bbb R^n)}-C_{28}\tau^{4n+6}\Vert
(\partial_{y_n}
w_{2}(\cdot,0),w_{2}(\cdot,0))\Vert^2
_{L^2(\Bbb R^n)\times H^{1,\widetilde s}(\Bbb R^n)}.\nonumber
\end{eqnarray}

By (\ref{barmalei})
\begin{eqnarray}\label{3.22y}
\varXi_{\lambda+2\mu}(w_{2,\nu})\ge C_{29}\int_{{\Bbb R^n}} (
{\vert \widetilde s\vert}(\lambda+2\mu)^2(y^*) \vert\nabla
w_{2,\nu}\vert^2+{\vert \widetilde s\vert}^3
(\lambda+2\mu)^2(y^*)\psi^3_{y_n}(y^*)
\vert w_{2,\nu}\vert^2)d \widetilde y\nonumber\\
+\vert \tilde s\vert\Vert \left(\partial_{y_n} w_{ k,n,\nu}-{\vert \widetilde s
\vert} \psi_{y_n}(y^*) w_{k,n,\nu}\right)(\cdot,0)\Vert^2_{L^2(\Bbb R^n)}
- \epsilon{\vert \widetilde s\vert}\Vert (\partial_{y_n}
w_{2,\nu}(\cdot,0),w_{2,\nu}(\cdot,0))\Vert^2
_{L^2(\Bbb R^n)\times H^{1,\widetilde s}(\Bbb R^n)}\nonumber\\-C_{30}(\tau^{4n+6}\Vert
(\partial_{y_n}
w_{2}(\cdot,0),w_{2}(\cdot,0))\Vert^2
_{L^2(\Bbb R^n)\times H^{1,\widetilde s}(\Bbb R^n)}+\vert \tilde s\vert\Vert \mbox{\bf r}\Vert^2_{L^2(\Bbb R^n)}).
\end{eqnarray}

From (\ref{3.22y}), (\ref{garmodon}) inequality (\ref{3.4}) for $V_\mu^+(i,j)(\cdot,0)$, we obtain the
estimate
\begin{eqnarray}\label{3.25}
{\vert \widetilde s\vert}\Vert \alpha_\mu(\tilde y,0, \tilde D,\tilde s,\tau) w_{k,j,\nu}(\cdot,0)\Vert^2
_{L^2({\Bbb R^n})}
\le C_{31}\biggl(\varXi_{\lambda+2\mu}(w_{2,\nu})\nonumber\\ +\epsilon{\vert \widetilde s\vert}\Vert (\partial_{y_n}
w_{2,\nu}(\cdot,0),w_{2,\nu}(\cdot,0))\Vert^2
_{L^2(\Bbb R^n)\times H^{1,\widetilde s}(\Bbb R^n)}
+\Vert P_{\mu}(y,D,\tilde s,\tau)
\mbox{\bf w}_{1,\nu} \Vert^2_{L^2(\mathcal Q)}+\tau^{4n+6}\Vert
\text{\bf w}\Vert^2_{H^{1,\widetilde s}(\mathcal Q)}\nonumber\\+\tau^{4n+6}\Vert
(\partial_{y_n}
\mbox{\bf w}(\cdot,0),\mbox{\bf w}(\cdot,0))\Vert^2
_{L^2(\Bbb R^n)\times H^{1,\widetilde s}(\Bbb R^n)}+
\Vert\tilde s\vert\Vert \mbox{\bf g} e^{\vert s\vert\varphi}\Vert^2_{L^2(\Bbb R^n)}\biggr).
\end{eqnarray}
Since by (\ref{napoleon}) $\nabla' \tilde \ell(0)=0$ then $\vert\alpha_\mu(y^*, \zeta^*,\tau)\vert=\vert\root \of{r_\mu(y^*,\zeta^*,\tau)}\vert\ne 0.$ By Lemma 8.2
$$
\Vert \alpha_\mu(\tilde y,0, \tilde D,\tilde s,\tau) w_{k,j,\nu}(\cdot,0)\Vert^2
_{L^2({\Bbb R^n})}=\int_{\Bbb R^n} A(\tilde y,\tilde D,\tilde s,\tau)w_{k,j,\nu}(\cdot,0)\overline w_{k,j,\nu}(\cdot,0)d\tilde y
$$
$$+\int_{\Bbb R^n}R \alpha_\mu(\tilde y,0,\tilde D,\tilde s,\tau)w_{k,j,\nu}(\cdot,0)\overline w_{k,j,\nu}(\cdot,0)d\tilde y.
$$
Here $A(\tilde y,\tilde D,\tilde s,\tau)$ is the pseudodifferential operator with the symbol $\vert\alpha_\mu(\tilde y,0, \zeta,\tau)\vert^2$ and $R\in \mathcal L( H_0^{1,\tilde s}(\Upsilon_\ell)) ; L^2(\Upsilon_\ell)))$ with the norm $\Vert R\Vert\le C_{32}\pi_{C^1(\Upsilon_\ell)}(\alpha_\mu)\le C_{33}\tau^2.$
Therefore
$$
\Vert \alpha_\mu(\tilde y,0, \tilde D,\tilde s,\tau) w_{k,j,\nu}(\cdot,0)\Vert^2
_{L^2({\Bbb R^n})}=\int_{\Bbb R^n} A(\tilde y,\tilde D,\tilde s,\tau)w_{k,j,\nu}(\cdot,0)\overline w_{k,j,\nu}(\cdot,0)d\tilde y-C_{34}\tau^2\Vert w_{k,j,\nu}(\cdot,0)\Vert_{L^2(\mathcal Q)}^2.
$$
Applying  the G\aa rding inequality (\ref{goblin}) we obtain
\begin{equation}\label{granata}
\Vert \alpha_\mu(\tilde y,0, \tilde D,\tilde s,\tau) w_{k,j,\nu}(\cdot,0)\Vert^2
_{L^2({\Bbb R^n})}\ge C_{35}\Vert w_{k,j,\nu}(\cdot,0)\Vert^2_{H^{1,\tilde s}(\mathcal Q)}-C_{36}\tau^{4n+6}\Vert w_{k,j,\nu}(\cdot,0)\Vert_{L^2(\mathcal Q)}^2.
\end{equation}
Inequalities (\ref{ixx}), (\ref{3.22}), (\ref{granata}) and (\ref{3.25})  imply
\begin{eqnarray}\label{3.26}
\varXi_{\lambda+2\mu}(w_{2,\nu})\ge C_{37}{\vert \widetilde s\vert}
\Vert
(\partial_{y_n} \mbox{\bf w}_{\nu}(\cdot,0),\mbox{\bf  w}
_{\nu}(\cdot,0))\Vert^2
_{L^2(\Bbb R^n)\times H^{1,\widetilde s}(\Bbb R^n)}\nonumber\\
- C_{38}(\Vert
 P_{\mu}(y,D,\tilde s,\tau) \mbox{\bf w}_{1,\nu}\Vert^2
_{L^2(\mathcal Q)}+{\vert \widetilde s\vert}\Vert  \mbox{\bf g}
e^{\vert s\vert\varphi}\Vert^2_{L^2({\Bbb R^n})}+\tau^{4n+6}\Vert \text{\bf w}
\Vert^2_{H^{1,\widetilde s}(\mathcal Q)})\nonumber\\
-C_{39}\tau^{4n+6}\Vert (\partial_{y_n}
w_{2}(\cdot,0),w_{2}(\cdot,0))\Vert^2
_{L^2(\Bbb R^n)\times H^{1,\widetilde s}(\Bbb R^n)},
\end{eqnarray}
where $C_{37}>0.$  From (\ref{3.19}) and (\ref{3.26}), we obtain
(\ref{pulemet}).
$\blacksquare$

\bigskip
\section{ Case $r_\mu(y^*,\zeta^*,\tau)\ne 0$ and
$r_{\lambda+2\mu}(y^*,\zeta^*,\tau)\ne 0.$}\label{section5}
\bigskip

In this section we consider the conic neighborhood $ \mathcal O(y^*,
\delta_1(y^*))$ of the point $(y^*,\zeta^*)$ such that
\begin{equation}\label{4.1}
\vert r_\mu(y^*,\zeta^*,\tau)\vert\ne 0\quad\text{and} \,\,
\vert r_{\lambda+2\mu}
(y^*,\zeta^*,\tau)\vert\ne 0. \end{equation}

In that case, thanks to (\ref{4.1}) and Proposition \ref{gorokx1},
factorization
(\ref{min}) holds true for $\beta=\mu$ and $\beta=\lambda+2\mu.$
Then Proposition \ref{gorokx} yields the a priori estimate
\begin{eqnarray}\label{4.4}
(1+\vert \tilde s\vert)(\sum_{k,j=1, k<j}^n\Vert
 V^+_{\mu}(k,j)(\cdot,0)\Vert^2_{L^2({\Bbb R^n})}+
\Vert
V^+_{\lambda+2\mu}(\cdot,0)\Vert^2_{L^2({\Bbb R^n})} )   \\
\le C_1(\Vert P_{\lambda+2\mu}(y,D,\tilde s,\tau) w_{2,\nu}\Vert^2
_{L^2(\mathcal
Q)}+\Vert P_\mu(y,D,\tilde s,\tau) \mbox{\bf w}_{1,\nu}
\Vert^2_{L^2(\mathcal Q)} +\tau^{4n+6}\Vert
\text{\bf w}\Vert_{H^{1,\widetilde s}(\mathcal Q)}^2).   \nonumber
\end{eqnarray}

Using (\ref{poko1}), we rewrite (\ref{ix1}) as
\begin{equation}\label{4.5}
\frac{\lambda+2\mu}{\mu}(y^*)\left(\partial_{y_j}  w_{2,\nu}
-\vert\widetilde  s\vert\psi_{y_j}
w_{2,\nu}\right)
- {i}\alpha^+_\mu(\widetilde y,0,\widetilde D,\tilde s,\tau)w_{j,n,\nu}
=V^+_{\mu}(i,n)(\cdot,0)-r_{j,n,\nu},
\end{equation}
where $i\in\{1,\dots, n-1\}$ and
\begin{equation}\label{4.6}
\sum_{k=1}^{n-1}\frac{\mu}{\lambda+2\mu}(y^*)\left(-\partial_{y_k}
w_{k,n,\nu}+{\vert \widetilde s\vert}\psi_{y_k}
w_{k,n,\nu}\right)
\end{equation}
$$
- {i}\alpha^+_{\lambda+2\mu}(\widetilde y,0,\widetilde D,\tilde s,\tau)
w_{2,\nu}=V^+_{\lambda+2\mu}(\cdot,0)-r_{2,\nu},
$$
where the function $\mbox{\bf r}=(r_{1,n,\nu},\dots, r_{n-1,n,\nu},r_{2,\nu})$
satisfies estimate (\ref{garmodon}).
 Let
$\text{\bf B}(\widetilde y,{\widetilde D},\tilde s,\tau)$ be the matrix
pseudodifferential operator with the symbol
\begin{equation}\label{logo}
\mbox{\bf B}(\tilde y,\zeta,\tau)=
\left( \begin{matrix}
-i\alpha^+_\mu(\tilde y,0,\zeta,\tau)&0 &\dots&\frac{\lambda+2\mu}{\mu}
(i\xi_1-\vert\tilde
s\vert\psi_{y_1}) \\
0&\dots&\dots &\dots\\
0&
-i\alpha^+_\mu(\tilde y,0,\zeta,\tau)&\dots &\frac{\lambda+2\mu}{\mu}
(i\xi_i-\vert\tilde
s\vert\psi_{y_i}) \\
\dots&\dots&\dots &\dots\\
 \frac{\mu}{\lambda+2\mu}(-i\xi_1+\vert\tilde s\vert\psi_{y_1}) &\dots
&\frac{\mu}{\lambda+2\mu}(-i\xi_{n-1}+\vert\tilde s\vert\psi_{y_{n-1}})
 &-i\alpha^+_{\lambda+2\mu}(\tilde y,0,\zeta,\tau)
\end{matrix}\right).
\end{equation}
%

We have
\begin{proposition}
The following formula is true:
\begin{eqnarray}\label{pardon}
det\,\text{\bf B}(y^*,\zeta,\tau)=(-{i})^n (\alpha^+_\mu(y^*,\zeta,\tau))^{n-1}
\alpha^+_{\lambda+2\mu}(y^*,\zeta,\tau)                           \nonumber\\
+ (-1)^{n-1}(-{i})^{n-2} (\alpha^+_\mu(y^*,\zeta,\tau))^{n-2}\sum_{j=1}
^{n-1}
(-{i}\xi_j+\vert s\vert\psi_{y_j}(y^*))^2.
\end{eqnarray}
\end{proposition}
{\bf Proof.} By $\text{\bf B}_{n}$ we denote the matrix determined by
(\ref{logo}) of the size $n\times n$ and $\text{\bf B}_{i,j,n}$ be the minor
obtained
from the matrix  $\text{\bf B}_{n}$ by crossing out the $i$-th row and
the $j$-th column.   Our proof is based on the  induction method.
Except the formula (\ref{pardon}), we claim
\begin{equation}\label{paka}
\vert \text{\bf B}_{1,n-1,n}(y^*,\zeta,\tau)\vert
= (-{i})^{n} (\alpha^+_\mu(y^*,\zeta,\tau))^{n-2}(-{i}\xi^*_1
+\vert \widetilde s\vert\psi_{y_1}(y^*))^2.
\end{equation}
For $n=2,3$, we can easily verify the formulae by direct computations.
Suppose that (\ref{pardon}) and (\ref{paka}) are true for $n-1$. Then
\begin{eqnarray}\label{pardon3}
det\,\text{\bf B}_{n-1}(y^*,\zeta ,\tau)=(-{i})^{n-1}
(\alpha^+_\mu(y^*,\zeta,\tau))^{n-2}\alpha^+_{\lambda+2\mu}(y^*,\zeta,\tau)
                                         \nonumber\\
+ (-1)^{n-1}(-{i})^{n-3} (\alpha^+_\mu(y^*,\zeta,\tau))^{n-3}
\sum_{j=1}^{n-2}(-{i}\xi_j+\vert \widetilde s\vert\psi_{y_j}(y^*))^2
\end{eqnarray}
and
\begin{equation}\label{pardon1}
\vert \text{\bf B}_{1,n-1,n}(y^*,\zeta,\tau)\vert
= (-{i})^{n-1}\frac{\mu}{\lambda+2\mu}(y^*)(\alpha^+_
\mu(y^*,\zeta,\tau))^{n-2}
(-{i}\xi^*_1+\vert \widetilde s\vert\psi_{y_1}(y^*))^2.
\end{equation}

Since $det\,\text{\bf B}_{n}(y^*,\zeta,\tau) = -{i}\alpha^+_
\mu(y^*,\zeta,\tau)
\vert \text{\bf B}_{1,1,n}\vert+ (-1)^{1+n}\frac{\lambda+2\mu}{\mu}(y^*)
({i}\xi_1-\vert\widetilde s\vert\psi_{y_1}(y^*))\vert \text{\bf B}
_{1,n,n}\vert$,
by (\ref{pardon}) and (\ref{pardon1}) we have
\begin{eqnarray}
\mbox{det}\,\text{\bf B}_{n}(y^*,\zeta^*,\tau)
= -{i}\alpha^+_\mu( y^*,\zeta,\tau)((-{i})^{n-1}
(\alpha^+_\mu(y^*,\zeta,\tau))^{n-2}\alpha^+_{\lambda+2\mu}(y^*,\zeta,\tau)
                                                       \nonumber\\
+ (-1)^{n-1}(-{i})^{n-3} (\alpha^+_\mu(y^*,\zeta,\tau))^{n-3}
\sum_{j=2}^{n-1}
(-{i}\xi_j+\vert s\vert\psi_{y_j}(y^*))^2            \nonumber\\
+(-1)^{1+n}\frac{\lambda+2\mu}{\mu}(y^*)(-{i}\alpha^+_
\mu(y^*,\zeta,\tau))^{n-2}
({i}\xi_1-\vert\widetilde
s^*\vert\psi_{y_1}(y^*))\frac{\mu}{\lambda+2\mu}(y^*)
(-{i}\xi_1+\vert\widetilde s\vert\psi_{y_1}(y^*))       \nonumber\\
= (-{i})^{n} (\alpha^+_\mu(y^*,\zeta,\tau))^{n-1}\alpha^+_{\lambda+2\mu}
(y^*,\zeta,\tau)                           \nonumber\\
+ (-{i})^{n-2}(-1)^{n-1} (\alpha^+_\mu(y^*,\zeta,\tau))^{n-2}\sum_{j=2}^{n-1}
(-{i}\xi_j+\vert s^*\vert\psi_{y_j}(y^*))^2)\nonumber\\
+ (-1)^{1+n}(-{i}\alpha^+_\mu(y^*,\zeta,\tau))^{n-2}(-{i}\xi_1
+\vert\widetilde s\vert\psi_{y_1}(y^*))^2.
\end{eqnarray}
The proof of the proposition is complete. $\blacksquare$
\\
\vspace{0.2cm}

By (\ref{3.35}) -(\ref{3.37}) if $det\,\text{\bf B}(y^*,\zeta^*,\tau)=0$ and
(\ref{4.1}) holds true, then
\begin{equation}\label{4.7}
\zeta^*\in \mathcal U= \left\{\zeta\in {\Bbb R}^{n+1}; \thinspace
\sum_{j=1}^{n-1} (\xi_j+{i}\vert\widetilde  s\vert\psi_{y_j}
(y^*))^2=\frac{\rho(y^*)(\xi_0+{i}\vert \widetilde s\vert
\psi_{y_0}(y^*))^2}
{(\lambda+3\mu)(y^*)}\right\}.
\end{equation}
Now we consider the two cases

{\bf Case A.} Let $\zeta^*\notin\mathcal U.$

In that case there exists a parametrix of the operator $\text{\bf
B}(\widetilde y,{\widetilde D},\tilde s,\tau)$  which we denote as
$\text{\bf B}^{-1}(
\widetilde y,{\widetilde D},\tilde s,\tau)$.  Then
\begin{eqnarray}\label{4.8}
\mbox{\bf w}_\nu(\cdot,0)
= \text{\bf B}^{-1}(\widetilde y,{\widetilde D},\tilde s,\tau)            \\
\times(V^+_{\mu}(1,n)
(\cdot,0)-r_{1,n,\nu},\dots ,V^+_{\mu}(n-1,n)(\cdot,0)-r_{n-1,n,\nu},
V^+_{\lambda+2\mu}(\cdot,0)-r_{2,\nu})
+ K\mbox{\bf w}_\nu(\cdot,0),\nonumber
\end{eqnarray}
where
\begin{equation}\label{4.9}
K: L^2(\Bbb R^n)\rightarrow H^{1,\widetilde s}(\Bbb R^n).
\end{equation}

By (\ref{4.4}), (\ref{4.8}) and (\ref{4.9}) for some positive constant $C_2$,
we have
\begin{eqnarray}\label{gol}
\sum_{k,j=1,k<j}^n \varXi_\mu(w_{k,j,\nu})+\varXi_{\lambda+2\mu}(w_{2,\nu})
\ge C_2{\vert \widetilde s\vert}\Vert\partial_{y_n} \mbox{\bf w}_\nu(\cdot,0)
\Vert^2_{L^2(\Bbb R^n)}
- C_3(\Vert
P_{\mu}(y,D,\tilde s,\tau) \mbox{\bf w}_{1,\nu}\Vert^2_{L^2(\mathcal Q)}
                                 \nonumber\\
+ \Vert P_{\lambda+2\mu}(y,D,\tilde s,\tau)w_{2,\nu} \Vert^2_{L^2(\mathcal Q)}
+{\vert \widetilde s
\vert}\Vert \mbox{\bf g}e^{\vert s\vert\varphi}
\Vert^2_{L^2({\Bbb R^n})}
+ \tau^{4n+6}\Vert \text{\bf w}\Vert^2_{H^{1,\widetilde s}(\mathcal Q)}).
\end{eqnarray}
By (\ref{4.9}), (\ref{gol}) and (\ref{klop}), we obtain (\ref{pulemet}).

{\bf Case B.} Let $\zeta^*\in\mathcal U.$

By (\ref{3.35})-(\ref{3.37}) there exists a constant $C_{4}>0$ such that
\begin{equation}\label{4.21}
\vert p_{\lambda+3\mu}(y,\widetilde \xi,0)-{\vert \widetilde s\vert}^2p
_{\lambda+3\mu}(y,\widetilde \nabla \psi,0)\vert
\le C_{4}\delta_1\vert\zeta\vert^2\quad \forall
(y,\zeta)\in \Upsilon_\ell\times \mathcal O(y^*,\delta_1(y^*)).
\end{equation}

If $\tilde s^*=0$, then using (\ref{4.21}), we obtain
\begin{eqnarray}
p_\beta(y,\widetilde \xi,0)-{\widetilde s}^2p_\beta(y,\widetilde
\nabla \psi,0)
                    \nonumber\\
=\sum_{j=1}^{n-1}(\lambda+3\mu-\beta)(y)(\xi^2_j-\widetilde s^2\psi^2_{y_j}
(y^*))+p_{\lambda+3\mu}(y,\widetilde \xi,0)-\widetilde s^2p_{\lambda+3\mu}
(y,\widetilde \nabla \psi,0)                 \nonumber\\
\ge \sum_{j=1}^{n-1}(\lambda+3\mu-\beta)(y)(\xi^2_j-\widetilde s^2
\psi^2_{y_j}(y^*))-
C_{5}\delta_1\vert\zeta\vert^2\quad\forall
(y,\zeta)\in \Upsilon_\ell\times \mathcal O(y^*,\delta_1(y^*)).\nonumber
\end{eqnarray}

Therefore for all sufficiently small $\delta_1$ there exists
$C_{6}$ such that
\begin{equation}\label{4.22}
p_\beta(y,\widetilde \xi,0)-\widetilde s^2p_\beta(y,\widetilde \nabla \psi,0)
\ge C_{6}\delta_1\vert\zeta\vert^2\quad\forall
(y,\zeta)\in \Upsilon_\ell\times \mathcal O(y^*,\delta_1(y^*)).
\end{equation}

Now let us consider the case when
$$
\widetilde s^*\ne 0.
$$
By the definition of the set $\mathcal U$, there exists  positive constants
$C_{7}, C_8$ such that
\begin{eqnarray}\label{4.185}
\vert
\rho(y^*)\xi_0^2-\widetilde s^2\rho(y^*)\psi^2_{y_0}(y^*)-(\lambda+3\mu)(y^*)
\sum_{j=1}^{n-1}(\xi^2_j-\widetilde s^2\psi^2_{y_j}(y^*)) \vert
\le C_{7}\delta\vert \zeta\vert^2,\nonumber\\
\vert
\xi_0(\rho\psi_{y_0})(y^*)-(\lambda+3\mu)(y^*)\sum_{j=1}^{n-1}\xi_j\psi_{y_j}
(y^*)\vert\le C_{8}\delta\vert \zeta\vert.
\end{eqnarray}

Let $\psi_{y_0}(y^*)=0.$
By (\ref{4.185})  we have
$$
\sum_{j=1}^{n-1}(\xi^2_j-\widetilde s^2\psi^2_{y_j}(y^*))\ge -C_{9}\delta
\vert \zeta\vert^2\quad \forall\zeta\in \mathcal O(y^*,\delta_1(y^*)).
$$
and
\begin{eqnarray}
\frak J_3(\mu,w_{k,n,\nu})\ge -{\vert \widetilde s\vert}C_{10}\delta
\Vert
(\partial_{y_n} \mbox{\bf w}_{\nu}(\cdot,0),\mbox{\bf  w}
_{\nu}(\cdot,0))\Vert^2
_{L^2(\Bbb R^n)\times H^{1,\widetilde s}(\Bbb R^n)}\nonumber\\-C_{11}\tau^{4n+6}
\Vert
(\partial_{y_n} \mbox{\bf w}(\cdot,0),\mbox{\bf w}
_{\nu}(\cdot,0))\Vert^2
_{L^2(\Bbb R^n)\times H^{1,\widetilde s}(\Bbb R^n)},                              \nonumber\\
\frak J_3(\lambda+2\mu,w_{2,\nu})
\ge -{\vert \widetilde s\vert}C_{12}\delta \Vert
(\partial_{y_n} \mbox{\bf w}_{\nu}(\cdot,0),\mbox{\bf  w}
_{\nu}(\cdot,0))\Vert^2
_{L^2(\Bbb R^n)\times H^{1,\widetilde s}(\Bbb R^n)}\nonumber\\
-C_{13}\tau^{4n+6}
\Vert
(\partial_{y_n} \mbox{\bf w}(\cdot,0),\mbox{\bf w}
_{\nu}(\cdot,0))\Vert^2
_{L^2(\Bbb R^n)\times H^{1,\widetilde s}(\Bbb R^n)}.
\end{eqnarray}

By (\ref{4.7}), we see
$$
\rho(y^*)(\tilde s^*\psi_{y_0}(y^*))^2=\frac12(\lambda+3\mu)(y^*)
$$
$$
\times \left\{
\sum_{j=1}^{n-1}(\tilde s^*\psi_{y_j}(y^*))^2
- \vert \widetilde\xi^*\vert^2
+ \sqrt{ \left(\sum_{j=1}^{n-1}(\tilde  s^*\psi_{y_j}(y^*))^2
- \vert \widetilde\xi^*\vert^2\right)^2
+ 4\left(\tilde s^*\sum_{j=1}^{n-1}\xi_j^*\psi_{y_j}(y^*)\right)^2}
\right\}
$$
$$
\le (\lambda+3\mu)(y^*)\sum_{j=1}^{n-1}\psi_{y_j}^2(y^*)(\tilde s^*)^2.
$$
This inequality and (\ref{giga}) imply
\begin{equation}\label{lobo}
\psi_{y_n}(y^*)>\root\of{\frac{(\lambda+2\mu)(y^*)}{\mu(y^*)}}\root\of{\sum
_{j=1}^{n-1}\vert
\psi_{y_j}(y^*)\vert^2}>\frac{\root\of{\rho(y^*)(\lambda+2\mu)(y^*)}}
{\root\of{\mu(\lambda+3\mu)(y^*)}}\vert \psi_{y_0}(y^*)\vert.
\end{equation}

Inequality (\ref{lobo}) yields
$$
p_\beta(y^*,\widetilde \xi,0)-(\widetilde s^*)^2p_\beta(y^*,\widetilde \nabla
\psi,0)
= \left (\frac{\rho(\lambda+2\mu-\beta)}{\lambda+3\mu}\right )(y^*)
((\xi_0^*)^2-(\widetilde s^*)^2(\psi_{y_0}(y^*))^2)
$$
$$
> \left(\frac{\rho(\lambda+2\mu-\beta)}{\lambda+3\mu}\right )(y^*)(\xi_0^*)^2
-\frac{(\mu(\lambda+2\mu-\beta))(y^*)}{(\lambda+2\mu)(y^*)}(\widetilde s^*)^2
\vert\psi_{y_n}(y^*)\vert^2> C_{14}\vert \zeta^*\vert^2.
$$
Hence there exists a positive constant $C_{15}$ such that
\begin{eqnarray}\label{vv}
\frak J_1(\mu,w_{k,n,\nu})+
\frak J_3(\mu,w_{k,n,\nu})
\ge {\vert \widetilde s\vert}C_{15}\delta \Vert
(\partial_{y_n} \mbox{\bf w}_{\nu}(\cdot,0),\mbox{\bf  w}
_{\nu}(\cdot,0))\Vert^2
_{L^2(\Bbb R^n)\times H^{1,\widetilde s}(\Bbb R^n)}\\
-C_{16}\tau^{4n+6}
\Vert
(\partial_{y_n} \mbox{\bf w}(\cdot,0),\mbox{\bf w}
_{\nu}(\cdot,0))\Vert^2
_{L^2(\Bbb R^n)\times H^{1,\widetilde s}(\Bbb R^n)},\nonumber
                                              \\
\frak J_1(\lambda+2\mu,w_{2,\nu})+\frak J_3(\lambda+2\mu,w_{2,\nu})
\ge {\vert \widetilde s\vert}C_{17}\delta
\Vert
(\partial_{y_n} \mbox{\bf w}_{\nu}(\cdot,0),\mbox{\bf  w}
_{\nu}(\cdot,0))\Vert^2
_{L^2(\Bbb R^n)\times H^{1,\widetilde s}(\Bbb R^n)}\nonumber\\
-C_{18}\tau^{4n+6}
\Vert
(\partial_{y_n} \mbox{\bf w}(\cdot,0),\mbox{\bf w}
_{\nu}(\cdot,0))\Vert^2
_{L^2(\Bbb R^n)\times H^{1,\widetilde s}(\Bbb R^n)}.                                    \nonumber
\end{eqnarray}

By $Re\, r_\beta(y^*,\zeta,\tau)>0$ for $\zeta^*\in \mathcal U$,
short computations yield
\begin{equation}\label{4.16}
\Gamma_\beta^+(y^*,\zeta^*,\tau)=\mbox{\bf
I}_\beta\,\root\of{Re\, r_\beta(y^*,\zeta^*,\tau)},\quad \mbox{\bf
I}_\beta=\mbox{sign}(lim_{\zeta\rightarrow \zeta^*}  \,\,Im \,
r_{\mu}(y^*,\zeta,\tau)/\vert \widetilde s\vert).
\end{equation}

Therefore for any $(\widetilde y,\zeta)$ from $ \Upsilon_\ell\times \mathcal O
(y^*,\delta_1(y^*))$
\begin{equation}
-\mbox{Re}\,\{ \Gamma_\beta^+(\widetilde y,0,\widetilde  \xi,\widetilde s,\tau)
(\nabla_{\widetilde \xi}p_{\beta}(y^*,\widetilde\xi ,0),\widetilde \nabla
\psi(y^*))\}\ge \epsilon(\delta_1,\delta_2)\vert(\widetilde \xi,\widetilde s)
\vert^2.
\end{equation}
By (\ref{4.16}), Lemmata \ref{Fops3} and \ref{Fops5} (the G\aa rding inequality
(\ref{goblin})) in Section \ref{Q!8}, we have
\begin{eqnarray}\label{4.19}
\frak J_2(\mu,w_{k,n,\nu})=-
\text{Re}\int_{{\Bbb R^n}}2\vert\widetilde s\vert\mu(y^*) {i}
\Gamma^+_{\mu}(\widetilde y,0,{\widetilde D},\tilde s,\tau)
w_{k,n,\nu}\overline{(
\nabla_{\widetilde \xi}p_\mu(y^*,\widetilde\nabla w_{k,n,\nu},0), \widetilde
\nabla \psi(y^*))}d \widetilde y
                                   \nonumber\\
-\text{Re}\int_{{\Bbb R^n}}2\vert \widetilde s\vert\mu(y^*)
{i}V^+_{\mu}(k,n)(\cdot,0)\overline{(\nabla_{\widetilde \xi}p_
\mu(y^*,\widetilde\nabla w_{k,n,\nu},0), \widetilde \nabla \psi(y^*))}d
\widetilde y
                                 \nonumber\\
\ge -C_{19}(\epsilon(\delta,\delta_1){\vert \widetilde s\vert+\tau^{4n+6})}
\Vert
(\partial_{y_n} \mbox{\bf w}_{\nu}(\cdot,0),\mbox{\bf  w}
_{\nu}(\cdot,0))\Vert^2
_{L^2(\Bbb R^n)\times H^{1,\widetilde s}(\Bbb R^n)}
                                                       \\
-C_{20}(\delta,\delta_1)(\Vert P_{\mu}(y,D,\tilde s,\tau) w_{k,n,\nu}\Vert^2
_{L^2(\mathcal Q)}+\Vert P_{\lambda+2\mu}(y,D,\tilde s,\tau)
w_{2,\nu}\Vert^2 _{L^2(\mathcal Q)}+\tau^{4n+6}\Vert \text{\bf
w}\Vert^2_{H^{1,\widetilde s}(\mathcal Q)})\nonumber
\end{eqnarray}
and
\begin{eqnarray}\label{4.24}
\frak J_2(\lambda+2\mu, w_{2,\nu})              \nonumber\\
= -\text{Re}\int_{{\Bbb R^n}}\vert
s\vert(\lambda+2\mu)(y^*) {i}  \Gamma^+_{\lambda+2\mu}
(\widetilde y,0,{\widetilde D},\tilde s,\tau) w_{2,\nu}
\overline{(\nabla_{\widetilde \xi}
p_{\lambda+2\mu}(y^*,\widetilde\nabla w_{2,\nu},0), \widetilde
\nabla \psi(y^*))}
d \widetilde y                     \nonumber \\
-\text{Re}\int_{{\Bbb R^n}}\vert \widetilde s\vert(\lambda+2\mu)(y^*)
{i} V^+_{\lambda+2\mu}(\cdot,0)\overline{(\nabla_{\widetilde \xi}
p_{\lambda+2\mu}(y^*,\widetilde\nabla w_{2,\nu},0), \widetilde \nabla
\psi(y^*))}d \widetilde y\nonumber\\
\ge -C_{21}\epsilon(\delta,\delta_1)\vert \widetilde s\vert\Vert
(\partial_{y_n} \mbox{\bf w}_{\nu}(\cdot,0),\mbox{\bf  w}
_{\nu}(\cdot,0))\Vert^2
_{L^2(\Bbb R^n)\times H^{1,\widetilde s}(\Bbb R^n)}\\
-C_{22}\tau^{4n+6}
\Vert
(\partial_{y_n} \mbox{\bf w}(\cdot,0),\mbox{\bf w}
_{\nu}(\cdot,0))\Vert^2
_{L^2(\Bbb R^n)\times H^{1,\widetilde s}(\Bbb R^n)}\nonumber\\
-C_{23}(\delta,\delta_1)(\Vert P_{\mu}(y,D,\tilde s,\tau) w_{k,n,\nu}\Vert^2
_{L^2(\mathcal Q)}+\Vert P_{\lambda+2\mu}(y,D,\tilde s,\tau)
w_{2,\nu}\Vert^2 _{L^2(\mathcal Q)}+\tau^{4n+6}\Vert \text{\bf
w}\Vert^2_{H^{1,\widetilde s}(\mathcal Q)}).\nonumber
\end{eqnarray}

By (\ref{vv}), (\ref{4.19}), (\ref{4.24}) and (\ref{2.1}), there exists
a constant $C_{24}>0$ such that
\begin{eqnarray}\label{4.23}
\sum_{k,j=1, k<j}^n\varXi_\mu(w_{k,j,\nu})+\varXi_\mu(w_{2,\nu})\ge
C_{24}{\vert \widetilde s\vert} \Vert
(\partial_{y_n} \mbox{\bf w}_{\nu}(\cdot,0),\mbox{\bf  w}
_{\nu}(\cdot,0))\Vert^2
_{L^2(\Bbb R^n)\times H^{1,\widetilde s}(\Bbb R^n)}-\nonumber\\
-C_{25}\tau^{4n+6}
\Vert
(\partial_{y_n} \mbox{\bf w}(\cdot,0),\mbox{\bf w}
_{\nu}(\cdot,0))\Vert^2
_{L^2(\Bbb R^n)\times H^{1,\widetilde s}(\Bbb R^n)}\\
-C_{26}(\delta,\delta_1)(\Vert P_{\mu}(y,D,\tilde s,\tau)
\mbox{\bf w}_{1,\nu}\Vert^2
_{L^2(\mathcal Q)}+\Vert P_{\lambda+2\mu}(y,D,\tilde s,\tau)
w_{2,\nu}\Vert^2 _{L^2(\mathcal Q)}+\tau^{4n+6}\Vert \text{\bf
w}\Vert^2_{H^{1,\widetilde s}(\mathcal Q)}).\nonumber
\end{eqnarray}
From (\ref{4.23}) and (\ref{klop}), we obtain (\ref{pulemet}).
Thus in each case, we established the estimate  (\ref{pulemet}).
\section{\bf End of the proof.}\label{Q2}
First we finish the proof of the proposition \ref{zoopa}.
Now let us take the covering of the sphere $\Bbb S^{n}$  by
conical neighborhoods $\mathcal O(\zeta^*,\delta_1(\zeta^*)).$
From this covering we take the finite subcovering $\cup_{m=1}^N\mathcal O
(\zeta_m^*,\delta_1(\zeta_m^*))$.  Let $\chi_m$ be the partition of
unity associated to this subcovering. Hence $\sum_{m=0}^N\chi_m
(\widetilde \xi,s)\equiv 1.$  Then by (\ref{pulemet})
\begin{eqnarray}\label{golifax}
{\vert \widetilde s\vert}\Vert(
\partial_{y_n} \mbox{\bf w}(\cdot,0),\mbox{\bf  w}(\cdot,0)
)\Vert^2_{L^2(\Bbb R^n)\times H^{1,\widetilde s}(\Bbb R^n)}
+ {\vert \widetilde s\vert}\tau\Vert
\mbox{\bf w}\Vert^2_{H^{1,\widetilde s}(\mathcal Q)}
                         \\
\le\sum_{m=0}^N \left({\vert \widetilde s\vert}\Vert \eta_\ell(
\partial_{y_n} (\chi_m\mbox{\bf w})(\cdot,0),\chi_m\mbox{\bf w}
(\cdot,0))\Vert^2_{L^2(\Bbb R^n)\times H^{1,\widetilde s}
(\Bbb R^n)}
+ {\vert \widetilde s\vert}\tau\Vert\eta_\ell
\chi_m\mbox{\bf w}\Vert^2_{H^{1,\widetilde s}(\mathcal Q)}\right)+\nonumber\\
 C_1\sum_{m=0}^N \left({\vert \widetilde s\vert}\Vert (1-\eta_\ell)(
\partial_{y_n} (\chi_m\mbox{\bf w})(\cdot,0),\chi_m\mbox{\bf w}
(\cdot,0))\Vert^2_{L^2(\Bbb R^n)\times H^{1,\widetilde s}
(\Bbb R^n)}
+ {\vert \widetilde s\vert}\tau\Vert(1-\eta_\ell)
\chi_m\mbox{\bf w}\Vert^2_{H^{1,\widetilde s}(\mathcal Q)}\right)\nonumber\\
\le C_{2}\biggl({\tau^{4n+6}}\Vert \text{\bf w}
\Vert^2_{H^{1,\widetilde s}(\mathcal Q)} +{\vert \widetilde s\vert}\Vert
\mbox{\bf g}e^{\vert s\vert \varphi}\Vert^2_{L^2(\Bbb R^n)}
+ \sum_{m=0}^N\Vert P_{\mu}(y,D,\tilde s,\tau)\mbox{\bf w}_{1,m}\Vert^2
_{L^2(\mathcal Q)}\nonumber\\
+ \Vert P_{\lambda+2\mu}(y,D,\tilde s,\tau)
w_{2,m}\Vert^2_{L^2(\mathcal Q)}
+ {\tau^{4n+6}}\Vert (
\partial_{y_n} \mbox{\bf w}(\cdot,0),\mbox{\bf  w}(\cdot,0)
)\Vert^2_{L^2(\Bbb R^n)\times H^{1,\widetilde s}
(\Bbb R^n)}\biggr)+                     \nonumber\\
C_3\sum_{m=0}^N \left({\vert \widetilde s\vert}\Vert (1-\eta_\ell)(
\partial_{y_n} (\chi_m\mbox{\bf w})(\cdot,0),\chi_m\mbox{\bf w}
(\cdot,0))\Vert^2_{L^2(\Bbb R^n)\times H^{1,\widetilde s}
(\Bbb R^n)}
+ {\vert \widetilde s\vert}\tau\Vert(1-\eta_\ell)
\chi_m\mbox{\bf w}\Vert^2_{H^{1,\widetilde s}(\mathcal Q)}\right)\nonumber .
\end{eqnarray}
By (\ref{nikolas})
\begin{eqnarray}\label{nikolas1}
\sum_{m=0}^N \left({\vert \widetilde s\vert}\Vert (1-\eta_\ell)(
\partial_{y_n} (\chi_m\mbox{\bf w})(\cdot,0),\chi_m\mbox{\bf w}
(\cdot,0))\Vert^2_{L^2(\Bbb R^n)\times H^{1,\widetilde s}
(\Bbb R^n)}
+ {\vert \widetilde s\vert}\tau\Vert(1-\eta_\ell)\chi_m
\mbox{\bf w}\Vert^2_{H^{1,\widetilde s}(\mathcal Q)}\right)\nonumber\\
\le
C_4\left(\tau^{4n+6}\Vert (
\partial_{y_n} \mbox{\bf w}(\cdot,0),\mbox{\bf w}
(\cdot,0))\Vert^2_{L^2(\Bbb R^n)\times H^{1,\widetilde s}
(\Bbb R^n)}
+ \tau^{4n+6}\Vert
\mbox{\bf w}\Vert^2_{H^{1,\widetilde s}(\mathcal Q)}\right)
\end{eqnarray}
Using  inequality (\ref{nikolas1}) in order to estimate the last terms in
(\ref{golifax}) we obtain
\begin{eqnarray}
{\vert \widetilde s\vert}\Vert (
\partial_{y_n} \mbox{\bf w}(\cdot,0),\mbox{\bf  w}(\cdot,0)
)\Vert^2_{L^2(\Bbb R^n)\times H^{1,\widetilde s}(\Bbb R^n)}
+ {\vert \widetilde s\vert}\tau\Vert
\mbox{\bf w}\Vert^2_{H^{1,\widetilde s}(\mathcal Q)}
                         \nonumber\\
\le C_{5}\biggl(\tau^{4n+6}\Vert \text{\bf w}
\Vert^2_{H^{1,\widetilde s}(\mathcal Q)} +{\vert \widetilde s\vert}\Vert
\mbox{\bf g}e^{\vert s\vert \varphi}\Vert^2_{L^2({\Bbb R^n})}
                                          \nonumber\\
+ \Vert P_{\mu}(y,D,\tilde s,\tau)\mbox{\bf
w}_{1}\Vert^2_{L^2(\mathcal Q)}+\Vert P_{\lambda+2\mu}(y,D,\tilde s,\tau)
w_{2}\Vert^2_{L^2(\mathcal Q)}\nonumber\\
+ {\vert \widetilde s\vert}\Vert(
\partial_{y_n} \mbox{\bf w}(\cdot,0),\mbox{\bf w}(\cdot,0)
)\Vert^2_{L^2(\Bbb R^n)\times H^{1,\widetilde s}
(\Bbb R^n)} \nonumber\\
+ \sum_{m=0}^N\Vert[\chi_m, P_{\mu}(y,D,\tilde s,\tau)]\mbox{\bf
w}_{1}\Vert^2_{L^2(\mathcal Q)}
+ \Vert [\chi_m, P_{\lambda+2\mu}(y,D,\tilde s,\tau)]
w_{2}\Vert^2_{L^2(\mathcal Q)}\biggr)  +           \nonumber\\
C_6\left(\tau^{4n+6}\left\Vert \left(
\partial_{y_n}\mbox{\bf w}(\cdot,0),\mbox{\bf w}
(\cdot,0)\right)\right\Vert^2_{L^2(\Bbb R^n)\times H^{1,\widetilde s}
(\Bbb R^n)}
+ \tau^{4n+6}\Vert
\mbox{\bf w}\Vert^2_{H^{1,\widetilde s}(\mathcal Q)}\right)\nonumber\\
\le C_{7}\biggl(\tau^{4n+6}\Vert \text{\bf w}
\Vert^2_{H^{1,\widetilde s}(\mathcal Q)} +{\vert \widetilde s\vert}\Vert
\mbox{\bf g}e^{\vert s\vert\varphi}\Vert^2_{L^2({\Bbb R^n})}
                                   \nonumber\\
+ \Vert P_{\mu}(y,D,\tilde s,\tau)\mbox{\bf
w}_{1}\Vert^2_{L^2(\mathcal Q)}+\Vert P_{\lambda+2\mu}(y,D,\tilde s,\tau)
w_{2}\Vert^2_{L^2(\mathcal Q)}\nonumber\\
+ {\vert \widetilde s\vert}\Vert (
\partial_{y_n} \mbox{\bf w}(\cdot,0),\mbox{\bf  w}(\cdot,0)
)\Vert^2_{L^2(\Bbb R^n)\times H^{1,\widetilde s}
(\Bbb R^n)}\biggr)+\nonumber\\
C_8\left(\tau^{4n+6}\Vert (
\partial_{y_n}\mbox{\bf w}(\cdot,0),\mbox{\bf w}
(\cdot,0))\Vert^2_{L^2(\Bbb R^n)\times H^{1,\widetilde s}
(\Bbb R^n)}
+ \tau^{4n+6}\Vert
\mbox{\bf w}\Vert^2_{H^{1,\widetilde s}(\mathcal Q)}\right).
\end{eqnarray}
Hence there exists $\tau_0>1$  independent of $s$ such that for all
$\tau\ge\tau_0$ and $ s\ge 1$, we see
\begin{eqnarray}\label{zmpolit}
{\vert \widetilde s\vert}\Vert(
\partial_{y_n} \mbox{\bf w}(\cdot,0),\mbox{\bf  w}(\cdot,0)
)\Vert^2_{L^2(\Bbb R^n)\times H^{1,\widetilde s}(\Bbb R^n)}
+ {\vert \widetilde s\vert}\tau\Vert
\mbox{\bf w}\Vert^2_{H^{1,\widetilde s}(\mathcal Q)}\nonumber\\
\le C_{9}\biggl(\tau^{4n+6}\Vert \text{\bf w}
\Vert^2_{H^{1,\widetilde s}(\mathcal Q)} +{\vert \widetilde s\vert}\Vert
\mbox{\bf g}e^{ s\varphi}\Vert^2_{L^2({\Bbb R^n})}\nonumber\\
+ \Vert P_{\mu}(y,D,\tilde s,\tau)\mbox{\bf
w}_{1}\Vert^2_{L^2(\mathcal Q)}+\Vert P_{\lambda+2\mu}(y,D,\tilde s,\tau)
w_{2}\Vert^2_{L^2(\mathcal Q)}  \biggr).
\end{eqnarray}

If in the inequality (\ref{zmpolit}) we return to the function $\mbox{\bf v}$,
then we obtain
\begin{equation}\label{zmpolit1}
\int_{\Sigma}(s\tau\varphi\vert \nabla \mbox{\bf  v}\vert^2 +s^3\tau^3\varphi
^3\vert  \mbox{\bf  v}\vert^2 )e^{2s\varphi}d\Sigma + \int_{Q}(
s\tau^2\varphi\vert
\nabla\mbox{\bf v}\vert^2+ s^3\tau^4\varphi^3\vert
\mbox{\bf v}\vert^2 )e^{2s\varphi}dx
\end{equation}
$$
\le C_{10}\left(s\tau\Vert\root\of{\varphi}
\mbox{\bf g}e^{s\varphi}\Vert^2_{L^2(\Sigma)}+\Vert \mbox{\bf q}
e^{s\varphi}\Vert^2_{L^2(Q)}
+\int_{\widetilde \Sigma}(s\tau\varphi\vert \nabla \mbox{\bf  v}\vert^2
+s^3\tau^2\varphi^3\vert  \mbox{\bf  v}\vert^2 )e^{2s\varphi}d\Sigma
\right).
$$

Proof of the proposition \ref{zoopa} is complete.

Now we show that the functions $\mbox{\bf v}$ given by (\ref{victory}) satisfies the estimate (\ref{3.2'1}).

First we show that function $\mbox{\bf v}$ satisfies the boundary conditions (\ref{poko1}).
By (\ref{gora1A}) and the zero Dirichlet boundary conditions for the function $\mbox{\bf u}$ we have
\begin{equation}\label{gora1AA}
- L_{\lambda,\mu}(x,D')\mbox{\bf u}
+\int_0^{x_0}L_{\widetilde \lambda, \widetilde\mu}
(x, \widetilde x_0, D')\mbox{\bf u}(\widetilde x_0, x')d\widetilde x_0
= \mbox{\bf F} \quad \mbox{on}\,\,\Sigma
\end{equation}
Next we move all terms containing the first  derivatives of the function $\mbox{\bf u}$ in the right hand side, divide the both sides by $\mu$  and denote the right hand side of obtained equality  as $g_k$. We have
\begin{eqnarray}\label{goraAAA}
-\Delta \mbox{\bf u}-\frac{(\lambda+\mu)}{\mu}\nabla'\mbox{div}\,\mbox{\bf u} +\int_0^{x_0}\left(\frac{\tilde \mu(x,\tilde x_0)}{\mu(x)}\Delta \mbox{\bf u}(\tilde x_0,x')+\frac{\tilde \lambda(x,\tilde x_0)+\tilde \mu(x,\tilde x_0)}{\mu(x)}\nabla'\mbox{div}\,\mbox{\bf u}\right) dx_0\nonumber\\=(g_1,\dots, g_n).
\end{eqnarray}
 For any index $\widehat i\in\{1,\dots, n\}$ the short computations imply
\begin{eqnarray}\label{gavnoed}
-\Delta u_{\widehat i} =\sum_{j=1, j\ne \widehat i}^n
-\partial_{x_j}(\partial_{x_j} u_{\widehat i}-\partial_{x_{\widehat i}} u_j)
-\partial^2_{x_{\widehat i}} u_{\widehat i} -\sum_{j=1,j\ne \widehat i}^n
\partial_{x_j}\partial_{x_{\widehat i} }u_j\nonumber\\
=-\sum^n_{j=1, j\ne \widehat i}\partial_{x_j}(\partial_{x_j} u_{\widehat i}
-\partial_{x_{\widehat i}} u_j)-\partial_{x_{\widehat i} }\text{div}\,
\mbox{\bf u }\quad \text{in}\,\,\Omega.
\end{eqnarray}
Using the equality (\ref{gavnoed}) we rewrite $ \widetilde j$-th equation in (\ref{goraAAA}) as
$$
b_{\hat j}(x,D)\mbox{\bf v}= -\sum^n_{j=1, j\ne \widehat i}\mbox{sign}(j-\hat i)\partial_{x_j} v_{\hat i j}-\frac{\lambda+2\mu}{\mu}\partial_{x_{\hat j}}v_2=g_{\hat j}
$$
Construction of the operator $\widetilde{\mbox{B}}(x,D)$ is complete. Next $v_{i,j}(x)=\nu_j(x')\frac{\partial u_i}{\partial\vec \nu}-\nu_i(x')\frac{\partial u_j}{\partial\vec \nu}, i<j.$ By()
$$
v_{i,j}(y^*)=0 \quad \mbox{for} \,\,j<n,\quad v_{i,n}(y^*)=-\frac{\partial u_j}{\partial x_n}(y^*),\quad v_2=-\frac{\partial u_j}{\partial x_n}(y^*).
$$
Set $\tilde{\mbox{\bf v}}=(v_{1,n},\dots, v_{2,n},\dots v_{n-1,n}, v_2).$ Obviously in the small neighborhood of  $y^*$ there exists a  smooth matrix $\mathcal C_0$ such that
$$
(\frac{\partial u_1}{\partial x_n},\dots, \frac{\partial u_n}{\partial x_n})=\mathcal C_0(x') \tilde{\mbox{\bf v}}\quad \forall x\in \Sigma\cap B(y^*,\delta).
$$  Then $\mathcal C(x')\mbox{\bf v}=\mbox{\bf v}-(\nu_2\frac{\partial u_1}{\partial x_2}-\nu_1 \frac{\partial u_2}{\partial x_1}, \dots , \nu_n \frac{\partial u_1}{\partial x_n}-\nu_1 \frac{\partial u_n}{\partial x_1}, \dots , \nu_n \frac{\partial u_{n-1}}{\partial x_n}-\nu_{n-1}\frac{\partial u_{n}}{\partial x_{n-1}}, \sum_{j=1}^n \nu_j\frac{\partial u_j}{\partial x_j})=0.$

Still we  can  not apply the proposition \ref{zoopa} directly since the traces of the function  $\mbox{\bf v}$ and its time derivative at moment $x_0=\pm T$ are  not equal zero. On the other hand $\mbox{\bf v}(\pm T,\cdot),\partial_{x_0}\mbox{\bf v}(\pm T,\cdot)\in H^2(\Omega).$
So for sufficiently small positive $\epsilon$ we extend the function $\mbox{\bf v}$ on $Q_\epsilon=(-T-\epsilon, T+\epsilon)\times \Omega$ by the formula $\mbox{\bf v}=\gamma(x_0)( \mbox{\bf v}(T,\cdot)+(x_0-T)\partial_{x_0}\mbox{\bf v}(T,\cdot))$ for $x_0\in (T,T+\epsilon)$ and $\mbox{\bf v}=\gamma(x_0)( \mbox{\bf v}(-T,\cdot)+(x_0+T)\partial_{x_0}\mbox{\bf v}(-T,\cdot))$ for $x_0\in (T-\epsilon,T),$
where $\gamma\in C^\infty_0[-T-\epsilon, T+\epsilon]$ and $\gamma\vert_{[-T-\epsilon/2, T+\epsilon/2]}=1.$ Provided that $\epsilon$ is sufficiently small we extend the functions $\rho,\mu,\lambda$ on $Q_\epsilon$ such that they are positive, keeping the same regularity and we extend the function $\psi$ on $Q_\epsilon$ such  that Condition 1.1 holds on $Q_\epsilon$, inequality (\ref{giga}) holds true on $[-T-\epsilon,T+\epsilon]\times \Gamma_0$
and
\begin{equation}\label{pinok222}
\partial_{x_0}\psi(x)<0 \quad \mbox{on}\,\,(0,T+\epsilon]\times \bar \Omega
\quad \mbox{and}\quad \partial_{x_0}\psi(x)>0 \quad \mbox{on}\,\,[-T-\epsilon,0)\times
\www\Omega.
\end{equation} To this extended function $\mbox{\bf v}$ we can now apply the proposition \ref{zoopa}.

In Klibanov \cite{Kl2}, the following inequality is proved:
let $R(x,\tilde x_0,\widetilde D)$ be a second-order differential
operator with smooth coefficients.
Then
\begin{equation}\label{qupol}
\left\Vert e^{s\varphi} \int_0^{x_0}R(x,\tilde x_0,\widetilde D)\mbox{ \bf u}
(\tilde x_0,\cdot)d\tilde x_0\right\Vert_{L^2[-T,T]}
\le o(\frac{1}{s})\sum_{\vert\alpha\vert\le 2}
\Vert  e^{s\varphi}\partial_x^{\alpha}\mbox{ \bf u}
\Vert_{L^2[-T,T]}.
\end{equation}
We set $\Sigma_\epsilon=(-T-\epsilon,T+\epsilon)\times \partial\Omega.$
By (\ref{pinok222})
\begin{eqnarray}
\label{lenin}s\tau\Vert \root\of{\varphi}\mbox{\bf g}e^{s\varphi}\Vert^2_{L^2(\Sigma_\epsilon)}\le
\frac{C_{11}}{s}   (\sum_{k=0}^1\Vert\root\of{\varphi}  (\mbox{\bf v}((-1)^k T,\cdot) ),\partial_{x_0}\mbox{\bf v}((-1)^k T,\cdot) )e^{s(\varphi((-1)^k T,\cdot)}\Vert^2_{H^{1}(\partial\Omega)}\nonumber\\+\Vert s\tau\varphi^\frac 32  (\mbox{\bf v}((-1)^k T,\cdot) ),\partial_{x_0}\mbox{\bf v}((-1)^k T,\cdot) )e^{s\varphi((-1)^k T,\cdot)}\Vert^2_{L^2(\partial\Omega)})\nonumber\\
+C_{12}(\int_{\widetilde \Sigma}(s\tau\varphi\vert \nabla \mbox{\bf  v}\vert^2
+s^3\tau^2\varphi^3\vert  \mbox{\bf  v}\vert^2 )e^{2s\varphi}d\sigma+s\tau\Vert \root \of{\varphi} \F e^{s\varphi}\Vert^2_{L^2(\Sigma)})
\end{eqnarray}
and
\begin{equation}\label{lenin1}
\Vert \mbox{\bf q}e^{s\varphi}\Vert^2_{L^2(Q_\epsilon)}\le \frac{C_{13}}{s} \sum_{k=0}^1\Vert\Delta (\mbox{\bf v}((-1)^k T,\cdot) ,\partial_{x_0}\mbox{\bf v}((-1)^k T,\cdot) )e^{s\varphi((-1)^k T,\cdot)}\Vert^2_{L^2(\Omega)}+\Vert \mbox{\bf q}e^{s\varphi}\Vert^2_{L^2(Q)}.
\end{equation}
By (\ref{gora2A}) and (\ref{victory})
\begin{equation}\label{victory1}
\mbox{\bf v}(\pm T,\cdot)=-\left ( \int_0^{\pm T}\frac{\widetilde
\mu(\pm T, x',\widetilde x_0)}{\mu(\pm T,x')} d\omega_{\mbox{\bf u}} \,d\widetilde x_0,
\int_0^{\pm T}\frac{(\widetilde \lambda
+2\widetilde \mu)(\pm T,x',\widetilde x_0)}{(\lambda+2\mu)(\pm T,x')} \text{div}\,
\mbox{\bf u}\, d\widetilde x_0\right),
\end{equation}
and
\begin{eqnarray}\label{victory1l}
\partial_{x_0}\mbox{\bf v}(\pm T,\cdot)=\\-\left ( \int_0^{\pm T}\partial_{x_0}\left(\frac{\widetilde
\mu(\pm T,x',\widetilde x_0)}{\mu(\pm T,x')} \right)d\omega_{\mbox{\bf u}} \,d\widetilde x_0,
\int_0^{\pm T}\partial_{x_0}\left (\frac{(\widetilde \lambda
+2\widetilde \mu)(\pm T,x',\widetilde x_0)}{(\lambda+2\mu)(\pm T,x')}\right) \text{div}\,
\mbox{\bf u}\, d\widetilde x_0\right).\nonumber
\end{eqnarray}

We have
\begin{proposition} \label{oxX} There exists a constant $C_{14}$ independent of $s$ and $\tau$ such that
\begin{eqnarray}
\Vert\Delta  \mbox{\bf v}(\pm T,\cdot)e^{\pm s\varphi(T,\cdot)}\Vert_{L^2(\Omega)}+\Vert\partial_{x_0}\Delta \mbox{\bf v}(\pm T,\cdot)e^{\pm s\varphi(T,\cdot)}\Vert_{L^2(\Omega)}\nonumber\\
\le C_{14}(\root\of{s\tau}\Vert\varphi^\frac 12\nabla  \mbox{\bf v}e^{s\varphi}\Vert_{L^2(Q)}+(s\tau)^\frac 32\Vert\varphi^\frac 32 \mbox{\bf v}e^{s\varphi}\Vert_{L^2(Q)}+
\Vert \mbox{\bf q} e^{s\varphi}\Vert_{L^2(Q)}) .
\end{eqnarray}
\end{proposition}

{\bf Proof.}
Let $\theta(x)$ be a smooth function. Integrating the first equation in (\ref{3.3}) in $x_0$ on $[0,\tilde T]$ we obtain
\begin{eqnarray}\label{ox}(\partial_{x_0} d\omega_{\mbox{\bf u}}\theta)(\tilde T,x')
-(\partial_{x_0} d\omega_{\mbox{\bf u}}\theta)(0,x')-\int_0^{\tilde T}\partial_{x_0}\theta\partial_{x_0} d\omega_{\mbox{\bf u}}dx_0\\ -\int_0^{\tilde T}\theta\mu  \Delta d\omega_{\mbox{\bf u}}d\tilde x_0-\int_0^{\tilde T}\theta\int_0^{x_0} \tilde \mu \Delta d\omega_{\mbox{\bf u}}d\tilde x_0dx_0=\int_0^{\tilde T}\theta\mbox{\bf q}_1dx_0.\nonumber
\end{eqnarray}
We set $\theta(x)=\tilde \mu/\mu, $ $Y(\tilde T,x')=\int_0^{\tilde T}\tilde \mu \Delta  d\omega_{\mbox{\bf u}}d x_0$
and $\alpha(\tilde T,x')=(\partial_{x_0} d\omega_{\mbox{\bf u}}\theta)(\tilde T,x')
-(\partial_{x_0} d\omega_{\mbox{\bf u}}\theta)(0,x')-\int_0^{\tilde T}\partial_{x_0}\theta\partial_{x_0} d\omega_{\mbox{\bf u}}dx_0-\int_0^{\tilde T}\theta\mbox{\bf q}_1dx_0.$
Then
$
Y(\tilde T,x')=-\int_0^{\tilde T}\frac{\tilde \mu}{\mu}(\tilde x_0,x') Y(\tilde x_0,x')d x_0+\alpha(\tilde T,x').
$
By Gronwall's inequality from (\ref{ox}) we have
$$
Y(\tilde T,x')\le \sup_{x_0\in[0,T]}\alpha(x_0,x')exp(\int_0^{\tilde T}\frac{\tilde \mu}{\mu}(\tilde x_0,x')d\tilde x_0).
$$
Then
$$\mbox{\bf v}(\pm T,\cdot)e^{\pm s\varphi(T,x')}=Y(T,\cdot)e^{\pm s\varphi(T,x')}\le\sup_{x_0\in[0,T]}\alpha(x_0,x')exp(\int_0^{ T}\frac{\tilde \mu}{\mu}(\tilde x_0,x')d\tilde x_0)e^{\pm s\varphi(T,x')}.
$$
So
$$\Vert\mbox{\bf v}(\pm T,\cdot)e^{\pm s\varphi(T,\cdot)}\Vert_{L^2(\Omega)}\le C_{15}(\sup_{x_0\in[0,T]}(\Vert \partial_{x_0}\mbox{\bf w}(x_0,\cdot)\Vert_{L^2(\Omega)}+s\tau\Vert \varphi\mbox{\bf w}(x_0,\cdot)\Vert_{L^2(\Omega)})+C_{16}\Vert \mbox{\bf q} e^{s\varphi}\Vert_{L^2(Q)}
$$
$$
\le  C_{17}(\root\of{s\tau}\Vert\varphi^\frac 12\nabla  \mbox{\bf v}e^{s\varphi}\Vert_{L^2(Q)}+(s\tau)^\frac 32\Vert\varphi^\frac 32 \mbox{\bf v}e^{s\varphi}\Vert_{L^2(Q)}+
\Vert \mbox{\bf q} e^{s\varphi}\Vert_{L^2(Q)}).
$$
Estimate for $\partial_{x_0}\Delta \mbox{\bf v}(\pm T,\cdot)$ proved similarly.
$\blacksquare$

Using proposition  \ref{zoopa} and (\ref{lenin}) -(\ref{lenin1}) we obtain

\begin{eqnarray}\label{zmpolit11}
\int_{\Sigma_\epsilon}(s\tau\varphi\vert \nabla \mbox{\bf  v}\vert^2 +s^3\tau^3\varphi
^3\vert  \mbox{\bf  v}\vert^2 )e^{2s\varphi}d\Sigma + \int_{Q_\epsilon}(
s\tau^2\varphi\vert
\nabla\mbox{\bf v}\vert^2+ s^3\tau^4\varphi^3\vert
\mbox{\bf v}\vert^2 )e^{2s\varphi}dx\\
\le C_{18}\left(s\tau\Vert\root\of{\varphi}
\mbox{\bf g}e^{s\varphi}\Vert^2_{L^2(\Sigma_\epsilon)}+\Vert \mbox{\bf q}
e^{s\varphi}\Vert^2_{L^2(Q_\epsilon)}
+\int_{\widetilde \Sigma_\epsilon}(s\tau\varphi\vert \nabla \mbox{\bf  v}\vert^2
+s^3\tau^2\varphi^3\vert  \mbox{\bf  v}\vert^2 )e^{2s\varphi}d\Sigma
\right)\nonumber\\
\le C_{19}\left(s\tau\Vert\root\of{\varphi}
\mbox{\bf g}e^{s\varphi}\Vert^2_{L^2(\Sigma)}+\Vert \mbox{\bf q}
e^{s\varphi}\Vert^2_{L^2(Q)}
+\int_{\widetilde \Sigma}(s\tau\varphi\vert \nabla \mbox{\bf  v}\vert^2
+s^3\tau^2\varphi^3\vert  \mbox{\bf  v}\vert^2 )e^{2s\varphi}d\Sigma
\right)\nonumber\\
+ \frac{C_{20}}{s} \sum_{k=0}^1\Vert\Delta (\mbox{\bf v}((-1)^k T,\cdot) ,\partial_{x_0}\mbox{\bf v}((-1)^k T,\cdot) )e^{s\varphi((-1)^k T,\cdot)}\Vert^2_{L^2(\Omega)}\nonumber\\
+\frac{C_{21}}{s}   (\sum_{k=0}^1\Vert\root\of{\varphi}  (\mbox{\bf v}((-1)^k T,\cdot) ),\partial_{x_0}\mbox{\bf v}((-1)^k T,\cdot) )e^{s\varphi((-1)^k T,\cdot)}\Vert^2_{H^{1}(\partial\Omega)}\nonumber\\+\Vert s\tau\varphi^\frac 32  (\mbox{\bf v}((-1)^k T,\cdot) ),\partial_{x_0}\mbox{\bf v}((-1)^k T,\cdot) )e^{s\varphi((-1)^k T,\cdot)}\Vert^2_{L^2(\partial\Omega)}).\nonumber
\end{eqnarray}
Applying the proposition  \ref{oxX} we have
\begin{eqnarray}\label{zmpolit11}
\int_{\Sigma_\epsilon}(s\tau\varphi\vert \nabla \mbox{\bf  v}\vert^2 +s^3\tau^3\varphi
^3\vert  \mbox{\bf  v}\vert^2 )e^{2s\varphi}d\Sigma + \int_{Q_\epsilon}(
s\tau^2\varphi\vert
\nabla\mbox{\bf v}\vert^2+ s^3\tau^4\varphi^3\vert
\mbox{\bf v}\vert^2 )e^{2s\varphi}dx\\
\le C_{22}\left(s\tau\Vert\root\of{\varphi}
\mbox{\bf g}e^{s\varphi}\Vert^2_{L^2(\Sigma)}+\Vert \mbox{\bf q}
e^{s\varphi}\Vert^2_{L^2(Q)}
+\int_{\widetilde \Sigma}(s\tau\varphi\vert \nabla \mbox{\bf  v}\vert^2
+s^3\tau^2\varphi^3\vert  \mbox{\bf  v}\vert^2 )e^{2s\varphi}d\Sigma
\right).\nonumber
\end{eqnarray}
Next we prove the following:
\begin{proposition}\label{govno} Let positive $\delta$ be sufficiently small.
There exist $\tau_0>0$ and $s_0>0$ such that for all $s\ge s_0$ and
$\tau\ge\tau_0$  there exists a constant $C_{23}>0$ independent of $s$ and $\tau$
such that
\begin{eqnarray}\label{govno1}
s\tau\Vert\varphi^\frac 12 \nabla\partial_{\vec \nu} \text{\bf u}
e^{s\varphi}\Vert_{{ L}^2( \Sigma_0)}^2
+s^3\tau^3\Vert\varphi^\frac 32 \partial_{\vec \nu}
\text{\bf u}
e^{s\varphi}\Vert_{L^2(\Sigma_0)}^2\nonumber\\\le C_{23} (s\tau\Vert\varphi
^\frac 12 \nabla \mbox{\bf v}
e^{s\varphi}\Vert_{{ L}^2(\Sigma_0)}^2
+s^3\tau^3\Vert\varphi^\frac 32  \mbox{\bf v}
e^{s\varphi}\Vert_{{ L}^2(\Sigma_0)}^2+s\tau\Vert \root\of{\varphi}\F e^{s\varphi}\Vert_{L^2(\Sigma_0)}^2)
\end{eqnarray}
and

\begin{equation}\label{pipi}
\Vert  s^\frac 12\tau\varphi^
\frac 12 \nabla d\omega_{\mbox{\bf u}}\,
e^{s\varphi}\Vert^2_{L^2(Q)} +\Vert ( s^\frac 12\tau\varphi^
\frac 12\nabla\mbox{div}\,\mbox{\bf u})
e^{s\varphi}\Vert^2_{L^2(Q)}\le
C_{23}(\int_{Q}(
s\tau^2\varphi\vert
\nabla\mbox{\bf v}\vert^2+ s^3\tau^4\varphi^3\vert
\mbox{\bf v}\vert^2 )e^{2s\varphi}dx).
\end{equation}
\end{proposition}

{\bf Proof.} The inequality (\ref{pipi}) follows from (\ref{qupol}), (\ref{pinok22}) and (\ref{victory}).  By (\ref{qupol}) there exists a constant $C_{24}$ independent of $s$ and $\tau$ such that
\begin{eqnarray}\label{zina1}
s\tau\Vert\varphi
^\frac 12 \nabla (d \omega_{\mbox {\bf u}},\mbox{div}\,\mbox {\bf u})
e^{s\varphi}\Vert_{{ L}^2(\Sigma_0)}^2
+s^3\tau^3\Vert\varphi^\frac 32  (d \omega_{\mbox {\bf u}},\mbox{div}\,\mbox {\bf u})
e^{s\varphi}\Vert_{{ L}^2(\Sigma_0)}^2\nonumber\\
\le C_{24}(
s\tau\Vert\varphi
^\frac 12 \nabla \mbox{\bf v}
e^{s\varphi}\Vert_{{ L}^2(\Sigma_0)}^2
+s^3\tau^3\Vert\varphi^\frac 32  \mbox{\bf v}
e^{s\varphi}\Vert_{{ L}^2(\Sigma_0)}^2).
\end{eqnarray}
Thanks to the zero Dirichlet boundary conditions on $ \Sigma_0$ there exists a smooth matrix $\mathcal M(x)$ such that
$\partial_{\vec \nu} \text{\bf u}=\mathcal M(x)(d \omega_{\mbox {\bf u}},\mbox{div}\,\mbox {\bf u}).$
Therefore
\begin{eqnarray}\label{zina}
s^3\tau^3\Vert\varphi^\frac 32 \partial_{\vec \nu}
\text{\bf u}
e^{s\varphi}\Vert_{L^2(\Sigma_0)}^2+s\tau\Vert\varphi^\frac 12 \partial_{\vec \nu}
\text{\bf u}
e^{s\varphi}\Vert_{H^1(\Sigma_0)}^2\\\le C_{25}(
s\tau\Vert\varphi
^\frac 12 \nabla (d \omega_{\mbox {\bf u}},\mbox{div}\,\mbox {\bf u})
e^{s\varphi}\Vert_{{ L}^2(\Sigma_0)}^2
+s^3\tau^3\Vert\varphi^\frac 32  (d \omega_{\mbox {\bf u}},\mbox{div}\,\mbox {\bf u})
e^{s\varphi}\Vert_{{ L}^2(\Sigma_0)}^2).\nonumber
\end{eqnarray}
From equation (\ref{gora1A}) on $\tilde \Sigma$ we have
\begin{equation}\label{eb1}
\partial_{\vec \nu}^2  \text{\bf u}=A(x,\widetilde D)\partial_{\vec \nu}
\text{\bf u} +B(x)(\text{\bf F}+\int_0^{x_0}\widetilde L_{\widetilde \lambda,
\widetilde\mu}(x, \widetilde x_0, D')\mbox{\bf u}(\widetilde x_0, x')
d\widetilde x_0),
\end{equation}
where $A(x,\widetilde D)$ is a first order differential operator and $B$ is
a $C^1-$ matrix function. From (\ref{eb1}) and (\ref{zina}) we have
\begin{eqnarray}\label{ep}
s\tau\Vert\varphi^\frac 12 \partial^2_{\vec \nu} \text{\bf u}
e^{s\varphi}\Vert_{L^2(\Sigma_0)}^2
 \le C_{26}(s\tau \Vert \varphi^\frac 12 A(x,\widetilde D)\partial_{\vec \nu}
\text{\bf u}\Vert_{L^2(\Sigma_0)}^2+s\tau\Vert \root\of{\varphi}\F e^{s\varphi}\Vert_{L^2(\Sigma_0)}^2\nonumber\\
 +s\tau\Vert \varphi^\frac 12B(x)\int_0^{x_0}\widetilde L_{\widetilde \lambda,
\widetilde\mu}(x, \widetilde x_0, D')\mbox{\bf u}(\widetilde x_0, x')
d\widetilde x_0)\Vert^2_{L^2(\Sigma_0)})\nonumber\\ \le C_{27}(s\tau\Vert\varphi
^\frac 12 \nabla (d \omega_{\mbox {\bf u}},\mbox{div}\,\mbox {\bf u})
e^{s\varphi}\Vert_{{ L}^2(\Sigma_0)}^2
+s^3\tau^3\Vert\varphi^\frac 32  (d \omega_{\mbox {\bf u}},\mbox{div}\,\mbox {\bf u})
e^{s\varphi}\Vert_{{ L}^2(\Sigma_0)}^2+\nonumber\\
\tau\Vert\varphi^\frac 12 \nabla\partial_{\vec \nu} \text{\bf u}
e^{s\varphi}\Vert_{{ L}^2( \Sigma_0)}^2
+s^2\tau^3\Vert\varphi^\frac 32 \partial_{\vec \nu}
\text{\bf u}
e^{s\varphi}\Vert_{L^2(\Sigma_0)}^2).
\end{eqnarray} In order to get the last inequality we used (\ref{qupol}).
Form (\ref{ep}), (\ref{zina1}) and (\ref{zina}) we obtain (\ref{govno1}).
The proof of Proposition \ref{govno} is complete. $\blacksquare$

Next we prove

\begin{proposition}\label{golos}
Let $\mbox{\bf u}\in H^1(Q), \mbox{\bf u}\vert_{\partial Q}
=0.$
There exist $\widehat\tau>1$  and $s_0>1$  such that
for any $\tau> \widehat\tau$ and for all
$s\ge s_0$
\begin{eqnarray}\label{3.2'}
\int_Q( s^2\tau^4\varphi^2\vert\nabla \mbox{\bf u}\vert^2
+ s^4\tau^6\varphi^4\vert \mbox{\bf u}\vert^2)e^{2s\varphi}dx
\le C_{28}(\Vert  s^\frac 12\tau\varphi^
\frac 12 \nabla d\omega_{\mbox{\bf u}}\,
e^{s\varphi}\Vert^2_{L^2(Q)} \\\nonumber+\Vert ( s^\frac 12\tau\varphi^
\frac 12\nabla\mbox{div}\,\mbox{\bf u})
e^{s\varphi}\Vert^2_{L^2(Q)}
+s^2\tau^3\Vert\varphi^\frac 12 \partial_{x_0}\partial_{\vec \nu} \text{\bf u}
e^{s\varphi}\Vert_{ L^2(\Sigma)}^2
+ s^2\tau^3\Vert \varphi\partial_{\vec \nu}  \text{\bf u}e^{s\varphi}
\Vert^2_{{ L}^2(\Sigma)}),
\end{eqnarray}
where $C_{28}$ is independent of $s$ and $\tau.$
\end{proposition}

{\bf Proof.} Let $x_0\in (-T,T)$ be arbitrarily fixed.
 For any index $\widehat i\in\{1,\dots, n\}$ consider the equality (\ref{gavnoed}).
Then the Carleman estimate with boundary  for the Laplace operator implies
\begin{eqnarray}\label{R}
\int_Q(s^2\tau^4\varphi^2 \vert \nabla '\mbox{\bf u}\vert^2+s^4\tau^6\varphi^4
\vert \mbox{\bf u}\vert^2)e^{2s\varphi}dx\le C_{29}(\Vert s^\frac 12\tau\varphi^
\frac 12\nabla ' \mbox{div}\, \mbox{\bf u}e^{s\varphi}\Vert^2_{L^2(Q)}
\nonumber\\
+\Vert s^\frac 12\tau\varphi^\frac 12 \nabla' d\omega_{\mbox{\bf u}}\,
e^{s\varphi}\Vert^2_{L^2(Q)}+\int_{\widetilde \Sigma}s^2\tau^3\varphi^2\vert
\partial_{\vec\nu}\mbox{\bf u}\vert^2 e^{2s\varphi}d\Sigma).
\end{eqnarray}
We differentiate both sides of equation (\ref{gavnoed})
with respect to the variable $x_0$ and take the $H^{-1}$-Carleman estimate
by authors in \cite{IY7}:

\begin{eqnarray}\label{RR}
\int_Q s^2\tau^4\varphi^2\vert \partial_{x_0} \mbox{\bf u}\vert^2e^{2s\varphi}
dx\le C_{30}(\Vert s^\frac 12\tau\varphi^\frac 12 \mbox{div}\, \partial_{x_0}
\mbox{\bf u}e^{s\varphi}\Vert^2_{L^2(Q)}\nonumber\\+\Vert s^\frac 12\tau
\varphi^\frac 12  \partial_{x_0} d\omega_{\mbox{\bf u}}\, e^{s\varphi}\Vert^2
_{L^2(Q)}+\int_{\widetilde \Sigma}s^2\tau^3\varphi^2\vert \partial
_{\vec\nu}\partial_{x_0}\mbox{\bf u}\vert^2 e^{2s\varphi}d\Sigma).
\end{eqnarray}
Combination of (\ref{R}) and (\ref{RR}) implies (\ref{3.2'}).
The proof of the proposition is complete.
                                    $\blacksquare$

By (\ref{pipi}),
Proposition \ref{govno} and Proposition \ref{golos}  from (\ref{zmpolit11}) we obtain the estimate
\begin{eqnarray}\label{zmpolit1o}
\int_{\Sigma}(s\tau\varphi\vert \nabla \mbox{\bf  u}\vert^2 +s^3\tau^2\varphi^3
\vert \mbox{\bf  u}\vert^2 )e^{2s\varphi}d\Sigma \nonumber\\
+ \int_{Q}( s\varphi\tau^2\vert\nabla d\omega_{\mbox{\bf u}}\vert^2
+(s\varphi)^3\tau^4\vert d\omega_{\mbox{\bf u}}\vert^2+ s\varphi\tau^2\vert
\nabla \mbox{div}\,\text{\bf u}\vert^2
+ (s\varphi)^3\tau^4\vert\mbox{div}\,\text{\bf u}\vert^2)e^{2s\varphi} dx
                                                   \nonumber\\
\le C_{31}\biggl(s\tau\Vert\root\of{\varphi}
\mbox{\bf F}e^{s\varphi}\Vert^2_{L^2({\Sigma})}+\Vert \mbox{\bf q} e^{s\varphi}
\Vert^2_{L^2( Q)}\\
+ \int_{\widetilde \Sigma}(s\tau\varphi\vert \nabla
d\omega_{\mbox{\bf u}}
\vert^2 +s^3\tau^2\varphi^3\vert  d\omega_{\mbox{\bf u}}\vert^2+s\tau\varphi
\vert \nabla  \mbox{div}\,{\mbox{\bf u}}\vert^2 +s^3\tau^2\varphi^3\vert
\mbox{div}\,{\mbox{\bf u}}\vert^2 )e^{2s\varphi}d x\biggr)           \nonumber
\end{eqnarray}
for all $s\ge s_0$ and for all $\tau\ge \tau_0$. By (\ref{nina}) and
(\ref{zmpolit1o}) we obtain the estimate (\ref{2.9'}).
Thus the proof of Theorem \ref{opa3} is finished. $\blacksquare$
\section{Proof of Theorem \ref{theorem 1.2}.}\label{Q1}

Denote $Q_\pm=(0,\pm T)\times \Omega.$
We  extend functions $\mbox{\bf y}(x),\F(x), \tilde \lambda(x), \tilde \mu(x)$ on $Q_-$ as $y(x)=y(-x_0,x'),\F(x)=\F(-x_0,x'), \tilde \lambda(x)=\tilde \lambda(-x_0,x'), \tilde \mu(x)=\tilde \mu(-x_0,x')$ for $Q_-.$
By (\ref{(1.18)}) the functions $\F,\partial_{x_0} \F$ belong to
$C([-T,T];H^1(\Omega))$  and $L^\infty (-T,T;H^1(\Omega))$ respectively and $\tilde \lambda,\tilde \mu\in C^1(\bar Q_+\times [0,T])\cap C^1(\bar Q_-\times [-T,0]).$
We set
$$
a_{\lambda,\mu}( \mbox{\bf z}(x_0,\cdot), \mbox{\bf v}(x_0,\cdot))=
\int_{\OOO} (\la(x')( \ddd \mbox{\bf z}(x),\ddd \mbox{\bf v}(x))
+ 2\mu(x') \sum_{\ell,j=1}^n \ep_{\ell j}(\mbox{\bf z}(x))\ep_{\ell j}(\mbox{\bf v}(x))) dx',
$$
$$
\ep_{\ell j}(\mbox{\bf y}) = \frac{1}{2}(\ppp_{x_\ell }y_j + \ppp_{x_j}y_\ell ), \quad
1\le \ell ,j \le n .
$$

Denote
$$
E_k({x_0})
=\int_{\OOO} (\la(x')\vert \partial_{x_0}^k\ddd \mbox{\bf y}(x)\vert^2
+ 2\mu(x') \sum_{\ell,j=1}^n \vert \ep_{\ell j}(\partial^k_{x_0}\mbox{\bf y}(x))\vert^2
+ \rho(x')\vert \ppp^{k+1}_{x_0}\mbox{\bf y}(x)\vert^2) dx'.
$$

The identity (6.194) established in \cite{DL1} page  612  for solution of the system (\ref{(1.1)}), (\ref{(1.222)}) yields

$$
\frac{d}{dx_0}E_k(x_0)=(\partial^k_{x_0}\int^{x_0}_0 \LLLL(x',\widetilde x_0, D')
\mbox{\bf y}({\widetilde x_0},x') d{\widetilde x_0}
+ \partial^{k}_{x_0}\F, \partial^{k+1}_{x_0} \y)_{L^2(\Omega)}.
$$
Integrating the right-hand side of this equality  for any $k\in\{0,1,2\}$ we have
$$
\frac{d}{dx_0}E_k(x_0)= (\partial^k_{x_0}\F, \partial^{k+1}_{x_0} \y)_{L^2(\Omega)}-\frac{d}{dx_0}\int_0^{x_0}  a_{\partial^k_{x_0}\tilde \lambda,\partial^k_{x_0}\tilde \mu}(\y(\tilde x_0,\cdot) ,\partial^k_{x_0}\y(x_0,\cdot))d\tilde x_0
$$
$$
-\frac{d}{dx_0}\sum_{p=0}^{k-1}C_p a_{\partial^{k-p-1}_{x_0}\tilde \lambda(x,x_0),\partial^{k-p}_{x_0}\tilde \mu(x,x_0)}(\partial^p_{x_0}\y( x_0,\cdot) ,\partial^k_{x_0}\y(x_0,\cdot))
$$
$$
+ \frac{d}{dx_0}\int_0^{x_0}  a_{\partial^k_{x_0}\tilde \lambda,\partial^k_{x_0}\tilde \mu}(\y(\tilde x_0,\cdot) ,\partial^k_{t}\y(t,\cdot))\vert_{t=x_0}d\tilde x_0
$$

$$
+\frac{d}{dx_0}\sum_{p=0}^{k-1}C_p a_{\partial^{k-p-1}_{x_0}\tilde \lambda(x,x_0),\partial^{k-p}_{x_0}\tilde \mu(x,x_0)}(\partial^p_{x_0}\y(\tilde x_0,\cdot) ,\partial^k_{t}\y(t,\cdot))\vert_{t=x_0},
$$
 where $C_p$ are some constants.
Integrating the above equality on  the interval $(0,x_0)$ we obtain
$$
E_k(x_0)\le E_k(0)+\Vert\partial^k_{x_0} \F\Vert^2_{L^2(Q_+)}+C_1\int_0^{x_0}E_k(\tilde x_0)d\tilde x_0-\int_0^{x_0}  a_{\partial^k_{x_0}\tilde \lambda,\partial^k_{x_0}\tilde \mu}(\y(\tilde x_0,\cdot) ,\partial^k_{x_0}\y(x_0,\cdot))d\tilde x_0
$$
$$
-\sum_{p=0}^{k-1}C_p a_{\partial^{k-p}_{x_0}\tilde \lambda(x,x_0),\partial^{k-p}_{x_0}\tilde \mu(x,x_0)}(\partial^p_{x_0}\y(x_0,\cdot) ,\partial^k_{x_0}\y(x_0,\cdot))
$$
$$+
\sum_{p=0}^{k-1}C_p a_{\partial^{k-p}_{x_0}\tilde \lambda(x,x_0),\partial^{k-p}_{x_0}\tilde \mu(x,x_0)}(\partial^p_{x_0}\y(x_0,\cdot) ,\partial^k_{x_0}\y(x_0,\cdot))\vert_{t=x_0}
$$
$$
+\int_0^{x_0}\frac{d}{dy_0} \int_0^{y_0}  a_{\partial^k_{y_0}\tilde \lambda(y_0,x',\tilde x_0),\partial^k_{y_0}\tilde \mu(y_0,x',\tilde x_0)}(\y(\tilde x_0,\cdot) ,\partial^k_{t}\y(t,\cdot))\vert_{t=y_0}d\tilde x_0dy_0
$$

$$
-\int_0^{x_0}\frac{d}{dy_0}\sum_{p=0}^{k-1}C_{p} a_{\partial^{k-p}_{y_0}\tilde \lambda(y_0,x',y_0),\partial^{k-p}_{y_0}\tilde \mu(y_0,x',y_0)}(\partial^p_{x_0}\y( y_0,\cdot) ,\partial^k_{t}\y(t,\cdot))\vert_{t=y_0}dy_0
$$
$$
\le  E_k(0)+\Vert\partial^k_{x_0} \F\Vert^2_{L^2(Q_+)}+C_2\int_0^{x_0}E_k(\tilde x_0)d\tilde x_0 +\frac 12 E_k(x_0).
$$
Thus the Gronwall inequality and the Korn inequality yield
\begin{lemma} \label{3}{\it
There exists a constant $C_3>0$ such that for any $x_0$ from $[0,T]$
\begin{equation}
\sum_{k=0}^2\Vert \nabla \partial^k_{x_0}\mbox{\bf y}(x_0,\cdot)\Vert^2_{L^2(\OOO)}
\le C_3 (\sum_{k=0}^2
 \Vert\partial^k_{x_0} \F\Vert^2_{L^2(Q_+)}+\Vert \F(0,\cdot)\Vert^2_{H^1_0(\Omega)}+\Vert \partial_{x_0}\F(0,\cdot)\Vert^2_{L^2(\Omega)}),
\end{equation}
for $\mbox{\bf y}$ satisfying (\ref{(1.1)})- (\ref{(1.222)}).}
\end{lemma}
From (\ref{(1.1)}), Lemma \ref{3} and the classical apriori estimates
for the stationary Lam\'e system, it follows that
\begin{equation}\label{pop}
\Vert \y\Vert_{L^2(0,T;H^2(\Omega))}+\Vert\partial_{x_0}\y\Vert_{L^2(0,T;H^2(\Omega))}\le C_4(\sum_{k=0}^2\Vert \nabla \partial^k_{x_0}\mbox{\bf y}\Vert_{L^2(Q_+)}+\Vert \F\Vert_{L^2(Q_+)}+\Vert \partial_{x_0}\F\Vert_{L^2(Q_+)}).
\end{equation}
%
%

Observe that by (\ref{logoped}) there exists a positive $\delta_1$ such that
\begin{equation}\label{oop}
\mbox{inf}_{x'\in \Omega}\varphi(0,x')>\mbox{sup}_{x\in ([T-\delta_1,T]\cup [-T,-T+\delta_1])\times \Omega}\varphi(x).
\end{equation}
Let $\tilde\gamma\in C^\infty_0[-T,T]$ satisfy
$\tilde\gamma\vert_{\vert t\vert\le T-\delta_1/2}=1$, and the
function $\mbox{\bf v}$ be given by  (\ref{victory}).
Then, taking the even extension of the function  $\mbox{\bf v}$ on  $Q_-$ we have

\begin{equation}\label{poko!}
\mbox{\bf P}(x,D) \mbox{\bf v}=(
\square_{\rho,\mu} (x,D)\mbox{\bf v}_1, \square_{\rho,\lambda+2\mu} (x,D) v_2)=\mbox{\bf q} \quad \mbox{in}\,\, Q, \quad \mathcal B(x',D)\mbox{\bf v}=\mbox{\bf g}.
\end{equation}
We claim $\partial_{x_0}\mbox{\bf q}\in L^2(Q).$ Indeed, since
$\F,\partial_{x_0} \F \in L^2(-T,T;H^1(\Omega))$, we have
$\partial_{x_0}d\omega_{\mbox{\bf F}},$ $\partial_{x_0}\mbox{div}\,\mbox{\bf F}\in L^2(Q).$
By (\ref{(1.222)}), the regularity  assumption on $\mathbf{y}$ and
time independence of the Lam\'e coefficients $\lambda,\mu$, we see that the
$x_0$-derivative of the extensions of the terms $\beta\partial_{x_j}y_i,
\beta\partial_{x_jx_k}y_i$ with $\beta\in\{\lambda, \mu, \nabla \mu,\nabla\lambda\}$ belong to $L^2(Q).$
By the same argument the  $x_0$- derivative of the extension of the functions $\int_0^{x_0}\beta\partial_{x_j}y_id\tilde x_0,\int_0^{x_0}\beta\partial_{x_jx_k}y_id\tilde x_0$ with $\beta\in\{\lambda, \mu, \nabla \mu,\nabla\lambda\}$ belongs to $L^2(Q).$ Hence our claim follows from (\ref{nokia}).

Setting $\widetilde {\mbox{\bf v}}=\tilde \gamma {\mbox{\bf v}}$ we have
\begin{equation}\label{poko!!}
\mbox{\bf P}(x,D) \widetilde{\mbox{\bf v}}=\widetilde {\mbox{\bf q}} + [\tilde\gamma, \mbox{\bf P}(x,D)]\mbox{\bf v} \quad \mbox{in}\,\, Q,  \quad \mathcal B(x',D)\mbox{\bf v}=\tilde\gamma\mbox{\bf g},\quad \widetilde{\mbox{\bf v}}(0,\cdot)=\partial_{x_0}\widetilde{\mbox{\bf v}}(0,\cdot)=0.
\end{equation}
Taking the time derivative of (\ref{poko!!}), we have
\begin{eqnarray}\label{poko!!l}
\mbox{\bf P}(x,D) \partial_{x_0}\widetilde{\mbox{\bf v}}=\partial_{x_0}\widetilde {\mbox{\bf q}} + \partial_{x_0}[\tilde\gamma, \mbox{\bf P}(x,D)]\mbox{\bf v} \quad \mbox{in}\,\, Q,\\ \quad \mathcal B(x',D)\mbox{\bf v}=\partial_{x_0}(\tilde\gamma\mbox{\bf g}),\quad  \partial^2_{x_0}\widetilde{\mbox{\bf v}}(0,\cdot)=(\rho d\omega_{\mbox{\bf F}/\rho}(0,\cdot),\mbox{div}\,(\mbox{\bf F}/\rho)(0,\cdot)).\nonumber
\end{eqnarray}
Since $\mbox{\bf f}\in H^1_0(\Omega)$, we have
\begin{equation}\label{gorbun}
\vert \mbox {\bf g}(x)\vert \le C\vert\nabla ' \mbox{\bf y}(x)\vert \quad {and}\quad \vert \partial_{x_0}\mbox {\bf g}(x)\vert \le C\vert\nabla ' \partial_{x_0}\mbox{\bf y}(x)\vert \quad \mbox{on}\,\,\Sigma.
\end{equation}
We claim that $\partial_{x_0}\widetilde{\mbox{\bf v}}\in H^1(Q).$
By (\ref{(1.222)}), the even extensions of the functions $d\omega_{\mbox{\bf y}},\mbox{div}\,\mbox{\bf y}$ on $Q_-$ belongs to $H^2(Q).$ Set $\tilde {\mbox{\bf w}}={\mbox{\bf w}}_+$ for $x_0>0$ and $\tilde{ \mbox{\bf w}}=\mbox{\bf w}_-$ for $x_0<0,$  $\mbox{\bf w}_\pm=(\mbox{\bf w}_{1,\pm},w_{2,\pm})$, where
 $$
\mbox{\bf w}_{1,\pm}=\pm\int_0^{\pm x_0}\frac{\widetilde
\mu(\pm x_0,x',\pm\widetilde x_0)}{\mu(x')} d\omega_{\mbox{\bf y}} \,d\widetilde x_0,
\quad w_{2,\pm}=\pm\int_0^{\pm x_0}\frac{(\widetilde \lambda
+2\widetilde \mu)(\pm x_0,x',\pm\widetilde x_0)}{(\lambda+2\mu)(x')}
\text{div}\,
\mbox{\bf y}\, d\widetilde x_0.
$$
By our regularity assumptions $\mbox{\bf w}_\pm, \nabla'\mbox{ \bf w}_\pm\in H^1(Q_\pm)$ and $\mbox{\bf w}_\pm(0,\cdot)=\nabla'\mbox{\bf w}_\pm(0,\cdot)=0.$
Therefore $\tilde{\mbox{\bf w}},  \nabla'\tilde {\mbox{\bf w}}\in H^1(Q).$ Note  that
$$
\partial_{x_0}\mbox{\bf w}_\pm=\pm(\frac{\widetilde
\mu(\pm x_0,x',\pm x_0)}{\mu(x')} d\omega_{\mbox{\bf y}},\frac{(\widetilde \lambda
+2\widetilde \mu)(\pm x_0,x',\pm x_0)}{(\lambda+2\mu)(x')} \text{div}\,
\mbox{\bf y})\pm
$$
$$
(\int_0^{\pm x_0} \frac{\partial}{\partial x_0}\left [\frac{\widetilde
\mu(\pm x_0,x',\pm\widetilde x_0)}{\mu(x')}\right ] d\omega_{\mbox{\bf y}} \,d\widetilde x_0,
\int_0^{\pm x_0}\frac{\partial}{\partial x_0}\left [\frac{(\widetilde \lambda
+2\widetilde \mu)(\pm x_0,x',\pm\widetilde x_0)}{(\lambda+2\mu)(x')}\right] \text{div}\,
\mbox{\bf y}\, d\widetilde x_0).
$$
Observe that by (\ref{(1.222)}) $\partial_{x_0}\mbox{\bf w}_\pm(0,\cdot)=\nabla '\partial_{x_0}\mbox{\bf w}_\pm(0,\cdot)=0$
and by our regularity assumptions on function $\mbox{\bf y}$ the functions  $\partial_{x_0}\mbox{\bf w}_\pm, \nabla'\partial_{x_0}\mbox{ \bf w}_\pm$ belong to  the Sobolev space $H^1(Q_\pm).$
Therefore $\partial_{x_0}\mbox{\bf w}\in H^1(Q).$
Then $\tilde{\mbox{\bf v}},\partial_{x_0}\tilde{\mbox{\bf v}} \in H^1(Q)$
by (\ref{victory}).
Then applying to the problems (\ref{poko!!l}) and (\ref{poko!!}) the Carleman estimate (\ref{3.2'1})  and using (\ref{gorbun}), we obtain
\begin{eqnarray}\label{Car}
\sum_{k=0}^1(
s\tau^2\Vert\varphi^\frac 12  \nabla \partial^k_{x_0}\widetilde{\mbox{\bf v}}\,
e^{s\varphi}\Vert^2_{L^2(Q)} +s^3\tau^4\Vert \varphi^\frac 32  \partial^k_{x_0}\widetilde{ \mbox{\bf v}}\,
e^{s\varphi}\Vert^2_{L^2(Q)} \nonumber\\
+\int_{ \Sigma}(s\tau\varphi \vert \nabla\partial^k_{x_0} \widetilde{ \mbox{\bf v}}\vert^2+s^3\tau^3\varphi^3\vert\partial^k_{x_0}\widetilde{ \mbox{\bf v}}\vert^2) e^{2s\varphi}d\Sigma)\nonumber\\
\le C_5\sum_{k=0}^1(\Vert \mbox{\bf P}(x,D)\partial^k_{x_0}\widetilde{\mbox{\bf v}} e^{s\varphi}\Vert^2_{L^2(Q)}+\nonumber\\
+\int_{\tilde \Sigma}(s\tau\varphi \vert \nabla  \partial^k_{x_0}\widetilde{\mbox{\bf v}}\vert^2+s^3\tau^3\varphi^3\vert \partial^k_{x_0}\widetilde{\mbox{\bf v}}\vert^2) e^{2s\varphi}d\Sigma).
\end{eqnarray}
Taking the scalar  product of  (\ref{poko!!l})  with $\partial^2_{x_0}\widetilde {\mbox{\bf v}} e^{2s\varphi}$ in $L^2(Q_+)$  and integrating by parts, using the inequality (\ref{Car}) we obtain the estimate
\begin{eqnarray}\label{P2}
\tau\int_{\Omega} \vert \partial^2_{x_0}\widetilde{\mbox{\bf v}}(0,\cdot)\vert^2 e^{2s\varphi(0,\cdot)}dx'
\le
 C_6(\Vert \mbox{\bf P}(x,D)\partial_{x_0}\widetilde{\mbox{\bf v}} e^{s\varphi}\Vert^2_{L^2(Q)}\nonumber\\
+\int_{\tilde \Sigma}(s\tau\varphi \vert \nabla  \partial_{x_0}\widetilde{\mbox{\bf v}}\vert^2+s^3\tau^3\varphi^3\vert \partial_{x_0}\widetilde{\mbox{\bf v}}\vert^2) e^{2s\varphi}d\Sigma).
\end{eqnarray}

By assumption (\ref{giga1}), $\nabla '\psi(0,x')$ is not equal to zero on $\bar \Omega.$
By Proposition \ref{golos}, (\ref{legion})  and (\ref{(1.18)}) we obtain
\begin{eqnarray}\label{P3}
\frac{\tau}{s}\int_\Omega\frac{1}{\varphi(0,x')}  \vert \nabla' \mbox{\bf f}\vert^2 e^{2s\varphi(0,x')}dx'\nonumber\\ \le C_7(
\tau\int_{\Omega} \vert \partial^2_{x_0}\widetilde{\mbox{\bf v}}(0,\cdot)\vert^2 e^{2s\varphi(0,\cdot)}dx'+\int_{ \Gamma} s\tau^2\varphi(0,x')\vert \partial_{\vec \nu}\partial^2_{x_0}\mbox{\bf y}(0,x')\vert^2 e^{2s\varphi(0,x')}d\sigma).
\end{eqnarray}

Using (\ref{(1.18)}) we have
\begin{eqnarray}\label{gugenot}\sum_{k=0}^1\Vert \partial^k_{x_0}(d\omega_{\mbox{\bf F}},\mbox{div}\,\mbox{\bf F})e^{s\varphi}\Vert_{L^2(Q)}
\\\le  \sum_{k=0}^1\Vert \partial^k_{x_0} (d\omega_{\mbox{\bf F}-\mbox{\bf F}(0,\cdot)},\mbox{div}\,(\mbox{\bf F}-\mbox{\bf F}(0,\cdot)))e^{s\varphi}\Vert_{L^2(Q)}+\sum_{k=0}^1\Vert \partial^k_{x_0}(d\omega_{\mbox{\bf F}}(0,\cdot),\mbox{div}\,\mbox{\bf F}(0,\cdot))e^{s\varphi}\Vert_{L^2(Q)}\nonumber\\=
\sum_{k=0}^1\Vert \partial^k_{x_0} (d\omega_{\mbox{\bf F}-\mbox{\bf F}(0,\cdot)},\mbox{div}\,(\mbox{\bf F}-\mbox{\bf F}(0,\cdot)))e^{s\varphi}\Vert_{L^2(Q)}+\Vert (d\omega_{\mbox{\bf F}}(0,\cdot),\mbox{div}\,\mbox{\bf F}(0,\cdot))e^{s\varphi}\Vert_{L^2(Q)}.\nonumber\end{eqnarray}

By (\ref{pinok22}) and assumption $\partial^2_{x_0}\psi(0,x')\ne 0$ on
$\bar\Omega$, we obtain
\begin{eqnarray}
\Vert (d\omega_{\mbox{\bf F}}(0,\cdot),\mbox{div}\,\mbox{\bf F}(0,\cdot))e^{s\varphi}\Vert_{L^2(Q)}=O(\frac{1}{s^{\frac 14}})\Vert (d\omega_{\mbox{\bf F}}(0,\cdot),\mbox{div}\,\mbox{\bf F}(0,\cdot))e^{s\varphi(0,\cdot)}\Vert_{L^2(\Omega)}\nonumber\\ =O(\frac{1}{s^{\frac 14}})\Vert \partial^2_{x_0}\widetilde{\mbox{\bf v}}(0,\cdot) e^{s\varphi(0,\cdot)}\Vert_{L^2(\Omega)}\quad \mbox{as}\,\,s\rightarrow +\infty.
\end{eqnarray}
Observe that by (\ref{(1.18)}) there exist functions $g_0,g_1\in C_0^1[-T,T]$
such that $g_0(0)=g_1(0)=0$ and
\begin{eqnarray}\label{P1}
\sum_{k=0}^1\Vert \partial^k_{x_0} (d\omega_{\mbox{\bf F}-\mbox{\bf F}(0,\cdot)},\mbox{div}\,(\mbox{\bf F}-\mbox{\bf F}(0,\cdot)))e^{s\varphi}\Vert_{L^2(Q)}\\
\le C_8(\int_\Omega\int_{-T}^{T}( g^2_1(x_0)\vert \nabla '\mbox{\bf f}\vert^2 +g^2_0(x_0)\vert \mbox{\bf f}\vert^2 ) e^{2s\varphi(x)} dx_0dx')^\frac 12=\nonumber\\
C_9(\int_\Omega\int_{-T}^{T}( g^2_1(x_0)\vert \nabla '\mbox{\bf f}\vert^2 +g^2_0(x_0)\vert \mbox{\bf f}\vert^2 ) \partial_{x_0}e^{2s\varphi(x)}/(2s\varphi\tau \partial_{x_0} \psi) dx_0dx')^\frac 12=\nonumber\\
C_{10}(\int_\Omega\int_{-T}^{T}\left (\frac{d}{dx_0}\left( \frac{g^2_1(x_0)}{2s\tau\varphi\partial_{x_0} \psi}\right)\vert \nabla '\mbox{\bf f}\vert^2 +\frac{d}{dx_0}\left( \frac{g^2_0(x_0)}{2s\varphi\tau \partial_{x_0}\psi}\right)\vert \mbox{\bf f}\vert^2\right ) e^{2s\varphi(x)} dx_0dx')^\frac 12=o(\frac{ 1}{\root\of{s}}).\nonumber
\end{eqnarray}
Combining (\ref{Car}), (\ref{P2}) and (\ref{gugenot})-(\ref{P1}), we see
that there exist functions $c_j,\tilde c_j, c,\tilde c \in L^\infty(Q)\cap C^1(\bar Q_\pm),$  and  $s_0$ independent of $\mbox{\bf y}$ such that
for all sufficiently large $s$
\begin{eqnarray}\label{P3}
\frac{\tau}{s} \int_\Omega  \vert \nabla' \mbox{\bf f}\vert^2 e^{2s\mbox{inf}_{x'\in \Omega}\varphi(0,x')}dx'\le
\frac{\tau}{s} \int_\Omega  \vert \nabla' \mbox{\bf f}\vert^2 e^{2s\varphi(0,x')}dx'\le C_{11}\sum_{k=0}^1 (\int_{{\supp \tilde\gamma '}\times \Omega }( \vert\partial_{x_0}^k\nabla \mbox{\bf v}\vert^2+ \vert\partial_{x_0}^k \mbox{\bf v}\vert^2\nonumber\\
+\sum_{j=1}^n\vert\partial_{x_0}^k \int_{0}^{x_0} c_j(x,\tilde x_0)\partial_{x_j}\mbox{\bf v}(\tilde x_0,x')d\tilde x_0 \vert^2+ \vert\partial_{x_0}^k \int_{0}^{x_0} c(x,\tilde x_0)\mbox{\bf v}(\tilde x_0,x')d\tilde x_0\vert^2\nonumber\\
+\sum_{k=0}^1 \vert\partial_{x_0}^k \nabla \mbox{\bf y}\vert^2+ \vert\partial_{x_0}^k \mbox{\bf y}\vert^2
+\sum_{j=1}^n\vert\partial_{x_0}^k \int_{0}^{x_0}\tilde c_j(x,\tilde x_0)\partial_{x_j}\mbox{\bf y}(\tilde x_0,x')d\tilde x_0 \vert^2\nonumber\\+ \vert\partial_{x_0}^k \int_{0}^{x_0} \tilde c(x,\tilde x_0)\mbox{\bf y}(\tilde x_0,x')d\tilde x_0\vert^2)e^{2s\sup_{x\in {\supp \tilde\gamma'}\times \Omega}\varphi}dx\nonumber\\
+\int_{ \Gamma} s\tau^2\varphi(0,x')\vert \partial_{\vec \nu}\partial^2_{x_0}\mbox{\bf y}(0,x')\vert^2 e^{2s\varphi(0,x')}d\sigma +\sum_{k=0}^2\int_{\tilde \Sigma}(s\tau\varphi \vert \nabla  \partial^k_{x_0}\mbox{\bf v}\vert^2+s^3\tau^3\varphi^3\vert \partial^k_{x_0}\mbox{\bf v}\vert^2) e^{2s\varphi}d\Sigma).\nonumber
\end{eqnarray}
By (\ref{oop}), we have
\begin{equation}\label{estimate}
\mbox{inf}_{x'\in \Omega}\varphi(0,x')>\mbox{sup}_{x\in{\supp \tilde\gamma '}\times \Omega}\varphi(x).
\end{equation}
Hence, by (\ref{estimate}), (\ref{pop}) and Lemma {\ref{3}} there exists $s_1$ such that
\begin{eqnarray}\label{P4}
\frac{\tau}{s} \int_\Omega  \vert \nabla' \mbox{\bf f}\vert^2 e^{2s\mbox{inf}_{x'\in \Omega}\varphi(0,x')}dx'
\le C_{12}(\sum_{k=0}^1 (\int_{{\supp \tilde\gamma '}\times \Omega }( \vert\partial_{x_0}^k\nabla \mbox{\bf v}\vert^2+ \vert\partial_{x_0}^k \mbox{\bf v}\vert^2\\
+\sum_{j=1}^n\vert\partial_{x_0}^k \int_{0}^{x_0} c_j(x,\tilde x_0)\partial_{x_j}\mbox{\bf v}(\tilde x_0,x')d\tilde x_0 \vert^2+ \vert\partial_{x_0}^k \int_{0}^{x_0} c(x,\tilde x_0)\mbox{\bf v}(\tilde x_0,x')d\tilde x_0\vert^2)e^{2s\sup_{x\in {\supp \tilde\gamma'}\times \Omega}\varphi}dx\nonumber\\
+\int_{\tilde \Sigma}(s\tau\varphi \vert \nabla  \partial^k_{x_0}\mbox{\bf v}\vert^2+s^3\tau^3\varphi^3\vert \partial^k_{x_0}\mbox{\bf v}\vert^2) e^{2s\varphi}d\Sigma)+\int_{ \Gamma} s\tau^2\varphi(0,x')\vert \partial_{\vec \nu}\partial^2_{x_0}\mbox{\bf y}(0,x')\vert^2 e^{2s\varphi(0,x')}d\sigma )\nonumber
\end{eqnarray}
for all $s\ge s_1$.
From (\ref{P4}), (\ref{pop}), (\ref{victory}) and (\ref{estimate}) we obtain
\begin{eqnarray}
\frac{\tau}{s} \int_\Omega  \vert \nabla' \mbox{\bf f}\vert^2 e^{2s\mbox{inf}_{x'\in \Omega}\varphi(0,x')}dx' \le C_{13}(\int_{ \Gamma} s\tau^2\varphi(0,x')\vert \partial_{\vec \nu}\partial^2_{x_0}\mbox{\bf y}(0,x')\vert^2 e^{2s\varphi(0,x')}d\sigma\nonumber\\
 +\sum_{k=0}^1\int_{\tilde \Sigma}(s\tau\varphi \vert \nabla  \partial^k_{x_0}\mbox{\bf v}\vert^2+s^3\tau^3\varphi^3\vert \partial^k_{x_0}\mbox{\bf v}\vert^2) e^{2s\varphi}d\Sigma).\nonumber
\end{eqnarray}
 The proof of Theorem \ref{theorem 1.2} is completed. $\blacksquare$

%
%
%
%
\section{Appendix}\label{Q!8}

In the appendix we prove several lemmata which are used for the proof of
Theorem \ref{opa3} and represent the standard  properties of the  pseudodifferential operators with symbols of limited smoothness, keeping the dependence of
the norms on parameter $\tau.$

Let $\{\omega_j\}_{j=1}^\infty$ be a sequence of the eigenfunctions of
the operator $-\Delta$ on $\Bbb S^{n
+2}$ and let $\{\lambda_j\}_{j=1}^\infty$ be a sequence of the corresponding
eigenvalues of $-\Delta$. Assume that
$$
(\omega_k,\omega_j)_{L^2(\Bbb S^{n+2})}=\delta_{kj}\quad \forall k,j\in
\Bbb N.
$$
The following asymptotic formula is known (e.g., Courant and Hilbert
\cite{CH}):
$$
\lambda_j=c j^\frac {2}{n} + o(j^{2/n})\quad\mbox{as}\,\,
j\rightarrow+\infty.
$$
For each $k$, thanks to the standard elliptic estimate for the
Laplace operator, we have
 \begin{equation}\label{Fzina}
\Vert \omega_j\Vert_{H^{2k}({\Bbb S^{n+2}})}\le C_k\lambda_j^k.
\end{equation}
Therefore the Sobolev embedding theorem yields
\begin{equation}\label{Fest1}
\Vert\omega_j\Vert_{C^0({\Bbb S^{n
+2}})}\le C \lambda_j^{n}, \quad j\in \Bbb N.
\end{equation}
We extend the function $\omega_j$ on the set $\{\vert\xi\vert\le
1\}$ as a smooth function and we set
$$
\omega_j(\xi)=\omega_j(\xi/\vert \xi\vert)\quad\mbox{for}\,\,\,
\vert\xi\vert\ge 1.
$$
We introduce the pseudodifferential operator
$$
\widetilde\omega_j(D)w=\int_{{\Bbb  R}^{n+1}}\omega_j(\xi) \widehat w(\xi)
e^{{i}<y,\xi>}d\xi,\quad \widehat w(\xi)=\frac{1}{(2\pi)^{\frac
{n+1}{2}}}\int_{{\Bbb  R}^{n+1}}w(y)e^{-{i}<\xi,y>}dy.
$$
Here we recall that, in order to distinguish the Fourier transforms
with respect to different variables, we will use the following notations
\begin{equation}\label{OO}
\widehat{u}(\www \xi) := F_{\widetilde y\rightarrow
\widetilde \xi}u=\frac{1}{(2\pi)^\frac{n}{2}}\int_{{\Bbb R}^{n}}
e^{-{i}\sum_{j=0}^{n-1}y_j\xi_j}u(y_0,\dots,y_{n-1})d\widetilde y,
\end{equation}
$$
\widehat{u}({ \xi}_n) :=
F_{y_n\rightarrow \xi_n}u
=\frac{1}{(2\pi)^\frac{1}{2}}\int_{{\Bbb  R}}e^{-{i}y_n\xi_n}u(y_{n})
dy_n.
$$
We recall that we set $\www\xi = (\xi_0, ...., \xi_{n-1})$ and
$\www y = (y_0, ..., y_{n-1})$.

First we define the operator $A(\widetilde y,{\widetilde D},s)$ for functions
in $C^\infty_0(\mathcal O):$
$$
A(\widetilde y,{\widetilde D},s)u=\int_{{\Bbb R}^{n-1}}a(\widetilde y,\xi',s)
F_{\widetilde y\rightarrow \widetilde\xi}u e^{{i}\sum_{j=0}
^{n-1}y_j\xi_j}d\widetilde y.
$$

For the symbol $a$, we introduce the semi-norm
\begin{eqnarray}
\pi_{C^k}(a)=\sum_{j=1}^{\widehat N}\sup_{\vert \beta\vert\le \widehat N}
\sup_{\vert(\widetilde \xi,s)\vert\ge 1}\left\Vert \frac{\partial^{\beta_0}}
{\partial\xi_0^{\beta_0}}\cdots
\frac{\partial^{\beta_{n-1}}}{\partial\xi_{n-1}^{\beta_{n-1}}}
\frac{\partial^{\beta_{n}}}{\partial s^{\beta_{n}}}
a_j(\cdot,\widetilde\xi,s)\right\Vert_{C^k(\bar{\mathcal
O})}/(1+\vert(\widetilde\xi,s)\vert)^{\kappa-j-\vert\beta\vert} \nonumber \\
+\sup_{\vert (\widetilde\xi,s)\vert\le 1}\Vert
a(\cdot,\widetilde\xi,s)\Vert_{C^\kappa(\bar{\mathcal O})}. \nonumber
\end{eqnarray}

The following lemma allows us to extend the definition of the
operator $A$ on Sobolev spaces.
\begin{lemma}\label{Fops0}
Let $a(\widetilde y,\widetilde\xi,s)\in C^0_{cl}S^{1,s}({\mathcal O}).$ Then
$A\in {\mathcal L}(H^{1,s}_0({\mathcal O});L^2({\mathcal O}))$ and $\Vert
A\Vert_{{\mathcal L}(H^{1,s}_0({\mathcal O});L^2({\mathcal
O}))}\le C(\pi_{C^0}(a)). $
\end{lemma}

\noindent {\bf Proof.} Thanks to the assumption {\bf C}, it suffices
to consider the case where
\begin{equation}\label{Fzvezda}
a(\widetilde y,{\tau}\widetilde \xi, {\tau} s)={\tau}
a(\widetilde y,\widetilde \xi,s)\quad \forall{\tau}>1.
\end{equation}
The operator
\begin{equation}
\widetilde A(\widetilde y,D)v=\int_{\{\vert\xi\vert\le 1\}}
a(\widetilde y,\xi)F_{ y\to \xi}v(\xi)e^{{i}<y,\xi>}d\xi
\nonumber
\end{equation}
is a continuous operator from $L^2({\mathcal
O}\times {\Bbb  R})$ into $L^2({\mathcal O}\times {\Bbb R})$ with the norm
estimated as
$$ \Vert\widetilde A(\widetilde y,D)\Vert_{\mathcal L(L^2({\mathcal O}
\times {\Bbb R}),L^2({\mathcal O}\times {\Bbb R}))}\le C(\pi_{C^0(\mathcal O)}
(a)).
$$
Consider the symbol
$b(\widetilde y,\xi)=a(\widetilde y,\xi)/\vert \xi\vert .$ Then (\ref{Fzvezda})
implies
$$
b(\widetilde y,\tau\xi)=b(\widetilde y,\xi)\quad \forall {\tau}\ge 1.
$$
We can represent the symbol $b$ as
$$
b(\widetilde y,\xi)=\sum_{j=1}^\infty
b_j(\widetilde y)\omega_j(\xi/\vert \xi\vert ),\quad b_j(\widetilde y)
=(b(\widetilde y,\xi),\omega_j(\xi))_{L^2({\Bbb S^{n+2}})}.
$$
Observe that $b_j(\widetilde y)=(\Delta^k_\xi
b(\widetilde y,\xi),\omega_j(\xi))_{L^2({\Bbb S^{n+2}})}/\lambda^k_j.$
Therefore
\begin{equation}\label{Fest2}
\Vert b_j\Vert_{C^0(\bar{\mathcal O})}\le C_m\lambda_j^{-m}\quad
\forall m\in\{1,\dots, \infty\}.
\end{equation}
By (\ref{Fest1}) and (\ref{Fest2}), we have
$$
\Vert B(\widetilde y,D)v\Vert_{L^2({\mathcal O}\times {\Bbb  R})}\le
\sum_{j=1}^\infty\Vert b_j\Vert_{C^0(\bar{\mathcal O})}\Vert \widetilde
\omega_j(D)\Vert_{\mathcal L(L^2({\mathcal O}\times {\Bbb  R}),
L^2({\mathcal O}\times {\Bbb  R}))} \Vert v\Vert_{L^2({\mathcal O}
\times {\Bbb  R})}
$$
$$\le
\sum_{j=1}^\infty C_m\lambda_j^{-m}\lambda_j^{n}\Vert
v\Vert_{L^2({\mathcal O}\times {\Bbb R})}.
$$
Taking $m=3n$ we have
$$
\Vert B(\widetilde y,D)v\Vert_{L^2({\mathcal O}\times {\Bbb R})}\le
\sum_{j=1}^\infty C_m\lambda_j^{-2n}\Vert v\Vert_{L^2({\mathcal
O}\times {\Bbb R})}.
$$
Therefore the operator
$$
A^b(\widetilde y,D)v=\int_{\{\vert \xi\vert \ge 1\}} a(\widetilde y,\xi)
F_{y\to \xi}v(\xi)e^{{i}<y,\xi>}d\xi
$$
is a continuous operator from $H^{ 1}_0({\mathcal O}\times {\Bbb R})$ into
$L^2({\mathcal O}\times {\Bbb R})$ with the norm satisfying
$$
\Vert A^b\Vert_{\mathcal L(H^{ 1}_0({\mathcal O}\times {\Bbb R}),
L^2({\mathcal O}\times {\Bbb R}))}\le \sum_{j=1}^\infty C_m\lambda_j^{-2n}.
$$
Next we observe that
\begin{eqnarray}\label{Fnona}
\Vert \widetilde A(\widetilde y,D)v\Vert_{L^2({\mathcal O}\times
{\Bbb R})}=\root\of{2\pi}\Vert A(\widetilde y,D_0,\dots,D_{n-1},\xi_n)u
F_{y_n\rightarrow \xi_n} w\Vert_{L^2({\mathcal O}\times
{\Bbb R})}\nonumber\\
 \le
C(\pi_{C^0(\mathcal O)}(a))\left(\int_{-\infty}^{\infty}\Vert u\Vert^2
_{H^{1,\xi_n}({\mathcal O})} \vert F_{y_n\rightarrow \xi_n}
w\vert^2d\xi_n\right)^\frac 12
\end{eqnarray}
for each function $v(y)=u(\widetilde y)w(y_n)$.

We take a sequence $\{w_j(x_n)\}_1^\infty$ such that
$F_{y_n\rightarrow \xi_n} w_j(\xi_n)$ has a compact support and
$\vert F_{y_n\rightarrow \xi_n} w_j\vert^2\rightarrow
\delta(\xi_n-s)$ for arbitrary $s\in{\Bbb R}$.
Since the function $\xi_n\rightarrow \Vert A(\widetilde y,\widetilde D,\xi_n)u
\Vert_{L^2({\mathcal O})}$ is continuous, we have
\begin{eqnarray}
&&\Vert A(\widetilde y,\widetilde D,\xi_n)u \widehat w_j\Vert^2
_{L^2({\mathcal O}\times
{\Bbb R})}\nonumber\\
= && \int_{{\Bbb R}}\Vert
A(\widetilde y,\widetilde D,\xi_n)u\Vert_{L^2({\mathcal O})}^2\vert
F_{y_n\rightarrow \xi_n} w_j\vert^2d\xi_n \rightarrow  \Vert
A(\widetilde y,\widetilde D,s)u\Vert_{L^2({\mathcal O})}^2.\nonumber
\end{eqnarray}
This fact and (\ref{Fnona}) imply
$$
\Vert A(\widetilde y,\widetilde D,s)u\Vert_{L^2({\mathcal O})}\le C\pi
_{C^0(\mathcal O)}(a) \Vert u\Vert_{H^{1,s}({\mathcal O})}
$$
for almost all $s$. Since the norm of the operator $A$ is a continuous
function of $s$, we have this inequality
for all $s$. $\blacksquare$

\noindent
The following theorem provides an estimate for a commutator of
a Lipschitz continuous function and the pseudodifferential operator
$\widetilde\omega_j.$

\begin{proposition}\label{Fzonaa}
Let $f\in W^1_\infty({\mathcal O})$ be a function with compact support.
Then
$$
\Vert [f,\widetilde \omega_j]\Vert_{\mathcal L (L^2({\mathcal
O}),H^{1,s}({\mathcal O}))}\le C\Vert
f\Vert_{W^1_\infty({\mathcal O})}\lambda^{4n}_j,
$$
where the constant $C$ is independent of $j$.
\end{proposition}
From this proposition we have immediately
\begin{corollary}\label{Fzonaa}
Let $f\in C^\ell({\mathcal O})$ be a function with compact support.
Then
$$
\Vert [f,\widetilde \omega_j]\Vert_{\mathcal L (L^2({\mathcal
O}),H^{\ell,s}({\mathcal O}))}\le C\Vert
f\Vert_{W^1_\infty({\mathcal O})}\lambda^{4n}_j,
$$
where the constant $C$ is independent of $j$.
\end{corollary}

\noindent
The proof of this theorem is similar to the proof of Corollary
in \cite{Stein}, p. 309.

Let $M(\xi)=\mu(\vert \xi\vert)\vert \xi\vert,$ where $\mu\in
C^\infty(\Bbb R)$, $\mu(t)=1 $ for $t\ge 1$ and $\mu(t)=0$ for
$t\in [0,\frac 12].$

\begin{lemma}\label{Fops1}
Let $a(\widetilde y,\widetilde\xi,s)\in C^\ell_{cl}S^{\ell,s}({\mathcal
O}).$ Then $A(\widetilde y,\widetilde D,s)^*
=A^*(\widetilde y,
\widetilde D,s)+R$, where $A^*$ is the
pseudodifferential operator with symbol
$\overline{a(\widetilde y,\widetilde\xi,s)}$ and $R\in {\mathcal
L}(H^{\ell-1, s}_0({\mathcal O}),L^2({\mathcal O}))$ satisfies
$$
\Vert R\Vert_{{\mathcal L}(H_0^{\ell-1,s}({\mathcal O}),L^2({\mathcal O}))}\le
C\pi_{C^\ell(\mathcal O)}(a).
$$
\end{lemma}
\noindent
{\bf Proof.} Thanks to the assumption {\bf C} it suffices to
consider the case when $$
a(\widetilde y,\tau\widetilde \xi, \tau s)=\tau^\ell
a(\widetilde y,\widetilde\xi,s)\quad \forall \tau>1.$$ The symbol
$a(\widetilde y,\xi)$ can be represented as
$$
a(\widetilde y,\xi)=\sum_{j=1}^\infty
a_j(\widetilde y)M^{\ell}(\xi) \omega_j(\xi).
$$
Consider the operator
$$
\widetilde A(\widetilde y,D)=\sum_{j=1}^\infty a_j(\widetilde y)M^{\ell}(D)
\widetilde \omega_j(D),\quad M^{\ell}(D)w=\int_{{\Bbb R}^{n+1}}M^{\ell}(\xi)
\widehat w e^{{i}<y,\xi>}d\xi.
$$
Then we find the formal adjoint operator:
\begin{eqnarray}
\widetilde A(\widetilde y,D)^*=\sum_{j=1}^\infty(a_j(\widetilde y)M^{\ell}(D)
\widetilde
\omega_j(D))^*=\sum_{j=1}^\infty M^{\ell}(D)\widetilde
\omega_j(D)\overline{a_j(\widetilde y)} \nonumber\\
=\sum_{j=1}^\infty\overline{a_j(\widetilde y)}M^{\ell}(D)\widetilde
\omega_j(D)+\sum_{j=1}^\infty[M^{\ell}(D)\widetilde
\omega_j(D),\overline{a_j(\widetilde y)}],             \nonumber
\end{eqnarray}
where $[A,B] := AB - BA$.

Observe that $\sum_{j=1}^\infty\overline{a_j(\widetilde y)}M^{\ell}(D)
\widetilde
\omega_j(D)$ is the operator with symbol
$\overline{a(\widetilde y,\xi_0,\dots,\xi_{n})}\in C^\ell_{cl}S^{\ell,s}
({\mathcal O}).$  Let us estimate the norm of the operator
$\sum_{j=1}^\infty[M^{\ell-1}(D)\widetilde \omega_j(D),
\overline{a_j(\widetilde y)}].$
Proposition \ref{Fzonaa} implies
\begin{eqnarray}\label{dik}
\left\Vert\sum_{j=1}^\infty[\overline{a_j(\widetilde y)},
M^{\ell}(D)\widetilde
\omega_j(D)]\right\Vert_{\mathcal L(H^{\ell-1}_0({\mathcal O}\times \Bbb R),
L^2(\mathcal O\times \Bbb R))}\nonumber\\
\le C_m\sum_{j=1}^\infty\Vert
a_j\Vert_{C^\ell(\bar{\mathcal O})}\lambda_j^{\widetilde \kappa(n)}\le
C_m\sum_{j=1}^\infty \lambda_j^{\widetilde
\kappa(n)}\pi_{C^\ell(\mathcal O)}(a)\lambda_j^{-m}
\le C\pi_{C^\ell(\mathcal O)}(a).
\end{eqnarray}
Denote $v=u(\widetilde y)w(y_n), \widetilde v=\widetilde
u(y_0,\dots,y_{n-1})\widetilde w(y_n).$ We have
$$
(\widetilde A(\widetilde y,D)v,\widetilde v)_{L^2({\mathcal O}\times {\Bbb R})}
=(v,\widetilde
A(\widetilde y,D)^*\widetilde v)_{L^2({\mathcal O}\times {\Bbb R})}
=(v,\widetilde
A^*(\widetilde y,D)\widetilde v)_{L^2({\mathcal O}\times {\Bbb R})}
+(v,R\widetilde
v)_{L^2({\mathcal O}\times {\Bbb R})}.
$$
By (\ref{dik}), we have
\begin{equation}
\Vert R\Vert_{{\mathcal L}(H^{\ell-1}_0({\mathcal O}\times \Bbb R),
L^2({\mathcal O}
\times \Bbb R))}\le
C\pi_{C^\ell(\mathcal O)}(a).
\end{equation}
On the other hand
\begin{align*}
& (\widetilde A(\widetilde y,D)v,\widetilde v)_{L^2({\mathcal O}\times
{\Bbb R})}=2\pi\int_{{\Bbb R}} (A(\widetilde y,\widetilde D,\xi_n)u,
\widetilde u)
_{L^2({\mathcal O})}w\overline{\widehat{\widetilde w}}d\xi_n\\
= &2\pi\int_{{\Bbb R}} (u,A(\widetilde y,\widetilde D,\xi_n)^*\widetilde
u)_{L^2(\mathcal O)}w\overline{\widetilde w}d\xi_n.
\end{align*}
Taking into account that $(v,A^*(\widetilde y,D)\widetilde v)_{L^2({\mathcal
O}\times {\Bbb R})}=\int_{{\Bbb R}} (u,A^*(\widetilde y,\widetilde D,\xi_n)
\widetilde u)_{L^2({\mathcal O})}w\overline{\widetilde w}d\xi_n$, we have
\begin{align*}
&\left\vert\int_{{\Bbb R}} (u,(A(\widetilde y,\widetilde D,\xi_n)^*-A^*
(\widetilde y,\widetilde D,\xi_n))\widetilde
u)_{L^2({\mathcal O})}w\overline{\widetilde w}d\xi_n \right\vert
= \vert(v,R\widetilde v)_{L^2({\mathcal O}\times {\Bbb R})}\vert\\
\le &C\pi_{C^\ell}(a)\Vert
v\Vert_{L^2({\mathcal O}\times {\Bbb R})}\Vert \widetilde
v\Vert_{H^{\ell-1}_0({\mathcal O}\times {\Bbb R})}.
\end{align*}
We take a sequence $\{w_j\}_{j=1}^\infty$, $j \in \Bbb N$ such that
$F_{y_n\rightarrow \xi_n}w_j$ have compact supports and $\vert
F_{y_n\rightarrow \xi_n} w_j\vert^2\rightarrow \delta(\xi_n-s)$
for arbitrary $s\in{\Bbb R}$.
Since the function
$\xi_n\rightarrow \Vert A(\widetilde y,\widetilde D,\xi_n)u\Vert
_{L^2({\mathcal O})}$
is continuous, we have
$$
\left\vert\int_{{\Bbb R}} (u,(A(\widetilde y,\widetilde D,\xi_n)^*-A^*
(\widetilde y,\widetilde D,\xi_n))\widetilde
u)_{L^2({\mathcal O})}\vert w_j\vert^2d\xi_n\right\vert\rightarrow
\vert(u,(A(\widetilde y,\widetilde D,s)^*-A^*(\widetilde y,\widetilde D,s))
\widetilde u)
_{L^2({\mathcal O})}\vert
$$
Since
$$
\vert(u,(A(\widetilde y,\widetilde D,s)^*-A^*(\widetilde y,\widetilde D,s))
\widetilde u)
_{L^2({\mathcal
O})}\vert\le C\pi_{C^1(\mathcal O)}(a)\Vert u\Vert_{L^2({\mathcal O})}\Vert
\widetilde u\Vert_{H^{\ell-1}({\mathcal O})},
$$
the lemma is proved.
$\blacksquare$

\begin{lemma}\label{Fops3}
Let $a(\widetilde y,\widetilde\xi,s)\in C^1_{cl}S^{\widehat
j,s}({\mathcal O})$ where $\widehat j = 0, 1$ and
$b(\widetilde y,\widetilde\xi,s)\in C^1_{cl}S^{\mu,s}({\mathcal
O}).$ Then $A(\widetilde y,\widetilde D,s)B(\widetilde y,\widetilde D,s)
=C(\widetilde y,\widetilde D,s)+R_0$ where $C(\widetilde y,\widetilde D,s)$
is the operator with symbol\newline
$a(\widetilde y,\widetilde\xi,s)b(\widetilde y,\widetilde\xi,s)$ and $R_0\in
{\mathcal L}(H_0^{\mu+{\widetilde s},s}({\mathcal
O}),H^{{\widetilde s}+1,s}({\mathcal O}))$ for any $\widetilde{s}\in
[-1,0]$ if $\widehat j=0$ and $R_0\in {\mathcal
L}(H_0^{\mu,s}({\mathcal O}),L^2({\mathcal O}))$ if
$\widehat j=1.$ Moreover we have
$$
\Vert R_{0}\Vert_{{\mathcal L}(H_0^{\mu,s}({\mathcal
O}),L^2({\mathcal O}))}\le C\pi_{C^1(\mathcal O)}(a)\pi_{C^1(\mathcal O)}(b)
\quad
\mbox{for}\,\,\widehat j=1,
$$
$$
\Vert R_{0}\Vert_{{\mathcal
L}(H_0^{\mu+{\widetilde s},s}({\mathcal
O}),H^{{\widetilde s}+1,s}({\mathcal O}))}\le
C\pi_{C^1(\mathcal O)}(a)\pi_{C^1(\mathcal O)}(b)\quad \mbox{for}\,\,\widehat j
=0.
$$
\end{lemma}
\noindent
{\bf Proof.} We set
$$
A(\widetilde y,D)=\sum_{j=1}^\infty a_j(\widetilde y)M^{\widehat j}(D)
\widetilde\omega_j(D),
\quad B(\widetilde y,D)=\sum_{j=1}^\infty b_j(\widetilde y) M^{\mu}(D)
\widetilde\omega_j(D).$$
Observe that
\begin{align*}
&A(\widetilde y,D)B(\widetilde y,D)\\
=& \sum_{m,k=1}^\infty a_m(\widetilde y)b_k(\widetilde y)M^{\widehat
j+\mu}(D)\widetilde \omega_m(D)\widetilde \omega_k(D)+\sum_{m,k=1}^\infty
a_m(\widetilde y)[M^{\widehat j}\widetilde\omega_m,b_k]M^\mu(D)
\widetilde\omega_k(D).
\end{align*}
Since $C(\widetilde y,D)=\sum_{m,k=1}^\infty a_m(\widetilde y)b_k(\widetilde y)
M^{\widehat
j+\mu}(D)\widetilde \omega_m(D)\widetilde \omega_k(D)$, and for $\widehat j=1$,
\begin{align*}
&\Vert R_{0}\Vert_{{\mathcal L}(H^{\mu}_0({\mathcal
O}),L^2({\mathcal O}))}
= \left\Vert\sum_{m,k=1}^\infty a_m(\widetilde y)[M\widetilde\omega_m,b_k]
M^{\mu}(D)\widetilde\omega_k(D)\right\Vert_{{\mathcal L}(H^{\mu}_0
({\mathcal O}),
L^2({\mathcal O}))}\\
\le &\sum_{m,k=1}^\infty \Vert a_m\Vert_{C^1(\bar{\mathcal O})}
\Vert [M\widetilde\omega_m,b_k]\Vert_{{\mathcal L}(L^2,L^2)} \Vert
\widetilde\omega_k(D)\Vert_{{\mathcal L}(L^2(\mathcal O),L^2(\mathcal O))} \\
\le &C_l\pi_{C^0(\mathcal O)}(a)\sum_{m,k=1}^\infty \lambda_m^{-l}\Vert
[M\widetilde\omega_m,b_k]\Vert_{{\mathcal L}(L^2(\mathcal O),L^2(\mathcal O))}
\lambda_k^{\widetilde \kappa(n)}.
\end{align*}
Applying Proposition \ref{Fzonaa}, we obtain
\begin{eqnarray}
\Vert R_{0}\Vert_{{\mathcal L}(H^{\mu}_0({\mathcal
O}),L^2({\mathcal O}))}
\le \sum_{m,k=1}^\infty\Vert a_m\Vert_{C^1(\bar{\mathcal O})}\Vert
[M\widetilde\omega_m,b_k]\Vert_{{\mathcal L}(L^2(\mathcal O),L^2(\mathcal O))}
\Vert\widetilde\omega_k(D)\Vert_{{\mathcal L}(L^2(\mathcal O),L^2(\mathcal O))}
\nonumber\\
\le C_l\pi_{C^1(\mathcal O)}(a)\sum_{m,k=1}^\infty \lambda_m^{-l} \Vert
b_k\Vert_{C^1(\bar{\mathcal
O})}\lambda_m^{\widetilde\kappa_1(n)}\lambda_k^{\widetilde
\kappa(n)}\\
\le C_{l,l_1}\pi_{C^1(\mathcal O)}(a)\pi_{C^1(\mathcal O)}(b)\sum_{m,k=1}
^\infty \lambda_m^{-l}
\lambda_k^{-l_1}\lambda_m^{\widetilde\kappa_1(n)}\lambda_k^{\widetilde
\kappa(n)}\nonumber\\\le C_{l,l_1}\pi_{C^1(\mathcal O)}(a)\pi_{C^1(\mathcal O)}
(b)\sum_{k=1}^\infty
\lambda_m^{-l}\lambda_m^{\widetilde\kappa_1(n)} \sum_{m=1}^\infty
\lambda_k^{-l_1} \lambda_k^{\widetilde \kappa(n)}<\infty.\nonumber
\end{eqnarray}
Let $v=v_j=u(\widetilde y)w_j(y_n).$  We take a sequence
$\{w_j\}_{j=1}^\infty$ such that $F_{y_n\rightarrow \xi_n}w_j$, $j\in \Bbb N$,
have compact supports and $\vert F_{y_n\rightarrow \xi_n}
w_j\vert^2\rightarrow \delta(\xi_n-s)$ for arbitrary $s\in{\Bbb R}$.

Then
\begin{eqnarray}\label {Fvera}\Vert
A(\widetilde y,D)B(\widetilde y,D)v_j-C(\widetilde y,D)v_j\Vert^2
_{L^2({\mathcal O}\times
{\Bbb R})}\\
= 2\pi\int_{{\Bbb R}}\Vert
(A(\widetilde y,\widetilde D,\xi_n)B(\widetilde y,\widetilde D,\xi_n)
-C(\widetilde y,\widetilde D,\xi_n))u\Vert^2_{L^2({\mathcal
O})}\vert F_{y_n\rightarrow \xi_n} w_j\vert^2d\xi_n \le C\Vert
v_j\Vert^2_{H_0^{\mu}({{\Bbb R}^{n+1}})}               \nonumber
\end{eqnarray}
for any $u\in H_0^{1+\mu}({\mathcal O})$.

Passing to the limit in (\ref{Fvera}) as $j\rightarrow +\infty$, we
obtain
\begin{equation}
\Vert (A(\widetilde y,\widetilde D,s)B(\widetilde y,\widetilde D,s)
-C(\widetilde y,\widetilde D,s))u\Vert^2_{L^2({\mathcal O})}\le
C\Vert u\Vert^2_{H_0^{\mu,s}({\mathcal O})}.                      \nonumber
\end{equation}
Let $\widehat j=0.$ Then
\begin{equation}\label{Fzopa}
\Vert R_{0}\Vert_{{\mathcal
L}(H^{\mu+{s}}_0({\mathcal
O}),H^{s}_0({\mathcal O}))}\le
\sum_{m,k=1}^\infty\Vert a_m
[\widetilde\omega_m,b_k]M^{-{s}}\Vert_{{\mathcal L}(L^2({\mathcal
O}),H^{{s}+1}_0({\mathcal O}))} \Vert
\widetilde\omega_k(D)\Vert_{{\mathcal L}(L^2({\mathcal O}),L^2({\mathcal
O}))}.
\end{equation}
In order to estimate $\Vert a_m [\widetilde\omega_m,b_k]M^{-{s}}\Vert
_{{\mathcal L}(L^2({\mathcal O}),H^{{s}+1}_0({\mathcal O}))}$,
we observe that
$M^{{s}+1}a_m[\widetilde \omega_m,b_k]M^{-{s}}=a_m M^{{s}+1}[\widetilde
\omega_m,b_k]M^{-{s}}+[M^{{s}+1},a_m][\widetilde
\omega_m,b_k]M^{-{s}}.$ For the second term in this equality, we
have
\begin{equation}\label{Fzopa1}
\Vert[\widetilde \omega_m,b_k]M^{-{s}}\Vert_{{\mathcal L}
(L^2(\mathcal O),L^2(\mathcal O))}\le
\Vert b_k\Vert_{C^1(\bar{\mathcal O})}\Vert \widetilde
\omega_m\Vert_{{\mathcal L}(L^2({\mathcal O}),L^2({\mathcal O}))},
\end{equation}
$$
\Vert [M^{{s}+1},a_m]\Vert_{{\mathcal L}(L^2(\mathcal O),L^2(\mathcal O))}\le
C\Vert a_m\Vert_{C^1(\bar{\mathcal O})}.
$$

In order to estimate the first term, we observe that $ [\widetilde
\omega_m,b_k]^*=-[\widetilde \omega_m,b_k].$ Then
$$ [\widetilde \omega_m,b_k]\in{\mathcal L} (L^2({\Bbb R}^n),
H^{1}({\Bbb R}^n)),\quad [\widetilde \omega_m,b_k]\in{\mathcal L} (H^{-1}
({\Bbb R}^n), L^2({\Bbb R}^n))
$$

An interpolation argument yields
\begin{equation}\label{Fzopa2}
[\widetilde \omega_m,b_k]\in{\mathcal L} (H^{-\gamma}({\Bbb R}^n), H^{1-\gamma}
({\Bbb R}^n))
\quad \forall \gamma\in [0,1]
\end{equation}
and
\begin{equation}\label{Fzopa3}
\Vert[\widetilde \omega_m,b_k]\Vert_{{\mathcal L} (H^{-\gamma}({\Bbb R}^n),
H^{1-\gamma}({\Bbb R}^n))}\le \Vert[\widetilde
\omega_m,b_k]\Vert_{{\mathcal L} (L^2({\Bbb R}^n),
H^{1}({\Bbb R}^n))}  \quad \forall \gamma\in [0,1].
\end{equation}
Applying (\ref{Fzopa1})-(\ref{Fzopa3}) to (\ref{Fzopa}) we obtain
\begin{eqnarray}\label{Fzopa4}
&&\Vert R_{0}\Vert_{{\mathcal
L}(H^{\mu+{s}}_0({\mathcal
O}),H^{{s}+1}_0({\mathcal O}))}
\le C_l\sum_{m,k=1}^\infty \lambda_m^{-l} \Vert
b_k\Vert_{C^1(\bar{\mathcal
O})}\lambda_m^{\widetilde\kappa_1(n)}\lambda_k^{\widetilde
\kappa(n)}\nonumber\\
&&\le C_{l,l_1}\pi_{C^1(\mathcal O)}(b)\sum_{m,k=1}^\infty \lambda_m^{-l}
\lambda_k^{-l_1}\lambda_m^{\widetilde\kappa_1(n)}\lambda_k^{\widetilde
\kappa(n)}\le C_{l,l_1}\pi_{C^1(\mathcal O)}(b)\sum_{k=1}^\infty
\lambda_m^{-l}\lambda_m^{\widetilde\kappa_1(n)} \sum_{m=1}^\infty
\lambda_k^{-l_1} \lambda_k^{\widetilde \kappa(n)}<\infty.\nonumber
\end{eqnarray}
We finish the proof of lemma using similar arguments as in case
$\widehat j=1$.$\blacksquare$

\noindent
The direct consequence of Lemma \ref{Fops3} is the following
commutator estimate.

\begin{lemma}\label{Fops2}
Let $a(\widetilde y,\widetilde\xi,s)\in C^1_{cl}S^{1,s}({\mathcal
O})$ and $b(\widetilde y,\widetilde \xi,s)\in C^1_{cl}S^{1,s}({\mathcal
O}).$ \\
Then the commutator $[A,B]$ belongs to the space $ {\mathcal L}(H^{1,s}(
{\mathcal O});L^2({\mathcal
O}))$ and
$$\Vert [A,B]\Vert_{{\mathcal L}(H^{1,s}({\mathcal O});L^2({\mathcal O}))}\le
C(\pi_{C^0(\mathcal O)}(a)\pi_{C^0(\mathcal O)}(b)+\pi_{C^0(\mathcal O)}(a)
\pi_{C^1(\mathcal O)}(b)+\pi_{C^1(\mathcal O)}(a)\pi_{C^0(\mathcal O)}(b)).
$$
\end{lemma}
\noindent
{\bf Proof.} By Lemma \ref{Fops3} we have
$$
A(\widetilde y,\widetilde D,s)B(\widetilde y,\widetilde D,s)=C(\widetilde y,
\widetilde D,s)+R_0,\quad
B(\widetilde y,\widetilde D,s)A(\widetilde y,\widetilde D,s)=C(\widetilde y,
\widetilde D,s)+\widetilde R_0,
$$
where $R_0,\widetilde R_0\in \mathcal{L}(H^{1,s}({\mathcal O}),L^2({\mathcal
O}))$ and $C(\widetilde y,\widetilde D,s)$ is the pseudodifferential operator
with the symbol $c(\widetilde y,\widetilde \xi,s)=a(\widetilde y,
\widetilde \xi,s)b(\widetilde y,\widetilde \xi,s).$ Since $[A,B]=R_0
-\widetilde R_0$, we immediately obtain the
statement of the lemma. $\blacksquare$

\begin{lemma}\label{Fops4}
Let $a(\tilde y,\tilde \xi,s)\in C^1_{cl}S^{1,s}({\mathcal O})$ be a
symbol with compact support in ${\mathcal O}.$ Let  $\mathcal O_i
\subset\subset\mathcal O$ and $\mathcal O_1\cap\mathcal O_2=\emptyset .$
Suppose that $u\in H^{1,s}({\mathcal O})$ and $\mbox{supp}\,u
\subset {\mathcal O_1}.$  Then there exists a constant
$C$ such that
\begin{equation}\label{gabon}
\Vert A(\tilde y,\tilde D,s)u\Vert_{H^{1,s}({\mathcal O_2})}
\le \frac{C \pi_{C^1(\mathcal O)}(a)}{\mbox{dist} (\mathcal O_2,\mathcal O_1)^{2n+3}}
\Vert u\Vert_{H^{1,s}({\mathcal O})}.
\end{equation}
\end{lemma}
\noindent {\bf Proof.} By lemma  it suffices to prove  the inequality
(\ref{gabon}) only for $u\in C^\infty_0(\mathcal O_1).$ Let $b(t)\in C^\infty_
0(-2,2)$, and $b\vert_{[-1,1]}=1$  and  $\tilde y\in\mathcal O_2.$ We have
\begin{eqnarray}
A(\tilde y, \tilde D,s) u=\lim_{\epsilon\rightarrow +0}\int_{\Bbb R^n}
b(\epsilon \vert \tilde\xi\vert)a(\tilde y,\tilde\xi,s)e^{-i<\tilde y,\tilde
\xi>}\hat u(\tilde \xi)d\tilde\xi=\nonumber\\
\frac{1}{(2\pi)^\frac n2}\lim_{\epsilon\rightarrow +0}\int_{\Bbb R^n}
\int_{\Bbb R^n}b(\epsilon \vert \tilde\xi\vert)a(\tilde y,\tilde\xi,s)
e^{i<\tilde x-\tilde y, \xi>}u(\tilde x)d\tilde\xi d\tilde x=\nonumber\\
\frac{1}{i^{2k}(2\pi)^\frac n2}\lim_{\epsilon\rightarrow +0}\int_{\Bbb R^n}
\int_{\Bbb R^n}\frac{b(\epsilon \vert \tilde\xi\vert)a(\tilde y,\tilde\xi,s)}
{\vert\tilde x-\tilde y\vert^{2k}}(\Delta^k_\xi
e^{i<\tilde x-\tilde y,\tilde \xi>})u(\tilde x)d\tilde\xi d\tilde x=
\nonumber\\
\frac{1}{i^{2k}(2\pi)^\frac n2}\lim_{\epsilon\rightarrow +0}\int_{\Bbb R^n}
\int_{\Bbb R^n}\frac{\Delta^k_{\tilde\xi} (b(\epsilon \vert \tilde\xi\vert)
a(\tilde y,\tilde\xi,s)
)}{\vert \tilde x-\tilde y\vert^{2k}}e^{i<\tilde x-\tilde y, \tilde\xi>}
u(\tilde x)d\tilde\xi d\tilde x=\nonumber\\
\frac{1}{i^{2k}(2\pi)^\frac n2}\int_{\mathcal O_1}\int_{\Bbb R^n}
\frac{\Delta^k_{\tilde\xi} a(\tilde y,\tilde\xi,s)
}{\vert\tilde x-\tilde y\vert^{2k}}e^{i<\tilde x-\tilde y, \tilde\xi>}
u(\tilde x)d\tilde \xi d\tilde x.
\end{eqnarray}
Let $k=n.$ Then
\begin{eqnarray}
\vert \partial_{y_j}A(\tilde y, \tilde D,s) u\vert
= \left \vert\frac{1}{i^{2k}(2\pi)^\frac n2}\int_{\mathcal O_1}\int_{\Bbb R^n}
\partial_{y_i}\left (\frac{\Delta^k_{\tilde\xi} a(\tilde y,\tilde\xi,s)
}{\vert \tilde x-\tilde y\vert^{2k}}e^{i<\tilde x-\tilde y, \tilde\xi>}
\right)u(\tilde x)d\tilde\xi d\tilde x\right\vert\\
\le C\pi_{C^1(\mathcal O)}(a)\int_{\mathcal O_1}\int_{\Bbb R^n}\frac{\vert u(\tilde x)
\vert}{\vert \tilde x-\tilde y\vert^{2n+3}(1+\vert\tilde \xi\vert^{2n+1}
+s^{2n+1})}d\tilde\xi d\tilde x\nonumber\\
\le \frac{C \pi_{C^1(\mathcal O)}(a)}
{(s^{n+1}+1)\mbox{dist} (\mathcal O_2,\mathcal O_1)^{2n+3}}
\Vert u\Vert_{H^{1,s}({\mathcal O})}.\nonumber
\end{eqnarray}
Proof of the lemma is complete.
 $\blacksquare$

\noindent We shall use the following variant of the G\aa rding
inequality:
\begin{lemma}\label{Fops5}
Let $p(\tilde y,\tilde\xi,s)\in C^2_{cl}S^{2,s}({\mathcal O})$ be a
symbol with compact support in ${\mathcal O}.$
Let $u\in H^{1,s}({\mathcal O})$ and $\mbox{supp}\,u\subset\mathcal O_1.$ Let
$\mathcal O_1\subset\subset\mathcal O_2\subset\subset\mathcal O_3
\subset\subset \mathcal O$  and $\tilde\gamma\in C^\infty_0(\mathcal O_3)$
be a function
such that $\tilde\gamma\vert_{\mathcal O_2}=1$  be such that
$\mbox{Re}\,p(\tilde y,\tilde\xi,s)> \hat C\vert(\tilde \xi,s)\vert^2$ for any
$\tilde y\in\mathcal O_3$.  Then
\begin{eqnarray}\label{goblin}
\mbox{Re}(P(\tilde y,\tilde D,s)u,u)_{L^2(\mathcal O)}\ge \frac{\hat C}{2}
\Vert u\Vert^2_{H^{1,s}({\mathcal O})}\\-C_1\left (\left(
\sum_{k=0}^2(\pi_{C^k(\mathcal O_3)}(p)+1)\pi_{C^{2-k}(\mathcal O_3)}
(\tilde \gamma)+(\pi^2_{C^k(\mathcal O_3)}(p)+1)\pi^2_{C^{2-k}(\mathcal O_3)}
(\tilde \gamma))\right)^2\right.\nonumber\\
\left.+ \frac{1}{dist (\mathcal O_1,\Bbb R^n\setminus\mathcal O_2)^{2n+3}}
\right)\Vert u\Vert^2_{L^2({\mathcal O})}.\nonumber
\end{eqnarray}
\end{lemma}
\noindent {\bf Proof.} Consider the
pseudodifferential operator $A(\tilde y,\tilde D,s)$ with symbol
$A(\tilde y,\tilde \xi,s)=$\newline$(\tilde\gamma\mbox{Re}\, p(\tilde y,\tilde \xi,s)
-\tilde\gamma\frac{\hat C}{2} M^2(\tilde \xi,s))^\frac 12\in C^2_{cl}S^{1,s}
({\mathcal O}).$ Then, according to Lemma \ref{Fops3}
$$
A(\tilde y,\tilde D,s)^*A(\tilde y,\tilde D,s)
=\tilde \gamma\mbox{Re}\, p(\tilde y,\tilde D,s)-\tilde\gamma\frac{\hat
C}{2}M^2(\tilde D,s)+R,
$$
where $R\in {\mathcal L}(H^{1,s}({\mathcal O});L^2({\mathcal O}))$ and
\begin{eqnarray}
\Vert R\Vert_{{\mathcal L}(H^{1,s}({\mathcal O});L^2({\mathcal O}))}
\le C(\pi_{C^2(\mathcal O_3)}(a)+\pi^2_{C^2(\mathcal O_3)}(a))\nonumber\\
\le C\sum_{k=0}^2(\pi_{C^k(\mathcal O_3)}(p)+1)\pi_{C^{2-k}(\mathcal O_3)}
(\tilde \gamma)+(\pi^2_{C^k(\mathcal O_3)}(p)+1)\pi^2_{C^{2-k}(\mathcal O_3)}
(\tilde \gamma)).
\end{eqnarray}

Therefore
\begin{eqnarray}
\mbox{Re}(P(\tilde y,\tilde D,s)u,u)_{L^2({\mathcal O})}
=\Vert A(\tilde y,\tilde D,s)\Vert^2_{L^2({\mathcal O})}-((1-\tilde\gamma)
M^2(\tilde D,s)u,u)_{L^2({\Bbb R}^n)}\nonumber\\
+\frac{\hat C}{2}\Vert u\Vert^2_{H^{1,s}({\mathcal
O})}+(Ru,u)_{L^2({\mathcal O})}.\nonumber
\end{eqnarray}
Observing that \begin{eqnarray}\vert(Ru,u)_{L^2({\mathcal O})}\vert \le C\sum_{k=0}^2(\pi_{C^k(\mathcal O_3)}(p)+1)\pi_{C^{2-k}(\mathcal O_3)}
(\tilde \gamma)\\
C\sum_{k=0}^2(\pi_{C^k(\mathcal O_3)}(p)+1)\pi_{C^{2-k}(\mathcal O_3)}
(\tilde \gamma)+(\pi^2_{C^k(\mathcal O_3)}(p)+1)\pi^2_{C^{2-k}(\mathcal O_3)}
(\tilde \gamma))\Vert u\Vert_{L^2({\mathcal O})}
\Vert u\Vert_{H^{1,s}({\mathcal O})},\nonumber
\end{eqnarray}
and since by Lemma \ref{Fops4}
 we have
$$
\vert((1-\tilde\gamma)M^2(\tilde D,s)u,u)
_{L^2({\Bbb R}^n)}\vert \le \frac{C}{dist (\mathcal O_1,\Bbb R^n\setminus
\mathcal O_2)^{2n+3}}\Vert
u\Vert^2_{L^2({\mathcal O})},
$$
we obtain the statement of the lemma.
$\Box$

\end{document}